\documentclass[a4paper,leqno]{amsart}

\usepackage{amsmath}
\usepackage{amsthm}
\usepackage{amssymb}
\usepackage{amsfonts}
\usepackage[applemac]{inputenc}
\usepackage{mathrsfs}
\usepackage{esint}
\usepackage{color}
\usepackage{dsfont}
\usepackage[body={15cm,21.5cm},centering]{geometry} 

\newcommand{\overbar}[1]{\mkern 1.5mu\overline{\mkern-1.5mu#1\mkern-1.5mu}\mkern 1.5mu}

\def\le{\leqslant}
\def\ge{\geqslant}

\DeclareMathOperator*{\oscDaniel}{osc}

\newcommand \dps{\displaystyle }

\newcommand{\loc}{\mathrm{loc}}

\newcommand{\R}{\mathbb{R}}
\newcommand{\Z}{\mathbb{Z}}
\newcommand{\N}{\mathbb{N}}

\newcommand{\Id}{\text{Id}}
\newcommand{\e}{\varepsilon}

\newcommand{\calF}{\mathcal{F}}

\newcommand{\calY}{\mathcal{Y}}

\mathchardef\emptyset="001F

\newcommand{\ext}[1]{\overline{#1}}

\newcommand{\dig}[1]{\mathrm{diag}\left[ #1\right]}

\newcommand{\Expec}[1]{\left\langle #1 \right\rangle}
\newcommand{\expec}[1]{\left\langle #1 \right\rangle}
\newcommand{\step}[1]{\noindent \textit{Step} #1.}

\newcommand{\osc}[2]{\underset{\dps #1}{\mathrm{osc}} \,#2\,} 
\newcommand{\Sup}[2]{\underset{\dps #1}{\mathrm{sup}} \,#2\,}

\theoremstyle{plain}
\newtheorem{theorem}{Theorem}[section]
\newtheorem{lemma}[theorem]{Lemma}
\newtheorem{corollary}[theorem]{Corollary}
\newtheorem{proposition}[theorem]{Proposition}

\theoremstyle{definition}
\newtheorem{definition}[theorem]{Definition}

\newtheorem{remark}[theorem]{Remark}
\newtheorem*{remark*}{Remark}

\numberwithin{equation}{section}

\begin{document}

\title[Annealed Green functions and uncertainty quantification]
{Annealed estimates on the Green functions and uncertainty quantification}
\author[A.\ Gloria]{Antoine Gloria}
\address[A.\ Gloria]{Universit\'e Libre de Bruxelles (ULB) \\ Brussels, Belgium \\ and Team MEPHYSTO \\  Inria Lille - Nord Europe \\ Villeneuve d'Ascq, France}
\email{antoine.gloria@ulb.ac.be}
\author[D.\ Marahrens]{Daniel Marahrens}
\address[D.\ Marahrens]{Max-Planck-Institute for Mathematics in the Sciences\\
Inselstrasse 22, 04103 Leipzig\\ Germany}
\email{daniel.marahrens@mis.mpg.de}

\begin{abstract}
We prove Lipschitz bounds for linear elliptic equations in divergence form whose measurable coefficients are random stationary and satisfy a logarithmic
Sobolev inequality, extending to the continuum setting results by Otto and the second author for discrete elliptic equations. This improves 
the celebrated De Giorgi-Nash-Moser theory
 in the large (that is, away from the singularity) for this class of coefficients.
This regularity result is obtained as a corollary of optimal decay estimates on the derivative and mixed second derivative of the elliptic Green functions on $\R^d$.
As another application of these decay estimates we derive optimal estimates on the fluctuations of solutions of linear elliptic PDEs with ``noisy" diffusion coefficients.
\end{abstract}

\date{\today}


\maketitle
\setcounter{tocdepth}{1}
\tableofcontents

\section{Introduction}

For scalar linear elliptic equations in divergence form it is well-known that the best regularity theory one can hope for is that of De Giorgi, Nash, and Moser. In particular, solutions are H\"older continuous for some exponent $1\ge \alpha>0$ that depends only on the ellipticity contrast of the coefficient field ($\alpha=1$ for constant coefficients), see \cite{Han-Lin-97}. In view of explicit examples from quasiconformal mappings, see \cite[Theorem~12.3]{Gilbarg-Trudinger-98}, $\alpha<1$ for non-constant coefficients in general. 

\medskip

In the case when the coefficient field is periodic (and H\"older-continuous), Avellaneda and Lin proved in \cite{Avellaneda-Lin-87,Avellaneda-Lin-91} that $\alpha=1$ as well. (Indeed, the known counterexamples to optimal regularity cannot be periodic.) Their proof is based on a Campanato iteration (and the availability of periodic correctors) to lift the regularity of the associated homogenized equation to the non constant coefficients equation at large scales (whereas the small-scale behavior is controlled by the H\"older-regularity assumption on the coefficients via the Schauder theory). This also allows them to prove that the associated Green function has essentially the same behavior as for the Laplace equation.

\medskip

To extend the results by Avellaneda and Lin to the random setting, we face a ``lack of compactness" (it is no longer possible to rely on correctors, which are not necessarily well-behaved a priori). 
In their first contribution \cite{Gloria-Otto-09} to quantitative stochastic homogenization, Otto and the first author proved that the corrector gradient has bounded finite moments --- a Lipschitz-type regularity --- under 
a quantitative ergodicity assumption on the coefficients. 
These are the first ``improved regularity" results for an elliptic equation with random coefficients.
The interpretation of these results in terms of ``improved regularity" and their extension to more general equations than the corrector equation first appeared in the work \cite{Marahrens-Otto-13} by Otto and the second author for discrete elliptic equations. In this work, the authors proceed in a different way than Avellaneda and Lin, and start with the optimal control of the finite moments of the Green functions at large scales. In turn this allows them to improve the H\"older regularity exponent $\alpha$ for this class of coefficients. 
Besides the structure of their proof, the Green functions bounds they obtain are particularly relevant to stochastic homogenization.
Indeed, a key ingredient to \cite{Gloria-Otto-09,Gloria-Otto-10b,Gloria-Neukamm-Otto-15} is a so-called sensitivity estimate, which naturally involves Green's functions (see for instance Lemma~\ref{lem:osc_u} below). Their optimal control leads to the optimal control of several quantities of interest, like the error in the two-scale expansion (see \cite{Gloria-Neukamm-Otto-14}) or the fluctuations in elliptic equations with noisy coefficients (see \cite{Gloria-12c,Marahrens-Otto-13}).

\medskip

The aim of the present article is to extend the results by Otto and the second author in \cite{Marahrens-Otto-13} to the continuum setting of linear (non-necessarily self-adjoint) elliptic PDEs.
First, we develop a Lipschitz regularity theory for linear elliptic equations whose coefficients satisfy a quantitative ergodicity assumption in the form of a logarithmic-Sobolev inequality, see Definition~\ref{def:LSI}
and Theorem~\ref{thm:main2}.
Second, we obtain optimal bounds on the gradient and second-mixed gradient of the associated Green function, see Theorem~\ref{thm:main}.
Last we improve the fluctuation estimates of both \cite{Gloria-12c} and \cite{Marahrens-Otto-13}, and we unravel the central limit theorem scaling of a weak measure of the fluctuations, see Theorems~\ref{thm:strong_error} and \ref{thm:weak_error}.

\medskip 

We conclude this introduction by mentioning the independent and inspiring work by Armstrong and Smart.
In \cite{Armstrong-Smart-14}, the authors obtain a similar Lipschitz regularity theory,
with however better moment bounds and for nonlinear equations, under the assumption that the coefficients have finite range of dependence.
Their approach is much closer to the approach by Avellaneda and Lin, and rely on a Campanato iteration
using a quantitative homogenization result (to replace the compactness argument).

\section{Statement of the main results}


\subsection{Notation and assumptions on the coefficient field}

We let  $\lambda\in(0,1]$ denote an ellipticity constant
which is fixed  throughout the paper, and set
\begin{align}
  \Omega_0:=\Big\{\,A_0\in\R^{d\times d}\,:\,&A_0\text{ is bounded, i.~e. }|A_0\xi| \leq |\xi|\text{ for all $\xi\in\R^d$,}\nonumber \\
  &A_0\text{ is  elliptic, i.~e. } \lambda |\xi|^2\leq \xi\cdot A_0\xi \text{ for all $\xi\in\R^d$}\,\Big\}.\label{def:alpha-beta} 
\end{align}
We equip $\Omega_0$ with the usual
topology of $\R^{d\times d}$. A \textit{coefficient field},
denoted by $A$, is a Lebesgue-measurable function on $\R^d$ taking values in $\Omega_0$. 
We then define 
\begin{equation*}
\Omega:=\{\text{measurable maps }A:\R^d\to \Omega_0\},
\end{equation*}
which we equip with the $\sigma$-algebra $\calF$ 
that makes the evaluations $A\mapsto \int_{\R^d} A_{ij}(x)\chi(x)dx$ measurable for all $i,j\in \{1,\dots,d\}$ and 
all smooth functions $\chi$ with compact support.  This makes  $\calF$ countably generated.

\medskip

Following the convention in
statistical mechanics, we describe a \textit{random coefficient field} by equipping $(\Omega,\calF)$
with an ensemble $\expec{\cdot}$ (the expected value). Following \cite{Papanicolaou-Varadhan-79}, we shall assume that $\expec{\cdot}$ is stochastically continuous: For all $\delta>0$ and $x\in \R^d$,
$$
\lim_{|h|\downarrow 0}  \expec{\mathds{1}_{\{A\,:\,|A(x+h)-A(x)|>\delta\}}}\,=\,0
$$
We shall always assume that $\expec{\cdot}$ is
\textit{stationary}, i.~e. for all translations $z\in\R^d$ the
coefficient fields $\{\R^d\ni x\mapsto  A(x)\}$ and $\{\R^d\ni
x\mapsto  A(x+z)\}$ have the same joint distribution under $\expec{\cdot}$. 
Let
$\tau_z:\Omega\to\Omega,\; A(\cdot)\mapsto A(\cdot+z)$ denote the shift
by $z$, then $\expec{\cdot}$ is stationary if and only if $\tau_z$ is
$\expec{\cdot}$-preserving for all shifts $z\in\R^d$.
The stochastic continuity assumption ensures that the map $ \R^d \times \Omega \to \Omega, (x,\omega)\mapsto \tau_x \omega$
is measurable (where $\R^d$ is equipped with the $\sigma$-algebra of Lebesgue measurable sets).

\medskip

A random variable is a measurable function on $(\Omega,\calF)$.
A \textit{random field} $\tilde\zeta$ is a measurable function on $\R^d\times \Omega$.
In this article the random field under study is the Green function.
We are interested in the behaviour of the (massive) Green function
$G_\mu:\R^d\times\R^d\times \Omega \to \R$, which is defined for all $\mu>0$ 
and for all $y\in \R^d$ as the unique 
distributional solution in $W^{1,1}(\R^d)$ which is continuous away from the diagonal $x=y$ of the elliptic equation
\begin{equation}\label{Green}
 \mu G_\mu(x,y;A)-\nabla_x \cdot (A(x) \nabla_x G_\mu(x,y;A)) = \delta(x-y).
\end{equation}
For the existence, uniqueness and properties of $G_\mu$, see Definition~\ref{def:Green}.
Note that by definition of the $\sigma$-algebra, $G_\mu$ is measurable.

We make a quantitative ergodicity assumption in the form of the following  logarithmic Sobolev inequality.

%
\begin{definition}[Logarithmic Sobolev inequality (LSI)]\label{def:LSI}
 We say that the ensemble $\langle \cdot \rangle$ satisfies a logarithmic Sobolev inequality if there exist constants $\rho, \ell > 0$,
which we shall respectively call amplitude and correlation-length, such that 
\begin{equation}\label{LSI}
 \Expec{ \zeta^2 \log \frac{\zeta^2}{\langle \zeta^2 \rangle} }\, \le \, \frac{2}{\rho} \Expec{ \int_{\R^d}
\Big( \osc{A|_{B_\ell(z)}}{\zeta} \Big)^2 \;dz }
\end{equation}
for all measurable functions $\zeta : \Omega \to \R$, where the expectation in the RHS is an outer expectation (the oscillation is not necessarily measurable). Here the expression $\osc{A|_{B_\ell(z)}}{\zeta}$ denotes the oscillation of $\zeta$ with respect to all coefficient fields that coincide with $A$ outside of $B_\ell(z)$,
where $B_\ell(z)$ is the ball of radius $\ell$ centered at $z\in\R^d$, that is,
\begin{eqnarray}
\left(\osc{A|_{B_\ell(z)}}{\zeta}\right)(A)
&=&\left(\sup_{\dps A|_{B_\ell(z)}}\zeta\right)(A)
-\left(\inf_{\dps A|_{B_\ell(z)}}\zeta\right)(A)\nonumber\\
&=&\sup\left\{\zeta(\tilde A)|\tilde A\in \Omega,\;
\tilde A|_{\mathbb{R}^d\setminus B_\ell(z)}=A|_{\mathbb{R}^d\setminus B_\ell(z)}\right\}\nonumber\\
&&-
\inf\left\{\zeta(\tilde A)|\tilde A\in\Omega,\;
\tilde A|_{\mathbb{R}^d\setminus B_\ell(z)}=A|_{\mathbb{R}^d\setminus B_\ell(z)}\right\}.\label{Lc.4}
\end{eqnarray}
\qed
\end{definition}
An example of coefficient field which satisfies (LSI) is the Poisson inclusions process (and variants of it), see in particular~\cite{Bella-Otto-13}.
Without loss of generality, we assume in this article that $\ell\ge 1$.
\begin{remark}
The fact that outer expectations appear in the RHS of \eqref{LSI} is not a difficulty since in the rest of the article we shall always estimate the RHS of \eqref{LSI} by the expectation of measurable quantities (for which outer expectation and expectation coincide).
\qed
\end{remark}

\subsection{Lipschitz-regularity theory}

One way to formulate the De Giorgi-Nash-Moser theory is as follows: There exists $0<\alpha \le 1$ 
depending only on the ellipticity ratio $\lambda$ such that for all $p>\frac{d}{2}$, $\kappa>0$, $R>0$, and $\mu\ge 0$ with $R^2\mu \le \kappa$, 
if $u$ satisfies
\begin{equation*}
\mu u-\nabla \cdot A \nabla u\,=\,f \mbox{ in }B_{2R},
\end{equation*}
for some $f\in L^p(B_{2R})$, then
\begin{equation}\label{DGNM-pde}
 R^\alpha \Sup{x,y\in B_R}{\frac{|u(x)-u(y)|}{|x-y|^\alpha} }\, \lesssim \, \Big(\fint_{B_{2R}} u^2\Big)^{\frac{1}{2}} + \Big(\fint_{B_{2R}} |R^2 f|^p \Big)^{\frac{1}{p}} ,
\end{equation}
see for instance \cite[Theorem~8.24]{Gilbarg-Trudinger-98}. (Note that this follows from the statement for $R=1$ since by \eqref{eq:elliptic-u}, $f$ is replaced by $R^2 f$ when performing a change of variables $x\leadsto R^{-1}x$.) In the supremum above, we have set by convention $\frac{0}{0} := 0$.
This result has two aspects: a regularity \emph{in the small} and a regularity \emph{in the large}. 
In particular we may split the statement into two parts: \emph{in the small}, that is for $|x|\lesssim 1$, \eqref{DGNM-pde} quantifies the high frequencies of $u$ (local regularity),
\begin{equation}\label{DGNM-pde-small}
\Sup{B_1}{\frac{|u(x)-u(0)|}{|x|^\alpha} }\, \lesssim \, \Big(\fint_{B_2} u^2\Big)^{\frac{1}{2}} + \Big(\fint_{B_2} |f|^p \Big)^{\frac{1}{p}} ,
\end{equation}
and \emph{in the large}, \eqref{DGNM-pde} quantifies the low frequencies of $u$ (growth at large scales),
\begin{equation}\label{DGNM-pde-large}
\Sup{B_R\setminus B_1}{\frac{|u(x)-u(0)|}{|x|^\alpha} }\, \lesssim \,  R^{-\alpha} \Big(\fint_{B_{2R}} u^2\Big)^{\frac{1}{2}} + R^{-\alpha}\Big(\fint_{B_{2R}} |R^2f|^p \Big)^{\frac{1}{p}} .
\end{equation}
If we assume that the coefficients $A$ are uniformly H\"older-continuous, then we have an optimal regularity theory in the small, that is, \eqref{DGNM-pde-small} holds for the improved exponent $\alpha=1$ provided $p>d$ (see for instance \cite[Theorem~3.13]{Han-Lin-97}). However, the De Giorgi-Nash-Moser exponent \emph{cannot} be improved in the large by increasing the regularity of the coefficients, as classical examples from quasiconformal mappings show.
The improvement of the De Giorgi-Nash-Moser exponent \emph{in the large} 
is the aim of the following result for stationary coefficients that satisfy (LSI) and for periodic coefficients.
\begin{theorem}\label{thm:main2}
Let the ensemble be stationary and satisfy (LSI) with constants $\rho$ and $\ell$,
and let $\mu\ge 0$ and $d< p <\infty$.
Then for all $R\ge 2\ell$ and all $x\in B_R\setminus B_{2\ell}$, there exists a random variable $\calY_R(x)$ 
with bounded finite moments such that for all $u$ and $f\in L^p(B_{2R})$ related via
\begin{equation}\label{eq:elliptic-u}
\mu u-\nabla \cdot A \nabla u\,=\,f \mbox{ in }B_{2R},
\end{equation}
we have
\begin{equation}\label{eq:thm-de_giorgi1}
\fint_{B_\ell}\frac{|u(x+x')-u(x')|}{|x|}dx'
\,\leq \,\calY_R(x) \Big( R^{-1}\Big(\fint_{B_{2R}} u^2 \Big)^{\frac{1}{2}}  +R^{-1} \Big(\fint_{B_{2R}} |R^{2}f|^p \Big)^{\frac{1}{p} }\Big).
\end{equation}
In addition the random variables $\calY_R$ have the following boundedness property: For all $1\le q<\infty$, there exists $C_q<\infty$ depending only on $d,\lambda,p,q,\rho,\ell$
such that
\begin{equation}\label{eq:Y}
\sup_{R\ge 2\ell} \sup_{x\in B_R\setminus B_{2\ell}} \expec{\calY_R(x)^q}^\frac{1}{q}\,\leq\, C_q.
\end{equation}
\qed
\end{theorem}
\begin{remark}\label{cor:de_giorgi}
In the case of uniformly H\"older continuous coefficients in the sense 
that there exists a constant $\mathcal{C}_\gamma < \infty$ such that $\langle \cdot \rangle$-almost surely $[A]_{C^\gamma} \le \mathcal{C}_\gamma$, the regularity theory of Theorem~\ref{thm:main2} also holds 
in the small, as it should. 
In particular, \eqref{eq:thm-de_giorgi1} holds true for all $x\in B_{R}$ and \eqref{eq:Y} is replaced by 
\begin{equation*}
\sup_{R\ge 2\ell} \sup_{x\in B_R} \expec{\calY_R(x)^q}^\frac{1}{q}\,\leq\, C_q.
\end{equation*}
\end{remark}

\subsection{Bounds on the Green functions}

In general, the only optimal decay result which holds without further smoothness assumption is the following consequence of the celebrated De Giorgi-Nash-Moser theory (in dimensions $d>2$)
on the Green function itself: For all $A\in \Omega$, and all $\mu\ge 0$,
$$
 0  \,\le \, G_\mu (x,y;A) \,\le\, C \frac{e^{-c\sqrt{\mu}|x-y|}}{|x-y|^{d-2}} 
$$
for some constants $c,C>0$ depending only on $\lambda$ and $d$, see Definition~\ref{def:Green} below.
For the constant-coefficient operator, i.e.\ the massive Laplacian, we also have the following optimal gradient estimate: For all $\mu\ge 0$,
\begin{equation}\label{eq:grad-Green-lap}
|\nabla G_\mu (x,y;\Id)| \,\le\, C \frac{e^{-c\sqrt{\mu}|x-y|}}{|x-y|^{d-1}}.
\end{equation}
For variable-coefficients, the only \emph{generic} bound which holds for the gradient of the elliptic Green function is another consequence of the De Giorgi-Nash-Moser theory: There exists $0<\alpha\le 1$ depending only on $\lambda$ and $d$ (with $\alpha \uparrow 1$ as $\lambda \uparrow 1$) such that for all $x,y\in \R^d$ 
\begin{equation}\label{eq:grad-Green-DGNM}
\mbox{if }|x-y|\gtrsim 1, \mbox{ then}\quad 
 \int_{B_1(x)}|\nabla_x G_\mu (x,y;A)|dx \,\le\, C \frac{e^{-c\sqrt{\mu}|x-y|}}{|x-y|^{d-2+\alpha}},
\end{equation}
see Lemma~\ref{lem:quenched} below. As can be seen, there is a mismatch between the generic behavior and the fundamental solution of the Laplacian at the level of the gradient.
The behavior at the singularity $x=y$ can only be described for smooth coefficients (say, uniformly H\"older-continuous). In that case, the optimal scaling of \eqref{eq:grad-Green-lap} holds for $|x-y|\lesssim 1$, cf.~\cite[Theorem 3.3]{Grueter-Widman-82} for $\mu=0$. However, even for analytic coefficients, the estimate \eqref{eq:grad-Green-lap} cannot hold \emph{generically} \emph{in the large}, that is in the regime $|x-y|\uparrow +\infty$, for this would contradict the counterexamples from quasiconformal mappings already mentioned.  

\medskip

In order to deal with measurable coefficients we need to consider local square averages, and shall make use of the following notation:
For all $L>0$ and all $|x-y|\ge 3L$ we set
\begin{align}\label{nabla_G_L}
(\nabla G_\mu)_L(x,y) &:= \Bigg( \fint_{B_L(x)} |\nabla_{x'} G_\mu(x',y)|^2 \;dx' \Bigg)^{\frac{1}{2}}\\
 (\nabla \nabla G_\mu)_L(x,y) &:= \Bigg( \fint_{B_L(x)} \fint_{B_L(y)} |\nabla \nabla G_\mu(x',y')|^2 \;dy'dx' \Bigg)^{\frac{1}{2}},
  \label{nabla2_G_L}
\end{align}
where (here and in the whole article) $\nabla \nabla$ stands for the second \emph{mixed} derivative $\nabla_{x'}\nabla_{y'}$. 
\begin{theorem}\label{thm:main}
Let the ensemble be stationary and satisfy (LSI) with constants $\rho$ and $\ell$.
Then there exists a random field $\calY$ with bounded finite moments such that 
for all $x\in \R^d$ with $|x|\ge 3\ell$ and all $\mu\ge 0$ we have
\begin{align}\label{nabla_est}
(\nabla G_\mu)_\ell(x,0)&\le \,\calY(x) \frac{e^{-c\sqrt{\mu}|x|}}{|x|^{d-1}},\\
 (\nabla\nabla G_\mu)_\ell(x,0) &\le \, \calY(x) \frac{e^{-c\sqrt{\mu}|x|}}{|x|^{d}}\label{nabla2_est}.
\end{align}
In addition, the random field $\calY$ has the following boundedness property: For all $1\le q<\infty$ there exists $C_q<\infty$ depending only on  $\lambda,p,\rho,\ell$ such that
\begin{equation}\label{eq:Y2}
 \sup_{|x|\ge 3\ell} \expec{\calY(x)^q}^\frac{1}{q}\,\leq\, C_q.
\end{equation}
\qed
\end{theorem}
\begin{remark}\label{cor:C^gamma}
If in addition the coefficients are H\"older-continuous, then the estimates 
of Theorem~\ref{thm:main} hold at the singularity as well, that is, 
~\eqref{nabla_est} \&~\eqref{nabla2_est} hold true for all $x\in \R^d$ and \eqref{eq:Y2} is replaced by
\begin{equation*}
\sup_{x\in \R^d}  \expec{\calY(x)^q}^\frac{1}{q}\,\leq\, C_q.
\end{equation*}
\end{remark}
Note that by stationarity the above result implies a similar decay for $(\nabla G_\mu)_\ell(x,y)$ for arbitrary $x,y\in\R^d$.

\medskip

This result is based on and extends the annealed estimates by Delmotte and Deuschel \cite{Delmotte-Deuschel-05}, see Proposition~\ref{prop:DD} below.
It is the extended continuum version of the result by Otto and the second author in  \cite{Marahrens-Otto-13} for discrete elliptic equations. 
At the cost of a slightly smaller decay rate, one may take the random field $\calY$ independent of $x$
in~\eqref{nabla_est} \&~\eqref{nabla2_est}:
\begin{corollary}\label{cor:quenched}
Under the assumptions of Theorem~\ref{thm:main}, for all $\beta>0$ there exists a random variable $\calY_\beta$
with bounded finite moments such that for all $|x|\ge 3\ell$ we have
\begin{eqnarray*}
(\nabla G_\mu)_\ell(x,0) &\dps \le \,\mathcal{Y}_\beta \frac{e^{-c\sqrt{\mu}|x|}}{|x|^{d-1-\beta}},
\quad 
(\nabla \nabla G_\mu)_\ell(x,0)  &\le \, \mathcal{Y}_\beta \frac{e^{-c\sqrt{\mu}|x|}}{|x|^{d-\beta}}.
\end{eqnarray*}
\qed
\end{corollary}

\subsection{Estimates of fluctuations}

Combined with a sensitivity estimate, the optimal gradient bounds on the Green functions allow us to quantify the fluctuations of solutions of linear elliptic equations with ``noisy" diffusion coefficients (a quantification of the propagation of uncertainty
in elliptic PDEs).
More precisely we consider diffusion coefficients $A_\e$ on $\R^d$ of the form
$$
A_\e(x)\,:=\,\Id+B(\frac{x}{\e})
$$
where $B$ is a random perturbation which has order $1$, correlation-length unity (which we shall replace in the theorem by the (LSI) assumption),
and vanishing expectation.
Hence, $A_\e$ is a perturbation of the identity by some noise of correlation-length $\e$.
Let $f$ be some RHS, and consider the random solution $u_\e$ of 
\begin{equation*}
u_\e-\nabla \cdot A_\e \nabla u_\e\,=\,f \qquad \text{in }\R^d.
\end{equation*}
The question we are interested in is the characterization of the fluctuations of $u_\e$ in function of $\e$ and of the statistics of $B$, first in terms of scaling and second in terms of law. In this contribution we address the question of the scaling wrt $\e$, and give optimal estimates of both weak and strong measures of the fluctuation,
which generalize the bounds obtained for $B$ small (that is, in the regime of small ellipticity ratio) 
by the first author in \cite{Gloria-12c}\footnote{Note that the proof of \cite[Lemma~2.1]{Gloria-12c} is wrong under the general assumption of finite correlation-length. The assumption of  \cite[Theorem~3]{Gloria-12c} should be replaced by ``Assume that the stationary random field $B$ satisfies spectral gap", as it is the case for Poisson inclusions for instance. The optimal form of \cite[Theorem~3]{Gloria-12c} is given by Theorems~\ref{thm:strong_error} and~\ref{thm:weak_error} below --- the norms in  \cite[Theorem~3]{Gloria-12c}  have to be adapted accordingly.}. 
The natural norms which control these fluctuations are mixed norms $L^{p}_{\lambda,\e}(\R^d)$ which 
measure local fluctuations at scale $\e$ in $L^\lambda$ but large scale fluctuations in $L^p$. In particular,
for all $q,\lambda\ge 1$, $\e>0$ and $f\in L^1_\loc(\R^d)$ we set
\begin{equation}\label{eq:mixed-norm}
\|f\|_{L^q_{\lambda,\e}(\R^d)} \,:=\, \bigg(\int_{\R^d} \Big(\fint_{B_\e(x)} |f(y)|^\lambda dy\Big)^{\frac{q}{\lambda}}dx\bigg)^{\frac{1}{q}}.
\end{equation}
In particular it is bounded by the $L^q(\R^d)$-norm for $q\ge \lambda$ by Jensen's inequality.
We start with the estimate of the fluctuations in a strong norm.
\begin{theorem}\label{thm:strong_error}
Let  $A_\e=A(\frac{\cdot}{\e})$ be the $\e$-rescaling of the coefficient field $A\in \Omega$ distributed according to a stationary ensemble $\expec{\cdot}$ that satisfies (LSI).
Let $\mu\ge 0$.
For all $\e>0$, let $u_\e\in H^1(\R^d)$ be a distributional solution of 
\begin{equation}\label{eq:ueps-fluctu}
\mu u_\e-\nabla \cdot A_\e \nabla u_\e \,=\,f \qquad \text{in }\R^d.
\end{equation}
Then for all $\lambda>\frac{d}{2}$, $1\le \theta<\infty$, $2\le p <\infty$, $1\le r \le \frac{d}{d-1}$, and $q$ such that
\begin{equation}\label{params_strong}
 1 + \frac{1}{p} = \frac{1}{r} + \frac{1}{q},
\end{equation}
the fluctuations of $u_\e$ satisfy
\begin{equation*}
 \Big\langle \bigg( \int_{\R^d} |u_\e - \expec{u_\e} |^p \;dx \bigg)^\theta \Big\rangle^{\frac{1}{p\theta}} \, \lesssim \, \left\{
\begin{array}{lll}
d=2&:& |\ln (\mu\e^2)|^{\frac{1}{2}}+1\\
d>2&:&1
\end{array}
\right\} \e (\mu^{-\frac{1-d}{2}-\frac{d}{2r}}+1) \|f\|_{L^{q}_{\lambda,\e}(\R^d)}
\end{equation*}
where $\|f\|_{L^{q}_{\lambda,\e}(\R^d)}$ is given by~\eqref{eq:mixed-norm}.
In the border-line case $r=\frac{d}{d-1}$, we require in addition $q>1$.
\qed
\end{theorem}
We then turn to the estimate of weak norm of the fluctuations.
\begin{theorem}\label{thm:weak_error}
Let  $A_\e=A(\frac{\cdot}{\e})$ be the $\e$-rescaling of the coefficient field $A\in \Omega$ distributed according to a stationary ensemble $\expec{\cdot}$ that satisfies (LSI).
Let $\mu\ge 0$.
For all $\e>0$, let $u_\e\in H^1(\R^d)$ be a distributional solution of \eqref{eq:ueps-fluctu}.
Then 
for all $1\le \theta<\infty$, $2\le p <\infty$, $1\le r,\tilde r \le \frac{d}{d-1}$, $1\le q\le \frac{r}{r-1}$, and $1\le \tilde q\le \frac{\tilde r}{\tilde r-1}$ such that
\begin{equation}\label{params_weak}
 2 + \frac{1}{2} = \frac{1}{r} + \frac{1}{\tilde r} + \frac{1}{q} + \frac{1}{\tilde q}
\end{equation}
and for all $\lambda_1,\lambda_2\ge 1$ such that 
\begin{equation}\label{params_weak2}
 \frac{1}{\lambda_1} + \frac{1}{\lambda_2} < \frac{d+2}{d},
\end{equation}
the fluctuations of $u_\e$ satisfy for all $g\in L^1_\loc(\R^d)$,
\begin{equation*}
 \expec{ \Big| \int_{\R^d} (u_\e - \expec{u_\e}) g \;dx \Big|^\theta }^{\frac{1}{\theta}} \, \lesssim\, \e^\frac{d}{2} (\mu^{-(1-d)-\frac{d}{2}(\frac{1}{r} + \frac{1}{\tilde r})} +1)\|f\|_{L^{q}_{\lambda_2,\e}(\R^d)} \|g\|_{L^{\tilde q}_{\lambda_1,\e}(\R^d)}.
\end{equation*}
In the border-line case $r=\tilde r=\frac{d}{d-1}$, we require in addition $q,\tilde q >1$.
\qed
\end{theorem}
\begin{remark}
When the coefficients $A$ in Theorems~\ref{thm:strong_error} and~\ref{thm:weak_error} are uniformly H\"older continuous, then we can replace the mixed norms $L^{q}_{\lambda,\e}(\R^d)$ by the usual norms $L^q(\R^d)$ in the estimates. 
This shows that one can trade local integrability of $f$ and $g$ for regularity of $A$.
This is proved by replacing averaged bounds on the Green function by pointwise bounds, as in \cite{Marahrens-Otto-13}. We leave the details to the reader.
\qed
\end{remark}
\begin{remark}
Theorem~\ref{thm:weak_error} reveals the central limit scaling of the weak measure of the fluctuations.
While the most natural norms for the RHS on $\R^d$ are those which make the estimate independent of $\mu$,
the other estimates are valuable for $\mu>0$ since the massive term essentially localizes the equation to a bounded domain of size $\mu^{-\frac{1}{2}}$ (without boundary layers).
\qed
\end{remark}
These results generalize both \cite[Theorem~3]{Gloria-12c} and 
\cite[Corollaries 2~\& 3]{Marahrens-Otto-13} (cf. also \cite{Conlon-Naddaf-00} by Conlon and Naddaf in the case 
of discrete elliptic equations). 
Note that when the noise is in the zero-order term (that is, for $\mu$ replaced by $1+b_\e$ and $A_\e$ by $\Id$ in \eqref{eq:ueps-fluctu}), the CLT scaling (and in addition the characterization of the limiting law)
was established by Figari, Orlandi and Papanicolaou in \cite{Figari-Orlandi-Papanicolaou-82} for $d\ge 4$ and by Bal in \cite{Bal-08} for $d\le 3$. The arguments involved in the proof of Theorems~\ref{thm:strong_error} and~\ref{thm:weak_error} have a different flavor since the randomness is in the derivative of highest order.


\section{Structure of the proofs and auxiliary results}

We start with the definition and main properties of the elliptic Green function.
\begin{definition}[Green's function]\label{def:Green}
For all $A\in \Omega$ and every $\mu>0$, 
there exists a unique function $G_\mu(x,y;A)\ge 0$ with the following properties
\begin{itemize}
\item Qualitative continuity off the diagonal, that is,
\begin{equation}\label{gu.14}
\{(x,y)\in\mathbb{R}^d\times\mathbb{R}^d|x\not=y\}\ni(x,y)\mapsto G_\mu(x,y;A)\quad
\mbox{is continuous}.
\end{equation}
\item Upper pointwise bounds on $G_\mu$:
\begin{equation}\label{eq:ptwise-decay-estim}
G_\mu(x,y;A)\,\lesssim\,e^{-c\sqrt{\mu}|x-y|}
\left\{\begin{array}{ccc}
\ln(2+\frac{1}{\sqrt{\mu}|x-y|})&\mbox{for}&d=2\\
\frac{1}{|x-y|^{d-2}}&\mbox{for}&d>2
\end{array}\right\},
\end{equation}
where here and in the sequel the rate constant $c>0$ in the exponential
is generic and may change from term to term, but only depends on $d$ and $\lambda$.
\item Averaged bounds on $\nabla_x G_\mu$ and $\nabla_y G_\mu$:
\begin{eqnarray}
\left(R^{-d}\int_{R<|x-y|\le 2R}|\nabla_xG_\mu(x,y;A)|^2dx\right)^\frac{1}{2}&\lesssim& e^{-c\sqrt{\mu}R} R^{1-d},\label{gu.11}\\
\left(R^{-d}\int_{R<|y-x|\le 2R}|\nabla_yG_\mu(x,y;A)|^2dy\right)^\frac{1}{2}&\lesssim&  e^{-c\sqrt{\mu}R}R^{1-d}.\label{gu.12}
\end{eqnarray}
\item Differential equation:
We note that \eqref{eq:ptwise-decay-estim} and \eqref{gu.11} \& \eqref{gu.12} imply that the maps
$\mathbb{R}^d\ni x\mapsto(G_\mu(x,y;A),\nabla_xG_\mu(x,y;A))$ and 
$\mathbb{R}^d\ni y\mapsto(G_\mu(x,y;A),\nabla_yG_\mu(x,y;A))$ are (locally) integrable. Hence even for
discontinuous $A$, we may formulate the requirement
\begin{eqnarray}
\mu G_\mu-\nabla_x\cdot A(x)\nabla_xG_\mu=\delta(x-y)&&\mbox{distributionally in}\;\mathbb{R}_x^d,\label{eq:Green}\\
\mu G_\mu-\nabla_y\cdot A^*(y)\nabla_yG_\mu=\delta(y-x)&&\mbox{distributionally in}\;\mathbb{R}_y^d,\label{eq:Green-t}
\end{eqnarray}
where $A^*$ denotes the transpose of $A$.
\end{itemize}
We note that the uniqueness statement implies $G_\mu(x,y;A^*)=G_\mu(y,x;A)$ so that $G_\mu$ is symmetric 
when $A$ is symmetric.
\qed
\end{definition}
These standard properties of the massive Green functions are proved in \cite{Gloria-Otto-10b} (essentially following arguments of \cite{Grueter-Widman-82}).
\begin{remark}
All the main results of this article are stated for $\mu\ge 0$, whereas 
we shall only consider the case $\mu>0$ in the proofs. Indeed, one can pass to the limit as $\mu \downarrow 0$ in all our estimates, and local averages of $\nabla G_\mu$ and $\nabla \nabla G_\mu$ converge to local averages of $\nabla G$ and $\nabla \nabla G$, where $G$ is the Green function for $\mu=0$ (the existence of which is subtle for $d=2$).
\qed
\end{remark}

The improvement of the De Giorgi-Nash-Moser theory in the large is a consequence of the bounds on the Green function of Theorem~\ref{thm:main}.
As in the discrete case dealt with in \cite{Marahrens-Otto-13} the strategy is to upgrade to \emph{any} moment in probability the optimal bounds by Delmotte and Deuschel \cite{Delmotte-Deuschel-05} on the \emph{first} and \emph{second} moments of $\nabla \nabla G_\mu$ and $\nabla G_\mu$, respectively.
Yet, the bounds by Delmotte and Deuschel in \cite[Theorem~1.2]{Delmotte-Deuschel-05} are not enough at the level of the mixed second gradient, and we shall use the following result of \cite{Marahrens-Otto-14}
in its version with the massive term proved in \cite[Lemma~2.11]{Gloria-Otto-10b}:
\begin{proposition}\label{prop:DD}
If the ensemble is stationary, then the Green function satisfies for all $\mu> 0$, all $L\gtrsim 1$,
and all $x\in \R^d$ with $|x|\ge 2L$,
\begin{align}
 \expec{ (\nabla_x G_\mu)_L(x,0)^2}^{\frac{1}{2}} &\le C \frac{e^{-c\sqrt{\mu} |x|}}{|x|^{d-1}},\label{DD1}\\
 \expec{ (\nabla\nabla G_\mu)_L(x,0)} &\le C \frac{e^{-c\sqrt{\mu} |x|}}{|x|^{d}},\label{DD2}
\end{align}
for some constants $C$ and $c$ depending only on $\lambda$ and $d\ge 2$.
\qed
\end{proposition}
Estimate~\eqref{nabla2_est} of Theorem~\ref{thm:main} is a consequence of \eqref{DD2} and of the following reverse H\"older estimate valid for all $p\ge 1$ large enough:
\begin{multline}\label{eq:explain-strategy}
 \sup_{x,y:|x-y|\ge 6\ell} \Big\{ |x-y|^{d}e^{c\sqrt{\mu} |x-y|} \langle |(\nabla\nabla G_\mu)_\ell(x,y)|^{2p} \rangle^{\frac{1}{2p}} \Big\}\\
 \le C(d,\lambda,p,\rho,\ell) \sup_{x,y:|x-y|\ge 6\ell} \Big\{ |x-y|^{d}e^{c\sqrt{\mu} |x-y|}  \langle |(\nabla\nabla G_\mu)_\ell(x,y)| \rangle \Big\},
\end{multline}
and likewise for the first derivative.
This gain of integrability is achieved by the following lemma in the spirit of \cite[Lemma 4]{Marahrens-Otto-13}, where the assumption that $\expec{\cdot}$ satisfies (LSI) is crucial.
\begin{lemma}\label{lem:LSI_p}
Let $\langle\cdot\rangle$ satisfy (LSI) with constants $\rho,\ell>0$. Then for
arbitrary $\delta>0$ and $1\le p < \infty$ and for any random variable $\zeta$ 
we have
\begin{equation}
\expec{|\zeta|^{2p}}^\frac{1}{2p} \le C(d,\rho,p,\delta)\expec{|\zeta|}
+\delta\expec{\bigg(\int_{\R^d} \Big(\osc{A|_{B_\ell(z)}}{\zeta} \Big)^2 \;dz\bigg)^p}^\frac{1}{2p}
\label{LSI_p}
\end{equation}
for some finite constant $C(d,\rho,p,\delta)$, where we recall that the expectation in the RHS is an outer expectation.
\qed
\end{lemma}
Since $G_\mu$ is measurable on $\Omega$, one may apply this lemma to $\zeta = (\nabla\nabla G_\mu)_\ell(x,0)$ and $\zeta = (\nabla_x G_\mu)_\ell(x,0)$.
In order to prove the reverse H\"older
inequality \eqref{eq:explain-strategy}, it suffices to absorb the second RHS term of  \eqref{LSI_p} in the RHS. This is the content of the following lemma, which is essentially based on deterministic arguments.
\begin{lemma}[Absorption lemma]\label{lem:absorb}
Let $d\ge 2$. There exists $p_0\ge 1$ depending only on $\lambda$ and $d$ such that 
for all $L\sim 1$ and $p\ge p_0$, we have for the second derivative:
\begin{multline}
 \sup_{|x-y|\ge 6L} \bigg\{ |x-y|^{2pd} e^{2pc\sqrt{\mu} |x-y|}\expec{\bigg(\int_{\R^d} \Big(\osc{A|_{B_L(z)}}{(\nabla\nabla G_\mu)_L(x,y)} \Big)^2\;dz
\bigg)^p } \bigg\}\\
 \lesssim \sup_{|x-y|\ge 6L} \Big\{ |x-y|^{2pd}e^{2pc\sqrt{\mu} |x-y|} \expec{\big((\nabla\nabla G_\mu)_L(x,y)\big)^{2p} }\Big\} + 1,
\label{absorb}
\end{multline}
anf for the first derivative:
\begin{multline}
 \sup_{|x-y|\ge 6L} \bigg\{ |x-y|^{2p(d-1)} e^{2pc\sqrt{\mu} |x-y|}\expec{\bigg( \int_{\R^d} \Big(\osc{A|_{B_L(z)}}{(\nabla_x G_\mu)_L(x,y)} \Big)^2
\;dz\bigg)^p } \bigg\}\\
 \lesssim \, \sup_{|x-y|\ge 6L} \Big\{ |x-y|^{2p(d-1)} e^{2pc\sqrt{\mu} |x-y|}\expec{\big((\nabla_x G_\mu)_L(x,y)\big)^{2p}}  \Big\}\\ + \sup_{|x-y|\ge 6L} \Big\{
|x-y|^{2pd}e^{2pc\sqrt{\mu} |x-y|}  \expec{\big((\nabla\nabla G_\mu)_L(x,y)\big)^{2p}} \Big\} + 1,
\label{absorb1}
\end{multline}
where $\lesssim$ stands for $\le$ up to a multiplicative constant which depends on $d$, $\lambda$, and $p$.
\qed
\end{lemma}
A key ingredient to the proof of Lemma \ref{lem:absorb} are the following deterministic estimates.
\begin{lemma}\label{lem:quenched}
Let $d\ge 2$. There exist $q_0>1$ and $\alpha_0>0$ depending only on $d$ and $\lambda>0$ such that
for all $\mu>0$, $1\le q\le q_0$, 
and all $R \ge 4L \sim 1$,
\begin{align}
\int_{R\le|x-y|<2R} |\nabla_x G_\mu(x,y)|^{2q} \;dx &\lesssim {R}^{d+(1-d)2q}e^{-c\sqrt{\mu} R},\label{quenched_ap1}\\
 \int_{R\le|x-y|<2R} \int_{|y|<L} |\nabla \nabla G_\mu(x,y)|^{2q} \;dydx &\lesssim R^{-2q\alpha_0} e^{-c\sqrt{\mu} R} ,\label{quenched_ap}
\end{align}
where the multiplicative constants depend only on $d$ and $\lambda$.
In addition we have the following local boundedness estimate for all $L\sim 1$
\begin{equation}\label{loc_bound}
 \sup_{x,y\in\R^d:3L \le |x-y| < 6L} \big\{ (\nabla \nabla G_\mu)_L(x,y) \big\} \lesssim 1.
\end{equation}
\qed
\end{lemma}
%
\begin{remark}\label{rem:G_vs_GL}
(i) Our results beg the question if we can upgrade \eqref{nabla_est} and \eqref{nabla2_est} to a stronger version without space integrals as
in \eqref{DD1} and \eqref{DD2}. The answer is negative if $p>1$. Let us consider \eqref{DD1} (in the parabolic setting, \eqref{DD2}
directly follows from \eqref{DD1}). 
Using the De Giorgi-Nash-Moser theory, we may upgrade \eqref{nabla_est} to pointwise-estimates away from
the singularity if $p=1$, but not otherwise.
Indeed, the De Giorgi-Nash-Moser theory yields away from the singularity that
$$
 \Big\langle \int_{B_L(x)} \int_{B_L(y)} |\nabla_1 G(x',y')|^2
\;dy'dx' \Big\rangle \gtrsim \Big\langle \int_{B_L(x)} |\nabla_1 G(x',y)|^2\;dx' \Big\rangle.
$$
Now by stationarity, the left hand side equals
\begin{align*}
 \Big\langle \int_{B_L(0)} |\nabla_1 G(x+x',y)|^2 \;dx' \Big\rangle &= \Big\langle \int_{B_L(0)} |\nabla_1 G(x-y,-x')|^2 \;dx' \Big\rangle\\
 &\gtrsim \Big\langle |\nabla_x G(0,y-x)|^2 \Big\rangle,
\end{align*}
where the last inequality again follows from de Giorgi-Nash-Moser theory.
On the other hand, if $p>1$, pointwise bounds on $\langle |\nabla G|^{2p} \rangle$ cannot be expected since there is no local regularity to
control $\langle \int_{B_L} |\nabla G|^{2p} dx \rangle$. On the other hand, clearly energy methods allow to control locally the
$L^2$--norms of the gradient, which shows why $\langle |\nabla G|^2 \rangle$ may indeed be bounded. In other
words, the spatial integrals in \eqref{nabla_est} and \eqref{nabla2_est} are necessary to smooth out local effects when the coefficients
lack regularity if and only if $p>1$.\\
(ii) In a similar spirit, we observe that the restriction $|x|\gtrsim L$ is not necessary in \cite{Delmotte-Deuschel-05}, but cannot be avoided here. Indeed,
assuming Proposition \ref{prop:DD} only for $|x|\gtrsim 1$, we may remove this restriction by a simple scaling argument. The same is true
if we (could) replace $(\nabla\nabla G)_L$ by $\nabla\nabla G$ as discussed in (i). On the other hand, the presence of the averaging
operation $(\cdot)_L$ breaks the scaling invariance by introducing a length scale $L$. Therefore we cannot expect to obtain information on
the blow-up of $(\nabla \nabla G)_L(x,y)$ as the singularity enters the integral, i.e.\ as $|x-y|\downarrow 2L$.
\qed
\end{remark}

\medskip

We turn now to the fluctuation estimates.
By a scaling argument, it is enough to prove Theorems~\ref{thm:strong_error} and~\ref{thm:weak_error} for $\e=1$ and $\ell=\frac{1}{2}$.
We thus consider the solution $u\in H^1(\R^d)$ of 
\begin{equation}\label{u}
\mu u - \nabla \cdot A \nabla u\, =\, f,\qquad\mu\ge0.
\end{equation}
We shall only consider the case $\mu>0$ in the proofs. The results for $\mu=0$ are then obtained by letting $\mu\downarrow 0$ in the estimates.
The starting point is the following spectral gap estimate
\begin{lemma}[$q$-(SG)]\label{lem:var-estim} 
If $\expec{\cdot}$ satisfies (LSI) with amplitude $\rho>0$ and correlation-length $\ell<\infty$, then we have for all $q\ge 1$
and all random variables $\zeta$
\begin{equation}\label{eq:var-estim}
\expec{(\zeta-\expec{\zeta})^{2q}}^{\frac{1}{q}}\;\lesssim \;\left\langle\Big(\int_{\mathbb{R}^d}
\Big(\osc{A|_{B_{\tilde \ell}(z)}}{\zeta}\Big)^2dz\Big)^q\right\rangle^{\frac{1}{q}},
\end{equation}
with $\tilde \ell=2\ell$, where the multiplicative constant depends on $q$ and $\rho$.
\qed
\end{lemma}
This is a standard result. It is indeed enough to assume that $\expec{\zeta}=0$ and $\expec{\zeta^2}=1$.
To prove estimate~\eqref{eq:var-estim} for $q=1$ it suffices to apply (LSI) to the random variable $\chi=\sqrt{1-\alpha^2}+\alpha \zeta$ and make a Taylor expansion as $\alpha \downarrow 0$, this yields
the result for the correlation-length $\ell$. The estimate for $q>1$ is a consequence of the estimate for $q=1$ (up to increasing $\ell$ to $\tilde \ell=2\ell$),
see for instance \cite[Corollary~2.3]{Gloria-Otto-10b}.
Since we have assumed that $\ell=\frac{1}{2}$, \eqref{lem:var-estim} holds for $\tilde \ell=1$.

The following lemma is a sensitivity estimate which quantifies how much the solution $u$ of \eqref{u} depends on the coefficients $A$.
\begin{lemma}\label{lem:osc_u}
Let $\lambda_1, \lambda_2\in [1,+\infty]$ satisfy 
\begin{equation}\label{lambda_i}
 \frac{1}{\lambda_1} + \frac{1}{\lambda_2} < \frac{d+2}{d} \quad\Leftrightarrow\quad \frac{1}{\lambda_1'} + \frac{1}{\lambda_2'} > \frac{d-2}{d}.
\end{equation}
In particular, at most one of $\lambda_1'$, $\lambda_2'$ may be infinite.
Denote by $u(\cdot;A)\in H^1(\R^d)$ the solution of 
 \eqref{u}. We use the short-hand notation $\tilde u$ for $u(\cdot;\tilde A)$, $\tilde A\in \Omega$.
We then have that
\begin{equation}\label{eq:osc-estim-u}
 \sup_{\tilde A : \tilde A|_{\R^d \setminus B_\ell(z)} = A|_{\R^d \setminus B_\ell(z)}} \| u - \tilde u \|_{L^{\lambda_1'}(B_\ell(x))} \,\lesssim \, \mathcal{K}_{G_\mu,u}(x,z),
\end{equation}
where
\begin{equation}\label{K_G}
 \mathcal{K}_{G_\mu,u}(x,z) \,:= \,
\begin{cases}
(\nabla G_\mu)_{2\ell}(x,z) (\nabla u)_\ell(z) &\text{if } |x-z| > 6\ell,\\
\| f \|_{L^{\lambda_2}(B_{2\ell}(x))} + (\nabla u)_{9\ell}(z) &\text{if } |x-z| \le 6\ell.
\end{cases}
\end{equation}
If $\lambda_1' = +\infty$, we reformulate this result in the pointwise form
$$
 \osc{A|_{B_\ell(z)}}{u(x)} \,\lesssim \,\mathcal{K}_{G_\mu,u}(x,z).
$$
\qed
\end{lemma}
In the proof of Lemma~\ref{lem:osc_u} we shall make use of the following standard result.
\begin{lemma}\label{lem:DGNM}
Let $p,q\in [1,+\infty]$ satisfy 
$$
 \frac{1}{q'} + \frac{1}{p} > \frac{d-2}{d} \quad\Leftrightarrow\quad \frac{1}{q} < \frac{1}{p} + \frac{2}{d}.
$$
If $u$ is a solution of~\eqref{u} in $B_2=B_2(0)$, we have that
$$
 \| u \|_{L^p(B_1)} \lesssim \| u \|_{L^2(B_2)} + \| f \|_{L^q(B_2)},
$$
where the multiplicative constant depends on $\lambda$, $d$ and $q$, but not on $\mu\ge 0$.
\qed
\end{lemma}
This result is usually stated for $p=\infty$ only, cf.~\cite[Theorem~8.17]{Gilbarg-Trudinger-98}.
Although we think it should follow from the Nash-Aronson bounds (if $d>2$), Young's inequality and the well-known estimate with $p=+\infty$, we display a direct proof for $p<\infty$ using a (simplified) Moser-type iteration that works for $d=2$ and uses less machinery.


\section{Proofs of the estimates on the Green functions}

\subsection{Proof of Theorem~\ref{thm:main}}

The proof is a simple combination of Proposition~\ref{prop:DD}, Lemma~\ref{lem:LSI_p} and Lemma~\ref{lem:absorb}.

\medskip

\step{1} Proof of~\eqref{nabla2_est}. 

\noindent We apply~\eqref{LSI_p} of Lemma~\ref{lem:LSI_p} to $\zeta(A) = (\nabla\nabla G_\mu)_\ell(A;x,y)$ for some $x,y \in\R^d$ such that $|x-y| \ge 6\ell$ to the effect of
\begin{multline*}
 \expec{ (\nabla\nabla G_\mu)_\ell(x,y)^{2p}}^{\frac{1}{2p}} \\
 \le\,  C(d,\rho,\ell,p,\delta) \expec{ (\nabla\nabla G_\mu)_\ell(x,y) } + \delta \expec{ \bigg( \int_{\R^d} \Big(\osc{A|_{B_L(z)}}{(\nabla\nabla G_\mu)_L(x,y)} \Big)^2 dz \bigg)^p}^{\frac{1}{2p}}.
\end{multline*}
Combined with \eqref{DD2} in Proposition~\ref{prop:DD} this yields
\begin{multline*}
 |x-y|^d e^{c\sqrt{\mu} |x-y|} \expec{ (\nabla\nabla G_\mu)_\ell(x,y)^{2p} }^{\frac{1}{2p}}\, \le \, C(d,\lambda,\rho,\ell,p,\delta)\\ + \delta |x-y|^d e^{c\sqrt{\mu} |x-y|}\expec{ \bigg( \int_{\R^d} \Big(\osc{A|_{B_\ell(z)}}{(\nabla\nabla G_\mu)_\ell(x,y)} \Big)^2 \;dz \bigg)^p }^{\frac{1}{2p}}.
\end{multline*}
We take the supremum over all $x$ and $y$ such that $|x-y|\ge6\ell$ and insert \eqref{absorb} in Lemma~\ref{lem:absorb} to obtain that
\begin{multline*}
 \sup_{|x-y|\ge 6\ell} \Big\{ |x-y|^d e^{c\sqrt{\mu} |x-y|}\expec{ (\nabla\nabla G_\mu)_\ell(x,y)^{2p}}^{\frac{1}{2p}} \Big\} \,\le\, C(d,\lambda,\rho,\ell,p,\delta)\\ + C(d,\lambda,p,\ell) \delta \sup_{|x-y|\ge 6\ell} \Big\{ |x-y|^d e^{c\sqrt{\mu} |x-y|} \expec{ (\nabla\nabla G_\mu)_\ell(x,y)^{2p} }^{\frac{1}{2p}} +1\Big\}.
\end{multline*}
Choosing $\delta$ small enough, we may absorb the last RHS term in the LHS. This yields~\eqref{nabla2_est}.

\medskip

\step{2} Proof of~\eqref{nabla_est}.

\noindent We proceed as in Step~1: Take $\zeta(A) = (\nabla G_\mu)_\ell(A;x,y)$ in Lemma~\ref{lem:LSI_p} to deduce
$$
 \expec{ (\nabla G_\mu)_\ell(x,y)^{2p} }^{\frac{1}{2p}}\, \le\, C(d,\rho,\ell,p,\delta) \expec{ (\nabla G_\mu)_\ell(x,y) }  + \delta \expec{ \bigg( \int_{\R^d} \Big(\osc{A|_{B_\ell(z)}}{(\nabla G_\mu)_\ell(x,y)} \Big)^2 dz \bigg)^p }^{\frac{1}{2p}}.
$$
Combined with \eqref{DD1} in Proposition~\ref{prop:DD}, this turns into
\begin{multline*}
 |x-y|^{d-1} e^{c\sqrt{\mu} |x-y|} \expec{ (\nabla G_\mu)_\ell(x,y)^{2p} }^{\frac{1}{2p}} \,\le\, C(d,\lambda,\rho,\ell,p,\delta)\\ + \delta |x-y|^{d-1} e^{c\sqrt{\mu} |x-y|}  \expec{ \bigg( \int_{\R^d} \Big(\osc{A|_{B_\ell(z)}}{(\nabla G_\mu)_\ell(x,y)} \Big)^2 \;dz \bigg)^p}^{\frac{1}{2p}}.
\end{multline*}
After taking the supremum over all $x,y$ such that $|x-y|\ge6\ell$, the estimate~\eqref{absorb1} from Lemma~\ref{lem:absorb} yields
\begin{multline*}
 \sup_{|x-y|\ge6\ell} \bigg\{ |x-y|^{d-1} e^{c\sqrt{\mu} |x-y|}  \expec{\bigg( \int_{\R^d} \Big(\osc{A|_{B_\ell(z)}}{(\nabla G_\mu)_\ell(x,y)} \Big)^2 \;dz \bigg)^p }^{\frac{1}{2p}}\\ 
 \le \, C(d,\lambda,p,\ell) \Bigg(1 + \delta\sup_{|x-y|\ge 6\ell} \bigg\{ |x-y|^{2p(d-1)} e^{2pc\sqrt{\mu} |x-y|} \expec{\big|(\nabla_x G_\mu)_\ell(x,y)\big|^{2p} } \bigg\} \Bigg)\\ 
 + C(d,\lambda,p,\ell) \delta \sup_{|x-y|\ge 6\ell} \bigg\{ |x-y|^{2pd} e^{2pc\sqrt{\mu} |x-y|} \expec{ \big|(\nabla\nabla G_\mu)_\ell(x,y)\big|^{2p}} \bigg\}.
\end{multline*}
By \eqref{nabla2_est} (proved in Step~1), the last term is bounded by a constant $C(d,\lambda,\rho,\ell,p)\delta$. We then conclude by taking $\delta$ small enough so that we can absorb the remaining supremum on the LHS.
The desired estimates \eqref{nabla_est} and~\eqref{nabla2_est} then follow from the definition 
$$
\calY(x)\,:=\, \max\{(\nabla G_\mu)_\ell(x,0) |x|^{d-1}e^{c\sqrt{\mu}|x|}, (\nabla\nabla G_\mu)_\ell(x,0)|x|^{d} e^{c\sqrt{\mu}|x|}\}.
$$

\subsection{Proof of Corollary~\ref{cor:quenched}}

For every $x\in\R^d$, there exists some $x'\in \frac{\ell}{\sqrt{d}}\Z^d$ such that the difference $x-x'$ has max-norm $|x-x'|_\infty \le \frac{\ell}{2\sqrt{d}}$. Hence its Euclidean norm satisfies $|x-x'|\le \frac{\ell}{2}$. Consequently, we have that $|(\nabla \nabla G_\mu)_\ell(x',0)| \ge |(\nabla \nabla G_\mu)_{\frac{\ell}{4}}(x,0)|$ and it holds that
\begin{align*}
 \Big\langle \sup_{x\in\R^d\setminus B_{4\ell}} \big\{ |x|^{d-\beta} |(\nabla\nabla G_\mu)_{\frac{\ell}{4}}(x,0)| \big\}^p \Big\rangle &\le \Big\langle \sup_{x'\in(\frac{\ell}{\sqrt{d}}\Z)^d, |x'|\ge 4\ell} \big\{ |x'|^{d-\beta} |(\nabla\nabla G_\mu)_{\frac{\ell}{4}}(x',0)| \big\}^p \Big\rangle\\
 &\le \sum_{x'\in(\frac{\ell}{\sqrt{d}}\Z)^d, |x'| \ge 4\ell} |x'|^{-\beta p} \langle |x'|^{dp} |(\nabla \nabla G_\mu)_\ell(x',0)|^p \rangle\\
 &\le C(d,\lambda,\rho,\ell,\gamma,\beta)
\end{align*}
as long as $\beta p > d$, which we may assume without loss of generality since by Jensen's inequality we may always increase $p$.
The same remark applies to $(\nabla G_\mu)_\ell$.
The choice
$$
\mathcal{Y}_\beta\,:=\,\max\Big\{\sup_{x\in\R^d\setminus B_{4\ell}}\big\{ |x|^{d-\beta} |(\nabla\nabla G_\mu)_{\frac{\ell}{4}}(x,0)| \big\},\sup_{x\in\R^d\setminus B_{4\ell}}\big\{ |x|^{d-1-\beta} |(\nabla G_\mu)_{\frac{\ell}{4}}(x,0)| \big\}\Big\}
$$
concludes the proof.

\subsection{Proof of Remark \ref{cor:C^gamma}}\label{ssec:C^gamma}

We split the proof into two steps.

\medskip

\step{1} Near-field estimates.

\noindent The results of \cite[Theorem 3.3]{Grueter-Widman-82} yield
$$
 |\nabla G_\mu(x,0)| \lesssim |x|^{1-d} \quad\text{and}\quad |\nabla\nabla G_\mu(x,0)| \lesssim |x|^{-d}
$$
for all $|x|\le 3\ell$.
(The fact that $G_\mu$ does not vanish on $\partial B_{3\ell}$ can be dealt with by substracting the corresponding boundary value problem, which is clearly bounded by the classical Schauder estimates and the Nash-Aronson $L^\infty$-estimate on $G_\mu$ away from the origin. The arguments are uniform wrt $\mu\ge 0$. The estimate for $d=2$ can be deduced from the corresponding estimate for $d=3$
by using the elegant argument by Avellaneda and Lin \cite{Avellaneda-Lin-87}, see for instance Step~2 of the proof of Lemma~\ref{lem:quenched} below.)

\medskip

\step{2} Far-field estimates.

\noindent
It remains to treat the $|x|\ge 3\ell$. Let $u$ be a $(\mu-\nabla \cdot A\nabla)$-harmonic function in $\R^d\setminus B_\ell$. Our goal is to prove the following reverse H\"older inequality
\begin{equation}\label{L^2_L^infty}
 |\nabla u(x)|^2 \lesssim \int_{B_{\ell}(x)} |\nabla u(x')|^2 \;dx'+\mu \fint_{B_\ell(x)} |u(x)|dx,
\end{equation}
with a constant depending on $\ell$, $d$, $\lambda$, and $\gamma$ only.
Without the derivative, this is a consequence of the De Giorgi-Nash-Moser theory. Since we are interested in $\nabla u$, we require the H\"older-continuity of the coefficient field. 
In the following, we will nonetheless pursue a strategy similar to Moser iteration to achieve the desired bound in \eqref{L^2_L^infty}. 
Since $A$ is H\"older-continuous, the function $u$ satisfies $u\in C^{2,\gamma}(\R^d \setminus \overbar{B_\ell})$ by interior Schauder theory. 
Now consider some length $0<L\le \frac{\ell}{2}$, and denote by $\ext u_L$ the average of $u$ on $B_L(x)$.
Let $\eta\in C^\infty_0(B_L(x))$. By assumption, we have that $\eta (\mu u-\nabla \cdot A \nabla u) = 0$ in $\R^d$. Fix some $y'\in B_L(x)$. The product rule yields
\begin{multline}\label{eq:eta_u^gamma}
\mu(\eta (u-\ext u_L))(y) -\nabla \cdot A(y') \nabla (\eta (u-\ext u_L))(y) 
\\ =\, 
-\mu \ext u_L \eta(y) +
\nabla \cdot \big( ( A(y) - A(y') ) \eta(y) \nabla (u-\ext u_L)(y) \big)
\\
 - \nabla \cdot ((u-\ext u_L)(y) A(y') \nabla \eta(y)) - \nabla \eta(y) \cdot A(y) \nabla (u-\ext u_L)(y)
\end{multline}
for all $y\in\R^d$. This is a constant-coefficient elliptic equation in $y$ with a right hand side in $H^{-1}(\R^d)$ and associated Green function $G_0(\cdot) \equiv G_\mu(\cdot,0;A(y'))$. The Green function representation yields for all $x'\in B_{\frac{L}{2}}(x)$
\begin{multline}\label{eq:rep_eta_u^gamma}
 (\eta (u-\ext u_L))(x') = \int_{\R^d} \Big( \nabla G_0(x'-y) \cdot (A(y') - A(y)) \nabla (\eta (u-\ext u_L))(y)\\
 + (u-\ext u_L)(y) \nabla G_0(x'-y) \cdot A(y) \nabla \eta(y) \\
 + G_0(x'-y)\big( \nabla \eta(y) \cdot A(y) \nabla (u-\ext u_L)(y)-\mu \ext u_L \eta(y)\big)\Big) \;dy.
\end{multline}
This can be made rigorous by mollification of the RHS of~\eqref{eq:eta_u^gamma}. Indeed, since $u\in C^{2,\gamma}(\overbar{B_L(x)})$, the limit exists and is given by~\eqref{eq:rep_eta_u^gamma}. 
Assume now that $\eta$ is a cutoff function for $B_{\frac{2L}{3}}(x)$ in $B_{L}(x)$
such that $|\nabla \eta|\lesssim \frac{1}{L}$. 
We may also take the gradient in \eqref{eq:rep_eta_u^gamma} w.~r.~t.~$x'$ at the point $y' \in B_{\frac{2L}{3}}(x) $ to obtain
\begin{multline*}
 \nabla  u(y') = \int_{\R^d} \Big( \nabla \nabla G_0(x'-y) \cdot (A(y') - A(y)) \nabla (\eta (u-\ext u_L))(y)\\
 + (u-\ext u_L)(y) \nabla \nabla G_0(x'-y) \cdot A(y) \nabla \eta(y) \\
  + \nabla G_0(x'-y) \big(\nabla \eta(y) \cdot A(y) \nabla (u-\ext u_L)(y)-\mu \ext u_L \eta(y) \Big) \;dy.
\end{multline*}
As above, this can be justified by mollification of the RHS of~\eqref{eq:eta_u^gamma}.
Indeed, the limit is well-defined since the constant-coefficient Green function $G_0$ classically satisfies 
\begin{eqnarray*}
 |\nabla G_0(y)|=|\nabla G_\mu(y,0;A(y'))(y)| &\le& C(d,\lambda) |y|^{1-d} 
 \\
 |\nabla \nabla G_0(y)|=|\nabla\nabla G_\mu(y,0;A(y'))| &\le& C(d,\lambda) |y|^{-d}
\end{eqnarray*}
uniformly in $y,y' \in \R^d$, while by assumption, the coefficient field satisfies $|A(y') - A(y)| \le \mathcal{C}_\gamma |y'-y|^\gamma$. It then follows
\begin{multline}\label{eq:rep_nabla_eta_u^gamma}
 |\nabla u(y')| \,\lesssim \, \mu |\ext u_L|+\int_{B_L(x)} |\nabla \nabla G_0(y'-y)| |A(y') - A(y)| |\nabla u(y)| \;dy
 \\
 + L^{-1} \int_{\mathcal{A}_{\frac{2L}{3},L}(x)} \Big( |u(y)-\ext u_L| |\nabla \nabla G_0(y'-y)| + |\nabla G_0(y'-y)| |\nabla u(y)| \big) \Big) \;dy
\end{multline}
for all $y' \in B_{\frac{L}{2}}(x)$, where $\mathcal{A}_{L',L''}(x):= \{ y : L'\le |y-x| \le L'' \}$ denotes the annulus centered at $x$ and of radii $L'$ and $L''$. Since $L\sim 1$, we allow the constant in $\lesssim$ to depend on $L$.
The constant-coefficient bounds yield
$$
 \int_{\mathcal{A}_{\frac{2L}{3},L}(x)} |u(y)-\ext u_L| |\nabla \nabla G_0(y'-y)| \;dy \lesssim \int_{B_{L}(x)} |u(y)-\ext u_L|\;dy.
$$
Combined with Jensen's and Poincar\'e's inequalities, this turns into
\begin{equation}\label{eq:Hol-1.1}
\int_{\mathcal{A}_{\frac{2L}{3},L}(x)} |u(y)-\ext u_L| |\nabla \nabla G_0(y'-y)| \;dy \,\lesssim\, \bigg( \int_{B_{L}(x)} |\nabla u(y)|^2 \;dy \bigg)^\frac{1}{2}.
\end{equation}
Likewise we obtain
\begin{equation}\label{eq:Hol-1.2}
 \int_{\mathcal{A}_{\frac{2L}{3},L}(x)} |\nabla G_0(y'-y)| |\nabla u(y)| \;dy \lesssim \bigg( \int_{B_{L}(x)} |\nabla u(y)|^2 \;dy \bigg)^\frac{1}{2}.
\end{equation}
We are left with the second RHS term  of~\eqref{eq:rep_nabla_eta_u^gamma}, which we bound, by the decay of $\nabla\nabla G_0$ and the H\"older continuity of $A$, by
\begin{equation}\label{eq:Hol-1.3}
 \int_{B_L(x)} |\nabla \nabla G_0(y'-y)| |A(y') - A(y)| |\nabla u(y)| \;dy \lesssim \int_{B_L(x)} |x'-y|^{\gamma-d} |\nabla u(y)| \;dy.
\end{equation}
Let $p\ge 2$.
We then take the $p$-th power of~\eqref{eq:rep_nabla_eta_u^gamma}, use \eqref{eq:Hol-1.1}--\eqref{eq:Hol-1.3}, and integrate over $y'$ in $B_{\frac{L}{2}}(x)$. This yields
\begin{multline}\label{eq:Hol-1.5}
 \int_{B_{\frac{L}{2}}(x)} |\nabla u(y')|^p \;dy' \lesssim \, ( \mu |\ext u_L|)^p+\int_{B_{\frac{L}{2}}(x)} \bigg( \int_{B_L(x)} |y'-y|^{\gamma-d} |\nabla u(y)| \;dy \bigg)^p \;dy' \\ + \bigg( \int_{B_{L}(x)} |\nabla u(y)|^2 \;dy \bigg)^\frac{p}{2}.
\end{multline}
We are almost in position to apply Young's convolution inequality. Mimicking its proof, we let $r$ and $p'$ be such that 
\begin{equation}\label{young_iteration}
p\ge p'\ge 1,p>r\ge 1, \quad 1+\frac{1}{p}=\frac{1}{r}+\frac{1}{p'},
\end{equation}
and use H\"older's inequality with exponents $(p,\frac{pp'}{p-p'},\frac{rp}{p-r})$ on the integrand
$$
|y'-y|^{\gamma-d} |\nabla u(y)| \,=\, \Big(|y'-y|^{(\gamma-d)\frac{r}{p}} |\nabla u(y)|^{\frac{p'}{p}} \Big)
\Big(|\nabla u(y)|^{\frac{p-p'}{p}} \Big)
\Big( |y'-y|^{(\gamma-d)\frac{p-r}{p}} \Big).
$$
This yields
\begin{multline}\label{eq:Hol-1.4}
 \bigg( \int_{B_L(x)} |y'-y|^{\gamma-d} |\nabla u(y)| \;dy \bigg)^p \le \int_{B_L(x)} |y'-y|^{(\gamma-d)r} |\nabla u(y)|^{p'} \;dy \\ \times \bigg( \int_{B_L(x)} |\nabla u(y)|^{p'} \;dy \bigg)^{\frac{p}{p'}-1} \bigg( \int_{B_L(x)} |y'-y|^{(\gamma-d)r} \;dy \bigg)^{\frac{p}{r}-1}.
\end{multline}
As long as we choose $1\le r<\frac{d}{d-\gamma}<2$ (since $\gamma<1$ and $d\ge 2$), the last RHS term is bounded (depending on $L$). Let us fix such an $1<r<2\le p$, in which case the exponents $(p,r,p'=\frac{p r}{r+(r-1) p})$ satisfy \eqref{young_iteration}. Integrating \eqref{eq:Hol-1.4} over $y'\in B_{\frac{L}{2}}(x)$ then yields
\[
 \int_{B_{\frac{L}{2}}(x)} \bigg( \int_{B_L(x)} |y'-y|^{\gamma-d} |\nabla u(y)| \;dy \bigg)^p \;dy' \lesssim \bigg( \int_{B_L(x)} |\nabla u(y)|^{p'} \;dy \bigg)^\frac{p}{p'}
\]
Combined with \eqref{eq:Hol-1.5}, this gives for all $p\ge 2$, $1<r<2$, and $p'=\frac{p r}{r+(r-1) p}$,
\begin{equation}\label{eq:Hol-1.6}
 \| \nabla u \|_{L^p(B_{\frac{L}{2}}(x))}\, \lesssim \, \mu  \|u\|_{L^\infty(B_L(x))}+\| \nabla u \|_{L^{p'}(B_L(x))} + \| \nabla u \|_{L^{2}(B_{L}(x))}.
\end{equation}
We start from $p_0'=2$ (that is, with $p_0=\frac{2 r}{r-(r-1) 2}>2$) and $L_0=\ell$, and iterate using the following exponents and ball size:
$$
p'_{n+1}:=p_n, \ p_{n+1}:=\frac{p_{n} r}{r-(r-1) p_n} , L_{n+1}=\frac{L_{n}}{2}.
$$
So defined, $p_n$ is a monotonically increasing sequence, so that $(p_n,r,p_n')$ satisfies \eqref{young_iteration} for all $n\in \N_0$ such that $p_n<\frac{r}{r-1}$. In particular,  \eqref{eq:Hol-1.6} then yields
\begin{equation*}
 \| \nabla u \|_{L^{p_n}(B_{\frac{\ell}{2^{n}}}(x))} \,\lesssim \, \mu  \|u\|_{L^\infty(B_\ell(x))}+ \| \nabla u \|_{L^{2}(B_{\ell}(x))}.
\end{equation*}
In addition, $p_n$ satisfies $p_n \ge (\frac{r}{r-2(r-1)})^n 2$, so that after \emph{finitely many
} steps, $p_n$ is such that $\frac{p_{n} r}{r-(r-1) p_n}>\frac{r}{r-1}$, at which point we may choose $p_{n+1} = \infty$.
This proves~\eqref{L^2_L^infty}, and Corollary~\ref{cor:C^gamma} now follows directly from Theorem~\ref{thm:main} with $L=\ell$, noting that the deterministic estimate on the Green function itself yields
\begin{equation*}
0\, \le \,\mu G_\mu(x,y) \, \lesssim \, \mu \ln(2+\frac{1}{\sqrt{\mu}|x-y|})\frac{e^{-c\sqrt{\mu}|x-y|}}{|x-y|^{d-2}}
\lesssim \, \frac{e^{-c\sqrt{\mu}|x-y|}}{|x-y|^{d-1}} ,
\end{equation*}
for any $0<c'<c$, as desired.


\section{Proofs of the fluctuation estimates}

\subsection{Proof of Theorem~\ref{thm:strong_error}}

We first assume that the coefficient field $A$ and the right-hand side $f$ are smooth. Since the estimates do not depend on the smoothness of the parameters, we may at the end lift this restriction by approximation.
The triangle inequality yields
$$
 \expec{ \bigg( \int_{\R^d} |u(x) - \langle u(x) \rangle|^p \;dx \bigg)^\theta }^{\frac{1}{p\theta}} \le \bigg( \int_{\R^d} \big\langle |u(x) - \langle u(x) \rangle|^{p\theta} \big\rangle^{\frac{1}{\theta}} \;dx \bigg)^{\frac{1}{p}}.
$$
Appealing to the spectral gap estimate of Lemma~\ref{lem:var-estim} with exponent $\frac{p\theta}{2}\ge 1$ yields
$$
 \bigg( \int_{\R^d} \big\langle |u(x) - \langle u(x) \rangle|^{p\theta} \big\rangle^{\frac{1}{\theta}} \;dx \bigg)^{\frac{1}{p}} \lesssim \bigg( \int_{\R^d} \expec{\bigg( \int_{\R^d} \Big( \osc{A|_{B_\ell(z)}}{u(x)} \Big)^2 \;dz \bigg)^{\frac{p\theta}{2}} }^{\frac{1}{\theta}} \;dx \bigg)^{\frac{1}{p}},
$$
and by the triangle inequality
$$
 \bigg( \int_{\R^d} \big\langle |u(x) - \langle u(x) \rangle|^{p\theta} \big\rangle^{\frac{1}{\theta}} \;dx \bigg)^{\frac{1}{p}} \lesssim \Bigg( \int_{\R^d} \bigg( \int_{\R^d} \expec{\Big( \osc{A|_{B_\ell(z)}}{u(x)} \Big)^{p\theta}}^{\frac{2}{p\theta}} \;dz \bigg)^{\frac{p}{2}} \;dx \Bigg)^{\frac{1}{p}}.
$$
By the oscillation estimate of Lemma~\ref{lem:osc_u}, this turns into
\begin{equation}\label{proof_strong1}
 \expec{\bigg( \int_{\R^d} |u(x) - \langle u(x) \rangle|^p \;dx \bigg)^\theta }^{\frac{1}{p\theta}} \lesssim \Bigg( \int_{\R^d} \bigg( \int_{\R^d} \Big\langle \mathcal{K}_{G_\mu,u}(x,z)^{p\theta} \Big\rangle^{\frac{2}{p\theta}} \;dz \bigg)^{\frac{p}{2}} \;dx \Bigg)^{\frac{1}{p}}.
\end{equation}
We now estimate the RHS. By the Cauchy-Schwarz inequality and Theorem~\ref{thm:main}, we have
\begin{equation}\label{proof_strong2}
 \expec{\mathcal{K}_{G_\mu,u}(x,z)^{p\theta} }^{\frac{2}{p\theta}}\, \le\, K(x-z)^2 \expec{ (\nabla u)_{9\ell}(z)^{2p\theta}}^{\frac{1}{p\theta}} + \chi_{B_{6\ell}}(x-z) \| f \|_{L^{\lambda}(B_{2\ell}(x))}^2,
\end{equation}
where again $\chi_D$ denotes the characteristic function of the set $D\subseteq \R^d$ and $K$ is the kernel
$$
 K(x-z) \,= \,\frac{e^{-c\sqrt{\mu}|x-z|}}{1+|x-z|^{d-1}} .
$$
In the following, the constant $c>0$ in $K$ may change from line to line (and only depends on
$\lambda$ and $d$).
In order to correctly capture the decay of $(\nabla u)_{9\ell}(z)$, we write $u$ in terms of its Green function representation and split the sum into two contributions:
$$
 u(z) = \int_{\R^d} G_\mu(z,y) f(y)\;dy = \int_{\R^d\setminus B_{11\ell}(z)} G_\mu(z,y) f(y)\;dy + \int_{B_{11\ell}(z)} G_\mu(z,y) f(y)\;dy.
$$
Thus
\begin{align}\label{eq:pr-strong-0}
 \expec{|(\nabla u)_{9\ell}(z)|^{2p\theta}}^{\frac{1}{2p\theta}} &= \, \expec{\bigg|\int_{B_{9\ell}(z)} \bigg( \int_{\R^d} \nabla_{z'} G_\mu(z',y) f(y) \;dy \bigg)^2 \;dz' \bigg|^{p\theta} }^{\frac{1}{2p\theta}}\\
 &\lesssim \, \expec{\bigg|\int_{B_{9\ell}(z)} \bigg( \int_{\R^d\setminus B_{11\ell}(z)} \nabla_{z'} G_\mu(z',y) f(y) \;dy \bigg)^2 \;dz' \bigg|^{p\theta} }^{\frac{1}{2p\theta}}\nonumber\\
&\qquad\qquad\qquad + \expec{\bigg|\int_{B_{9\ell}(z)} \bigg( \int_{B_{11\ell}(z)} \nabla_{z'} G_\mu(z',y) f(y) \;dy \bigg)^2 \;dz' \bigg|^{p\theta}}^{\frac{1}{2p\theta}}.\nonumber
\end{align}
We start by estimating the second RHS term, and consider the function
$$
 v:z'\mapsto \int_{B_{11\ell}(z)} G_\mu(z',y) f(y) \;dy,
$$
which solves on $\R^d$
$$
 \mu v-\nabla\cdot A\nabla v = f \chi_{B_{11\ell}(z)}.
$$
Set $\bar v:=\fint_{B_{11\ell}(z)} vdy$.
An energy estimate combined with the Sobolev embedding on $B_{11\ell}(z)$ yields for $\lambda > \frac{d}{2} \ge \frac{2d}{d+2}$,
\begin{multline}\label{eq:pr-strong-1-0}
 \|\nabla v\|^2_{L^2(\R^d)} \, \lesssim \,  \int_{B_{11\ell}(z)}  f vdy
\,\leq \, \int_{B_{11\ell}(z)}  f (v-\bar v)dy+\bar v \int_{B_{11\ell}(z)}  fdy
\\
\lesssim\,  \| f \|_{L^\lambda(B_{11\ell}(z))} \|\nabla v\|_{L^2(\R^d)}+|\bar v|\| f \|_{L^1(B_{11\ell}(z))}  .
\end{multline}
It remains to estimate $\bar v$. By the triangle inequality and H\"older's inequality
with exponents $(\lambda',\lambda)$, we have using the pointwise bounds \eqref{eq:ptwise-decay-estim}
on $G_\mu$ in Definition~\ref{def:Green}
\begin{multline}\label{eq:pr-strong-1-1}
|\bar v|\,\leq\, \int_{B_{11\ell}(z)}\int_{B_{11\ell}(z)} |G_\mu(z',y)||f(y)| \;dydz'
\\
\lesssim\,\int_{B_{11\ell}(z)}  \Big(\int_{B_{11\ell}(z)} G_\mu(z',y)^{\lambda'}dy\Big)^{\frac{1}{\lambda'}} 
\Big(\int_{B_{11\ell}(z)} |f|^\lambda dy\Big)^{\frac{1}{\lambda}} dz'
\\
\lesssim \, \kappa_d(\mu) \| f \|_{L^\lambda(B_{11\ell}(z))} ,
\end{multline}
where $\kappa_d(\mu)=1$ if $d>2$ and $\kappa_d(\mu)=|\ln \mu|+1$ if $d=2$, since $1\le \lambda'<\frac{d}{d-2}$.
By \eqref{eq:pr-strong-1-0}, \eqref{eq:pr-strong-1-1}, and Young's inequality, we may 
thus bound the second RHS of \eqref{eq:pr-strong-0} by
\begin{equation}\label{eq:pr-strong-1}
\int_{B_{9\ell}(z)} \bigg( \int_{B_{11\ell}(z)} \nabla_{z'} G_\mu(z',y) f(y) \;dy\bigg)^2dz' \,=\, \|\nabla v\|^2_{L^2(B_{9\ell}(z))} \, \lesssim \,  \kappa_d(\mu) \| f \|^2_{L^\lambda(B_{11\ell}(z))}.
\end{equation}
We then turn to the first RHS term of \eqref{eq:pr-strong-0},
and take local averages using H\"older's inequality with exponents $(\lambda',\lambda)$ (with respect to $dy$):
\begin{multline*}
 \expec{\Bigg( \int_{B_{9\ell}(z)} \bigg( \int_{\R^d\setminus B_{11\ell}(z)} \nabla_{z'} G(z',y) f(y) \;dy \bigg)^2 \;dz' \Bigg)^{p\theta} }^{\frac{1}{2p\theta}}\\ \lesssim \,\expec{\Bigg(\int_{B_{9\ell}(z)} \bigg( \int_{\R^d\setminus B_{11\ell}(z)} \| \nabla_{z'} G(z',y') \|_{L^{\lambda'}_{y'}(B_\ell(y))} \| f \|_{L^\lambda(B_\ell(y))} \;dy \bigg)^2 \;dz' \Bigg)^{p\theta} }^{\frac{1}{2p\theta}}.
\end{multline*}
Combined with the triangle inequality in $L^2_{z'}(B_{9\ell}(z))$, this yields
\begin{multline*}
 \expec{\Bigg( \int_{B_{9\ell}(z)} \bigg( \int_{\R^d\setminus B_{11\ell}(z)} \nabla_{z'} G(z',y) f(y) \;dy \bigg)^2 \;dz' \Bigg)^{p\theta} }^{\frac{1}{2p\theta}}\\ \lesssim \,\expec{ \bigg( \int_{\R^d\setminus B_{11\ell}(z)} \| \nabla_{z'} G(z',y') \|_{L^{\lambda'}_{y'}(B_\ell(y),L^2_{z'}(B_{9\ell}(z)))} \| f \|_{L^\lambda(B_\ell(y))} \;dy \bigg)^{2p\theta} }^{\frac{1}{2p\theta}}.
\end{multline*}
From the De Giorgi-Nash-Moser theory in the form of Lemma~\ref{lem:DGNM} (with RHS zero), we then have
\begin{equation*}
 \| \nabla_{z'} G(z',y') \|_{L^{\lambda'}_{y'}(B_\ell(y),L^2_{z'}(B_{9\ell}(z)))} \,\lesssim\,  (\nabla G)_{9\ell}(z,y).
\end{equation*}
We then finally appeal to Theorem~\ref{thm:main} and  the triangle inequality with respect to $L^{2p\theta}_{\langle\cdot\rangle}$ to obtain the following estimate of the first RHS term of \eqref{eq:pr-strong-0}:
\begin{multline}\label{proof_strong4}
 \expec{\Bigg( \int_{B_{9\ell}(z)} \bigg( \int_{\R^d\setminus B_{11\ell}(z)} \nabla_{z'} G(z',y) f(y) \;dy \bigg)^2 \;dz' \Bigg)^{p\theta} }^{\frac{1}{2p\theta}}\\ \lesssim \int_{\R^d\setminus B_{11\ell}(z)} K(z-y) \| f \|_{L^\lambda(B_\ell(y))} \;dy.
\end{multline}
Since $K(z-y)\sim 1$ for $y\in  B_{11\ell}(z)$, the combination of~\eqref{eq:pr-strong-0}, \eqref{eq:pr-strong-1}, and~\eqref{proof_strong4} yields
\begin{equation}\label{proof_strong5}
 \expec{|(\nabla u)_{9\ell}(z)|^{2p\theta} }^{\frac{1}{2p\theta}}\, \lesssim \, \int_{\R^d} K(z-y) \| f \|_{L^\lambda(B_{11\ell}(y))} \;dy.
\end{equation}
In total, collecting~\eqref{proof_strong1},~\eqref{proof_strong2} and~\eqref{proof_strong5},
we then have
\begin{multline*}
 \expec{\bigg( \int_{\R^d} |u(x) - \langle u(x) \rangle|^p \;dx \bigg)^\theta}^{\frac{1}{p\theta}} 
\,\lesssim \,
\kappa_d(\mu)^{\frac{1}{2}}\bigg( \int_{\R^d} \| f \|_{L^\lambda(B_{2\ell}(z))}^{p} \;dz \bigg)^{\frac{1}{p}}\\ + \Bigg( \int_{\R^d} \bigg( \int_{\R^d} K(x-z)^2 \Big( \int_{\R^d} K(z-y) \| f \|_{L^\lambda(B_{11\ell}(y))} \;dy \Big)^2 \;dz \bigg)^{\frac{p}{2}} \;dx \Bigg)^{\frac{1}{p}}.
\end{multline*}
Since $q\le p$ and the integral of the RHS term is equivalent to a discrete sum over an appropriate lattice of size $\ell$, we have that
$$
 \bigg( \int_{\R^d} \| f \|_{L^\lambda(B_{2\ell}(z))}^{p} \;dz \bigg)^{\frac{1}{p}} \lesssim \bigg( \int_{\R^d} \| f \|_{L^\lambda(B_{2\ell}(z))}^{q} \;dz \bigg)^{\frac{1}{q}}\,\lesssim\,\bigg( \int_{\R^d} \| f \|_{L^\lambda(B_{\ell}(z))}^{q} \;dz \bigg)^{\frac{1}{q}}.
$$
The most important term is the last one. By the triangle inequality in $L^2_y(\R^d)$,
\begin{multline*}
 \Bigg( \int_{\R^d} \bigg( \int_{\R^d} K(x-z)^2 \Big( \int_{\R^d} K(z-y) \| f \|_{L^\lambda(B_{6\ell}(y))} \;dy \Big)^2 \;dz \bigg)^{\frac{p}{2}} \;dx \Bigg)^{\frac{1}{p}}\\ \le \Bigg( \int_{\R^d} \bigg( \int_{\R^d} \Big( \int_{\R^d} K(x-z)^2 K(z-y)^2 \| f \|_{L^\lambda(B_{11\ell}(y))}^2 \;dz \Big)^{\frac{1}{2}} \;dy \bigg)^{p} \;dx \Bigg)^{\frac{1}{p}}.
\end{multline*}
We bound the integral over $z$ as follows:
\begin{multline}\label{int_kernel}
 \int_{\R^d} \frac{e^{-c\sqrt{\mu}|x-z|}}{1+|x-z|^{2(d-1)}} \frac{e^{-c\sqrt{\mu}|z-y|}}{1+|z-y|^{2(d-1)}} dz\,
\lesssim \,
\begin{cases} 
\dps \frac{e^{-c\sqrt{\mu}|x-y|}}{1+|x-y|^{2(d-1)}}  &\text{ if } d>2,\\
\dps (|\ln\mu|+1) \frac{e^{-c\sqrt{\mu}|x-z|}}{1+|x-z|^{2}} &\text{ if } d=2.
\end{cases}
\end{multline}
In other words,
$$
 \int_{\R^d} K(x-z)^2 K(z-y)^2 \;dz \lesssim K(x-y)^2 \kappa_d(\mu),
$$
where we recall that $\kappa_d(\mu)=1$ for $d>2$ and $\mu_d(\mu)=|\ln \mu|+1$ for $d=2$.
We thus have
\begin{multline*}
 \Bigg( \int_{\R^d} \bigg( \int_{\R^d} \Big( \int_{\R^d} K(x-z)^2 K(z-y)^2 \| f \|_{L^\lambda(B_{6\ell}(y))}^2 \;dz \Big)^{\frac{1}{2}} \;dy \bigg)^{\frac{p}{2}} \;dx \Bigg)^{\frac{1}{p}}\\ \lesssim 
\kappa_d(\mu)^\frac{1}{2}\Bigg( \int_{\R^d} \bigg( \int_{\R^d} K(x-y) \| f \|_{L^\lambda(B_{6\ell}(y))} \;dy \bigg)^{p} \;dx \Bigg)^{\frac{1}{p}}.
\end{multline*}
Let us pick $1\le r\le \frac{d}{d-1}$ and $1\le q < +\infty$ such that~\eqref{params_strong} holds.
If $r<\frac{d}{d-1}$, Young's inequality yields
\begin{equation}\label{Young_strong}
 \Bigg( \int_{\R^d} \bigg( \int_{\R^d} K(x-y) \| f \|_{L^\lambda(B_{6\ell}(y))} \;dy \bigg)^{p} \;dx \Bigg)^{\frac{1}{p}} \lesssim \| K \|_{L^r(\R^d)} \bigg( \int_{\R^d} \| f \|_{L^\lambda(B_{6\ell}(y))}^{q} \;dx \bigg)^{\frac{1}{q}}.
\end{equation}
We easily check that 
\begin{equation}\label{K_Young}
 \| K \|_{L^r(\R^d)} = \bigg( \int_{\R^d} K(x)^r \;dx \bigg)^{\frac{1}{r}} \lesssim 1+\mu^{-\frac{(1-d)r+d}{2r}}.
\end{equation}
In the border-line case $r=\frac{d}{d-1}$, the Hardy-Littlewood-Sobolev inequality immediately yields
provided $q>1$
\begin{equation}\label{HLS_strong}
 \Bigg( \int_{\R^d} \bigg( \int_{\R^d} |x-y|^{1-d} \| f \|_{L^\lambda(B_{11\ell}(y))} \;dy \bigg)^{p} \;dx \Bigg)^{\frac{1}{p}} \lesssim \bigg( \int_{\R^d} \| f \|_{L^\lambda(B_{11\ell}(y))}^{q} \;dx \Bigg)^{\frac{1}{q}},
\end{equation}
where we have also used the elementary fact that $\frac{1}{1+|x-y|^{d-1}} \le \frac{1}{|x-y|^{d-1}}$.
Collecting~\eqref{Young_strong},~\eqref{K_Young} and~\eqref{HLS_strong} yields
\begin{align*}
 \Bigg( \int_{\R^d} \bigg( \int_{\R^d} K(x-y) \| f \|_{L^\lambda(B_{11\ell}(y))} \;dy \bigg)^{p} \;dx \Bigg)^{\frac{1}{p}} &\lesssim \,(1+\mu^{-\frac{(1-d)r+d}{2r}}) \kappa_d(\mu)^\frac{1}{2}\bigg( \int_{\R^d} \| f \|_{L^\lambda(B_{6\ell}(x))}^{q} \;dx \Bigg)^{\frac{1}{q}}\\
 &\,\lesssim (1+\mu^{-\frac{(1-d)r+d}{2r}}) \kappa_d(\mu)^\frac{1}{2}
\bigg( \int_{\R^d} \| f \|_{L^\lambda(B_{\ell}(x))}^{q} \;dx \Bigg)^{\frac{1}{q}},
\end{align*}
where $p$, $q$ and $r$ are related by~\eqref{params_strong}.
This concludes the proof of the theorem.

\subsection{Proof of Theorem~\ref{thm:weak_error}}

Since transposition is a linear local operator, if $A$ satisfies the assumptions of Theorem~\ref{thm:weak_error}, then $A^*$ does as well, so that  the statement of Theorem~\ref{thm:weak_error} is symmetric with respect to interchanging $f$ and $g$ provided $A$ is replaced by $A^*$. Hence we may without loss of generality assume that $\lambda_1 \le \lambda_2$. By~\eqref{lambda_i}, this implies that
\begin{equation}\label{lambda_1}
 \lambda_2 \ge \frac{2d}{d+2}.
\end{equation}
By Jensen's inequality in probability we may assume w.~l.~o.~g. that $\theta\ge 2$.
The spectral gap estimate of Lemma~\ref{lem:var-estim}  for $q=\frac{\theta}{2}\ge 1$ yields
$$
 \expec{\Big| \int_{\R^d} (u(x) - \langle u(x) \rangle) g(x) \;dx \Big|^\theta}^{\frac{1}{\theta}} 
\,\lesssim\,  \expec{\bigg( \int_{\R^d} \Big( \osc{A|_{B_\ell(z)}}{\int_{\R^d} u(x) g(x) \;dx} \Big)^{2} \;dz \bigg)^{\frac{\theta}{2}} }^{\frac{1}{\theta}}.
$$
By the triangle inequality, we may insert the unperturbed solution $u$ and estimate
\begin{multline*}
\expec{\bigg( \int_{\R^d} \Big( \osc{A|_{B_\ell(z)}}{\int_{\R^d} u(x) g(x) \;dx} \Big)^{2} dz \bigg)^{\frac{\theta}{2}} }^{\frac{1}{\theta}} 
\\
\leq \,
2\expec{\bigg( \int_{\R^d} \Big( \sup_{\tilde A|_{B_\ell(z)}} \int_{\R^d} \big| u(x) - \tilde u(x) \big| |g(x)| \;dx \Big)^{2} dz \bigg)^{\frac{\theta}{2}} }^{\frac{1}{\theta}}.
\end{multline*}
Taking local averages combined with H\"older's inequality with exponents $(\lambda_1',\lambda_1)$ yields
\begin{multline*}
  \expec{\Bigg| \int_{\R^d} (u(x) - \langle u(x) \rangle) g(x) \;dx \Big|^\theta}^{\frac{1}{\theta}}\\
 \lesssim \,\expec{\bigg( \int_{\R^d} \Big( \sup_{\tilde A|_{B_\ell(z)}} \int_{\R^d} \| u - \tilde u \|_{L^{\lambda_1'}(B_\ell(x))} \| g \|_{L^{\lambda_1}(B_\ell(x))} \;dx \Big)^{2} \;dz \bigg)^{\frac{\theta}{2}} }^{\frac{1}{\theta}}.
\end{multline*}
We then put the supremum inside the inner integral and appeal to the sensitivity estimate of Lemma~\ref{lem:osc_u} to obtain
\begin{multline*}
\expec{\bigg( \int_{\R^d} \Big( \sup_{A|_{B_\ell(z)}} \int_{\R^d} \| u - \tilde u \|_{L^{\lambda_1'}(B_\ell(x))} \| g \|_{L^{\lambda_1}(B_\ell(x))} \;dx \Big)^{2} \;dz \bigg)^{\frac{\theta}{2}} }^{\frac{1}{\theta}}\\
 \lesssim \,\expec{\bigg( \int_{\R^d} \Big( \int_{\R^d} \mathcal{K}_{G,u}(x,z) \| g \|_{L^{\lambda_1}(B_\ell(x))} \;dx \Big)^{2} \;dz \bigg)^{\frac{\theta}{2}} }^{\frac{1}{\theta}}.
\end{multline*}
It remains to estimate the RHS.
By the triangle inequality in $L^s_{\langle\cdot\rangle}$, first with $s = \frac{\theta}{2}\ge 1$ and then $s=2$, we 
have
\begin{multline*}
 \expec{\bigg( \int_{\R^d} \Big( \int_{\R^d} \mathcal{K}_{G,u}(x,z) \| g \|_{L^{\lambda_1}(B_\ell(x))} \;dx \Big)^{2} \;dz \bigg)^{\frac{\theta}{2}} }^{\frac{1}{\theta}}\\ 
\le \, \bigg( \int_{\R^d} \Big( \int_{\R^d} \expec{ \mathcal{K}_{G,u}(x,z)^\theta }^{\frac{1}{\theta}} \| g \|_{L^{\lambda_1}(B_\ell(x))} \;dx \Big)^{2} \;dz \bigg)^{\frac{1}{2}}.
\end{multline*}
We then make use of~\eqref{proof_strong2} in the proof of Theorem~\ref{thm:strong_error} with $\lambda = \lambda_2$:
\begin{multline*}
 \bigg( \int_{\R^d} \Big( \int_{\R^d} \expec{\mathcal{K}_{G,u}(x,z)^\theta }^{\frac{1}{\theta}} \| g \|_{L^{\lambda_1}(B_\ell(x))} \;dx \Big)^{2} \;dz \bigg)^{\frac{1}{2}}
\\ 
\lesssim \, \bigg( \int_{\R^d} \Big( \int_{\R^d} K(x-z) \expec{ (\nabla u)_{9\ell}(z)^{2\theta} }^{\frac{1}{2\theta}} \| g \|_{L^{\lambda_1}(B_\ell(x))} \;dx \bigg)^{2} \;dz \Bigg)^{\frac{1}{2}}
\\
+ \bigg( \int_{\R^d} \| f \|_{L^{\lambda_2}(B_{2\ell}(z))}^2 \| g \|_{L^{\lambda_1}(B_{7\ell}(z))}^2 \;dz \bigg)^{\frac{1}{2}}.
\end{multline*}
By H\"older's inequality with $\frac{1}{2} = \frac{1}{q_1} + \frac{1}{\tilde q_1}$, we bound the second RHS term by 
$$
 \bigg( \int_{\R^d} \| f \|_{L^{\lambda_2}(B_{2\ell}(z))}^2 \| g \|_{L^{\lambda_1}(B_{7\ell}(z))}^2 \;dz \bigg)^{\frac{1}{2}} \lesssim \| f \|_{L^{q_1}_{\lambda_2,1}(\R^d)} \| g \|_{L^{\tilde q_1}_{\lambda_1,1}(\R^d)}.
$$
By~\eqref{params_weak}, since $r,\tilde r\ge 1$, we may choose $q_1 \ge q$ and $\tilde q_1 \ge \tilde q$ so that
$$
\| f \|_{L^{q_1}_{\lambda_2,1}(\R^d)} \| g \|_{L^{\tilde q_1}_{\lambda_1,1}(\R^d)}\, \lesssim\, \| f \|_{L^{q}_{\lambda_2,1}(\R^d)} \| g \|_{L^{\tilde q}_{\lambda_1,1}(\R^d)}.
$$
From~\eqref{proof_strong5} (with $p=1$) in the proof of Theorem~\ref{thm:strong_error}, we learn that
\begin{multline*}
 \bigg( \int_{\R^d} \Big( \int_{\R^d} K(x-z) \expec{ (\nabla u)_{2\ell}(z)^{2\theta} }^{\frac{1}{2\theta}} \| g \|_{L^{\lambda_1}(B_\ell(x))} \;dx \Big)^{2} \;dz \bigg)^{\frac{1}{2}}\\
 \lesssim \, \bigg( \int_{\R^d} \Big( \int_{\R^d} \int_{\R^d} K(x-z) K(z-y) \| f \|_{L^{\lambda_2}(B_{11\ell}(y))} \| g \|_{L^{\lambda_1}(B_\ell(x))} \;dxdy \Big)^{2} \;dz \bigg)^{\frac{1}{2}},
\end{multline*}
which holds by our choice $\lambda_2\ge \lambda_1$ which implies $\lambda_2> \frac{2d}{d+2}$ by \eqref{params_weak2}.
Let $p, \tilde p \ge 1$ be two exponents to be specified later such that $\frac{1}{2} = \frac{1}{p} + \frac{1}{\tilde p}$. We then have that
\begin{multline*}
 \bigg( \int_{\R^d} \Big( \int_{\R^d} \int_{\R^d} K(x-z) K(z-y) \| f \|_{L^{\lambda_2}(B_{11\ell}(y))} \| g \|_{L^{\lambda_1}(B_\ell(x))} \;dxdy \Big)^{2} \;dz \big)^{\frac{1}{2}}
\\
 \lesssim \,\bigg( \int_{\R^d} \Big( \int_{\R^d} K(z-y) \| f \|_{L^{\lambda_2}(B_{11\ell}(y))} \;dy \Big)^{p} \;dz \bigg)^{\frac{1}{p}} \bigg( \int_{\R^d} \Big( \int_{\R^d} K(x-z) \| g \|_{L^{\lambda_1}(B_\ell(x))} \;dx \Big)^{\tilde p} \;dz \bigg)^{\frac{1}{\tilde p}}.
\end{multline*}
We treat the two factors of the RHS the same way. 
First we consider the non-borderline case $r<\frac{d}{d-1}$, in which case Young's convolution inequality with $1 + \frac{1}{p} = \frac{1}{r} + \frac{1}{q}$ yields
$$
 \bigg( \int_{\R^d} \Big( \int_{\R^d} K(z-y) \| f \|_{L^\lambda_2(B_{11\ell}(y))} \;dy \Big)^{p} \;dz \bigg)^{\frac{1}{p}} \,\lesssim\, \| K \|_{L^r(\R^d)} \| f \|_{L^{q}_{\lambda_2,1}(\R^d)} \,\lesssim\, \mu^{-\frac{1-d}{2} - \frac{d}{2r}} \| f \|_{L^{q}_{\lambda_2,1}(\R^d)}.
$$
In the borderline case  $r=\frac{d}{d-1}$, the result follows from the Hardy-Littlewood-Sobolev inequality
provided $q>1$. An identical estimate holds for the second factor with exponents $1 + \frac{1}{\tilde p} = \frac{1}{\tilde r} + \frac{1}{\tilde q}$ (provided $\tilde q>1$ in the borderline case). Gathering these two estimates yields
\begin{multline*}
 \bigg( \int_{\R^d} \Big( \int_{\R^d} \int_{\R^d} K(x-z) K(z-y) \| f \|_{L^{\lambda_2}(B_{11\ell}(y))} \| g \|_{L^{\lambda_1}(B_\ell(x))} \;dxdy \Big)^{2} \;dz \bigg)^{\frac{1}{2}}
\\ 
\lesssim \, \mu^{-(1-d) - \frac{d}{2} ( \frac{1}{r} + \frac{1}{\tilde r})} \| f \|_{L^{q}_{\lambda_2,1}(\R^d)} \| g\|_{L^{\tilde q}_{\lambda_1,1}(\R^d)},
\end{multline*}
with
$$
 2 + \frac{1}{2} = 1 + \frac{1}{p} + 1 + \frac{1}{\tilde p} = \frac{1}{r} + \frac{1}{\tilde r} + \frac{1}{q} + \frac{1}{\tilde q}.
$$
This completes the proof.


\section{Proof of the Lipschitz regularity theory}

\subsection{Proof of Theorem~\ref{thm:main2}}

As opposed to the corresponding proof in the discrete case, cf. \cite[Corollary 4]{Marahrens-Otto-13}, we have to take care of the singularity of the Green function. This prevents us to make use of Morrey's inequality when the coefficients are only measurable, and we propose a more direct approach which partly mimics the proof of Morrey's inequality. We assume w.~l.~o.~g. that $R>9L$.
In the first five steps we assume that $d>2$, and indicate the changes for $d=2$ in Step~6.
\medskip

\step{1} Representation formula for $u(x+x')-u(x')$, $x\in B_R\setminus B_{2L}$, $x'\in B_L$.

\noindent 
In order to make use of the annealed estimates of Theorem~\ref{thm:main}, we rewrite equation \eqref{eq:elliptic-u} on $\R^d$ as follows. Let $\eta:\R^d\to [0,1]$ be a cutoff-function for $B_{\frac{4R}{3}}$ in $B_{\frac{5R}{3}}$ such that $|\nabla \eta|\lesssim R^{-1}$. A direct calculation shows that $\eta u \in H^1(\R^d)$ satisfies
\begin{equation}\label{eta_u}
\mu \eta u - \nabla \cdot A \nabla (u\eta) = \mu \eta u - \eta \nabla\cdot A \nabla u - \nabla \eta \cdot A \nabla u - \nabla \cdot ( u A \nabla \eta ).
\end{equation}
The sum of the first two RHS terms equals $\eta f$ while the other two terms belong to $H^{-1}(\R^d)$ and have compact support. The Green reprensentation formula yields
$$
 (\eta u)(x) = \int_{B_{2R}} \Big( G_\mu(x,y) \big( \eta(y) f(y) - \nabla \eta(y) \cdot A(y) \nabla u(y) \big) +  u(y) \nabla_y G_\mu(x,y) \cdot A(y) \nabla \eta(y)  \Big) dy.
$$
Assume first that $f$ and $A$ are smooth (so that $u$ is smooth and the formula holds classically).
We argue by density.
Since $0\le G_\mu(x,y) \le |y-x|^{2-d}$, $\eta f\in L^p(\R^d)$ with $p>\frac{d}{2}$,  
$\nabla \eta=0$ on $B_{3R/2}$ (and in particular at the singularity of $G_\mu(x,\cdot)$), and $\nabla u \in L^2(B_{2R})$, the first term of the integral is well-defined at the limit. Recalling that $y\mapsto \nabla_y G_\mu(x,y)$ is locally square-integrable away from $y=x$, the second term of the integral is well-defined as well since $\nabla \eta$ vanishes in a neighborhood of the singularity of $y\mapsto \nabla_y G_\mu(x,y)$ and $u\in L^2(B_{2R})$. 
Since $u$ is uniformly H\"older continuous, one can also take the limit of the LHS, so that the Green representation formula holds by a density and regularization argument.

We thus have for all $x\in B_R\setminus B_{2L}$ and $x'\in B_L$, using in addition that $\nabla \eta$ vanishes on $B_{\frac{3R}{2}}$, 
\begin{multline}\label{eq:main-step1}
u(x+x')-u(x')\,=\,\int_{B_{2R}} (G_\mu(x+x',y)-G_\mu(x',y)) \eta(y) f(y) dy
\\-\int_{B_{2R}\setminus B_{\frac{4R}{3}}} \Big((G_\mu(x+x',y)-G_\mu(x',y))\nabla \eta(y) \cdot A(y) \nabla u(y) \\
-  u(y) \nabla_y (G_\mu(x+x',y)-G_\mu(x',y)) \cdot A(y) \nabla \eta(y)  \Big) dy.
\end{multline}

\medskip

\step{2} Estimate of the integral on $B_{3L}(x)\cup B_{3L}\subset B_{2L}$.

\noindent Since $R>9L$, $x\in B_{R}$ and $x'\in B_L$, $B_{3L}(x)\cup B_{3L}\subset B_{\frac{4R}{3}}$, only the first integral term of the RHS of \eqref{eq:main-step1} has a contribution. We shall argue that
\begin{equation}\label{eq:main-step2}
\Big|\int_{B_{3L}(x)\cup B_{3L}} (G_\mu(x+x',y)-G_\mu(x',y)) \eta(y) f(y) dy\Big| \,\lesssim \,
\bigg(\int_{B_{2R}} |f|^q dy\bigg)^{\frac{1}{q}}.
\end{equation}
Indeed, the deterministic pointwise estimates on $G_\mu$ for $d>2$ combined with H\"older's inequality with exponents $(\frac{q}{q-1},q)$ yield
\begin{multline*}
\Big|\int_{B_{3L}(x)\cup B_{3L}} (G_\mu(x+x',y)-G_\mu(x',y)) \eta(y) f(y) dy\Big| 
\\
\lesssim \, \bigg( \int_{B_{3L}(x)\cup B_{3L}}(|x+x'-y|^{\frac{q(2-d)}{q-1}}+|x'-y|^{\frac{q(2-d)}{q-1}})dy\bigg)^{\frac{q-1}{q}} 
\times \bigg(\int_{B_{2R}} |f|^q dy\bigg)^{\frac{1}{q}}.
\end{multline*}
Since $q>\frac{d}{2}$ implies $\frac{q(d-2)}{q-1}<d$, the first factor is of order $1$, and 
\eqref{eq:main-step2} follows.

\medskip

\step{3} Representation formulas for $G_\mu(x+x',y)-G_\mu(x',y)$ and $\nabla_y G_\mu(x+x',y)-\nabla_y G_\mu(x',y)$, $x\in B_R\setminus B_{2L}$, $x'\in B_L$, $y\notin B_{3L}(x)\cup B_{3L}$.

\noindent When $y$ is not at the singularity of the Green function, we may write the difference of Green functions as the directional integral of its gradient: for all $y\notin [x',x'+x]$,
\begin{equation}\label{eq:main-step3-1}
G_\mu(x+x',y)-G_\mu(x',y)
\,=\, \int_0^1 \nabla_xG_\mu(tx+x',y)\cdot xdt,
\end{equation}
and for all $i\in \{1,\dots,d\}$,
\begin{equation}\label{eq:main-step3-2}
\nabla_{y_i} G_\mu(x+x',y)-\nabla_{y_i} G_\mu(x',y)
\,=\, \int_0^1 \nabla_x \nabla_{y_i} G_\mu(tx+x',y)\cdot xdt.
\end{equation}
When $y$ is close to $[x',x'+x]$, we have to refine this decomposition.
To this end, we define two points $x^+$ and $x^-$ and two sets $B^+$ and $B^-$ as follows:
$$
x^+:=\frac{x}{2}+(\frac{|x|}{2}+L)e_1,\quad x^-:=\frac{x}{2}-(\frac{|x|}{2}+L)e_1,
$$
where $e_1$ is the first unit vector of the canonical basis of $\R^d$, and 
\begin{eqnarray*}
B^+&:=&\{y\in B_{2R}\setminus (B_{3L}(x)\cup B_{3L}), (y-x)\cdot e_1\le 0\},\\
B^-&:=&\{y\in B_{2R}\setminus (B_{3L}(x)\cup B_{3L}), (y-x)\cdot e_1>0\}.
\end{eqnarray*}
Note that $B^+\cup B^-=B_{2R}\setminus (B_{3L}(x)\cup B_{3L})$.
For $x\in B^+$ we write $G_\mu(x+x',y)-G_\mu(x',y)=G_\mu(x+x',y)-G_\mu(x^++x',y)+G_\mu(x^++x',y)-G_\mu(x',y)$, so that 
\begin{multline}\label{eq:main-step3-3}
G_\mu(x+x',y)-G_\mu(x',y)\\
\,=\, \int_0^1 \nabla_x G_\mu(x^++t(x-x^+)+x',y)\cdot(x-x^+)dt+\int_0^1 \nabla_xG_\mu(tx^++x',y)\cdot x^+dt.
\end{multline}
We proceed correspondingly for $y\in B^-$.

\medskip

In the following step we estimate the RHS of \eqref{eq:main-step1}.
In view of Step~2, it only remains to estimate the integrals on $\tilde B:=B_{2R}\setminus (B_{3L}(x)\cup B_{3L})$. 

\medskip

\step{4} Estimates of the integrals on $\tilde B$.

\noindent We shall prove three estimates.
First,
\begin{multline}\label{eq:main-step4}
\Big|\int_{B^\pm} (G_\mu(x+x',y)-G_\mu(x',y)) \eta(y) f(y) dy\Big| \, 
\\
\lesssim \, |x|\bigg( \int_{B_{2R}} |f|^q dy\bigg)^{\frac{1}{q}}
\int_0^1 \bigg(\int_{B^\pm}\Big(|\nabla_x G_\mu(x^\pm+t(x-x^\pm)+x',y)|^{\frac{q}{q-1}}\\
+|\nabla_xG_\mu(tx^\pm+x',y)|^{\frac{q}{q-1}} \Big)dy\bigg)^{\frac{q-1}{q}} dt,
\end{multline}
where $B^\pm$ is a shorthand notation we use when the inequality holds both on $B^+$ and $B^-$. We only prove the claim for $B^+$. Since $|x|>L$, by construction $|x-x^+|\lesssim |x|$ and $|x^+|\lesssim |x|$, so that
\eqref{eq:main-step4} follows from \eqref{eq:main-step3-3} and H\"older's inequality with exponents $(\frac{q}{q-1},q)$.

\medskip

The second estimate is:
\begin{multline}\label{eq:main-step5}
\Big|\int_{\tilde B\setminus B_{\frac{4R}{3}}} (G_\mu(x+x',y)-G_\mu(x',y))\nabla \eta(y) \cdot A(y) \nabla u(y) dy \Big| \, 
\\
\lesssim \, |x|\bigg( R^{-1} \Big( \int_{B_{2R}} u^2 dy\Big)^{\frac{1}{2}}
+\Big( \int_{B_{2R}} f^2dy\Big)^{\frac{1}{2}} \bigg)\\
\times
\int_0^1 \bigg(\int_{\tilde B\setminus B_{\frac{4R}{3}}}|\nabla_x G_\mu(tx+x',y)|^{2}dy\bigg)^{\frac{1}{2}} dt.
\end{multline}
We proceed as for the proof of \eqref{eq:main-step4} and use in addition the following consequence of the definition of $\eta$ and Caccioppoli's inequality:
$$
\bigg(\int_{B_{2R}} |\nabla \eta|^2|\nabla u|^2dy\bigg)^{\frac{1}{2}} \, \lesssim \,R^{-1}\bigg(\int_{B_{2R}} u^2dy\bigg)^{\frac{1}{2}}+\bigg(\int_{B_{2R}} f^2dy\bigg)^{\frac{1}{2}}.
$$
Indeed, since $\nabla \eta$ has support in $B_{\frac{5R}{3}}\setminus B_{\frac{4R}{3}}$,
$$
\int_{B_{2R}} |\nabla \eta|^2|\nabla u|^2 dy\, \lesssim \,R^{-2}\int_{B_{\frac{5R}{3}}\setminus B_{\frac{4R}{3}}}
|\nabla u|^2dy. 
$$
Testing equation~\eqref{eq:elliptic-u} with test-function $\tilde \eta^2 u\in H^1_0(B_{2R})$, where $\tilde \eta$ has support in $B_{2R}$ and is such that $\tilde \eta|_{B_{\frac{5R}{3}}}\equiv 1$ and $|\nabla \tilde \eta|\lesssim \frac{1}{R}$,
 yields the Caccioppoli estimate
$$
\int_{B_{2R}} \nabla (\tilde \eta u)\cdot A \nabla (\tilde\eta u)dy\,\le \, \int_{B_{2R}} u^2 \nabla \tilde\eta\cdot A\nabla \tilde\eta dy+\int_{B_{2R}} \tilde\eta^2fu dy,
$$
which, by definition of $\tilde\eta$ and Young's inequality on the last term, we may use in the form
\begin{equation}\label{eq:Caccio-useful}
\int_{B_{\frac{5R}{3}}} |\nabla u|^2 dy\,\lesssim\, R^{-2} \int_{B_{2R}} u^2dy+R^2\int_{B_{2R}} f^2dy.
\end{equation}
Finally, we prove 
\begin{multline}\label{eq:main-step6}
\Big|\int_{\tilde B\setminus B_{\frac{4R}{3}}} u(y) \nabla_y (G_\mu(x+x',y)-G_\mu(x',y)) \cdot A(y) \nabla \eta(y)  dy \Big| \, 
\\
\lesssim \, |x|R^{-1} \bigg( \int_{B_{2R}} u^2 dy\bigg)^{\frac{1}{2}}
\int_0^1 \bigg(\int_{\tilde B\setminus B_{\frac{4R}{3}}}|\nabla\nabla G_\mu(tx+x',y)|^{2} 
dy\bigg)^{\frac{1}{2}} dt.
\end{multline}
This estimate follows from \eqref{eq:main-step3-2}, the bound $|\nabla \eta|\lesssim R^{-1}$, and Cauchy-Schwarz' inequality.

\medskip

\step{5} Conclusion for $d>2$.

\noindent The combination of \eqref{eq:main-step1}, \eqref{eq:main-step2}, \eqref{eq:main-step4}, \eqref{eq:main-step5}, and \eqref{eq:main-step6} yields, using that $|x|>L\sim 1$
and that $\tilde B\setminus B_{\frac{4R}{3}}=B_{2R}\setminus  B_{\frac{4R}{3}}$,
\begin{multline*}
\frac{|u(x+x')-u(x')|}{|x|}\,\lesssim \, R^{-1}\bigg(  \Big(\fint_{B_{2R}} |R^2f|^qdy \Big)^{\frac{1}{q}}+ \Big(\fint_{B_{2R}} (R^2f)^2dy \Big)^{\frac{1}{2}}+\Big(\fint_{B_{2R}} u^2dy \Big)^{\frac{1}{2}}\bigg)\\
\times\bigg\{ R^{-1+\frac{d}{q}}\int_0^1 \bigg(\int_{B^+}\Big(|\nabla_x G_\mu(x^++t(x-x^+)+x',y)|^{\frac{q}{q-1}}
+|\nabla_xG_\mu(tx^++x',y)|^{\frac{q}{q-1}} \Big)dy\bigg)^{\frac{q-1}{q}} dt \\
+R^{-1+\frac{d}{q}}\int_0^1 \bigg(\int_{B^-}\Big(|\nabla_x G_\mu(x^-+t(x-x^-)+x',y)|^{\frac{q}{q-1}}
+|\nabla_xG_\mu(tx^-+x',y)|^{\frac{q}{q-1}} \Big)dy\bigg)^{\frac{q-1}{q}} dt \\
 +R^{-1+\frac{d}{q}}+R^{-1+\frac{d}{2}}\int_0^1 \bigg(\int_{B_{2R}\setminus B_{\frac{4R}{3}}}|\nabla_x G_\mu(tx+x',y)|^{2}dy\bigg)^{\frac{1}{2}} dt\\
+R^{\frac{d}{2}}\int_0^1 \bigg(\int_{B_{2R}\setminus B_{\frac{4R}{3}}}|\nabla\nabla G_\mu(tx+x',y)|^{2}dy\bigg)^{\frac{1}{2}} dt
\bigg\}.
\end{multline*}
Dividing both sides of the inequality by the first RHS term and averaging over $x'\in B_L$ yield using Jensen's inequality and that $q>d$ (so that $R^{-1+\frac{d}{q}}\lesssim 1$):
\begin{multline*}
\frac{R \fint_{B_L}\frac{|u(x+x')-u(x')|}{|x|}dx' }{\Big(\fint_{B_{2R}} u^2 dy\Big)^{\frac{1}{2}}+ \Big(\fint_{B_{2R}} |R^2f|^qdy \Big)^{\frac{1}{q}}}\\
\lesssim \, 
\int_0^1 \bigg(\fint_{B_L}\int_{B^+}\Big(|\nabla_x G_\mu(x^++t(x-x^+)+x',y)|^{\frac{q}{q-1}}
+|\nabla_xG_\mu(tx^++x',y)|^{\frac{q}{q-1}} \Big)dydx'\bigg)^{\frac{q-1}{q}} dt \\
+\int_0^1 \bigg(\fint_{B_L}\int_{B^-}\Big(|\nabla_x G_\mu(x^-+t(x-x^-)+x',y)|^{\frac{q}{q-1}}
+|\nabla_xG_\mu(tx^-+x',y)|^{\frac{q}{q-1}} \Big)dydx'\bigg)^{\frac{q-1}{q}} dt \\
 +R^{-1+\frac{d}{2}}\int_0^1 \bigg(\fint_{B_L}\int_{B_{2R}\setminus B_{\frac{4R}{3}}}|\nabla_x G_\mu(tx+x',y)|^{2}dydx'\bigg)^{\frac{1}{2}} dt\\
+R^{\frac{d}{2}}\int_0^1 \bigg(\fint_{B_L}\int_{B_{2R}\setminus B_{\frac{4R}{3}}}|\nabla\nabla G_\mu(tx+x',y)|^{2}dydx'\bigg)^{\frac{1}{2}} dt+1
\,=:\,\calY_R(x).
\end{multline*}
This proves estimate \eqref{eq:thm-de_giorgi1}. It remains to prove
the moment bounds \eqref{eq:Y} on $\calY_R(x)$, which formally follow from taking the expectation of the $p$-th power of the RHS of this inequality and bounding
$|\nabla G_\mu(x,y)|$ by $|x-y|^{1-d}$ and $|\nabla \nabla G_\mu(x,y)|$ by $|x-y|^{-d}$. It remains to show that it is enough to use bounds on large moments of local square averages
of $|\nabla G_\mu(x,y)|$ and  $|\nabla \nabla G_\mu(x,y)|$ instead, which we control optimally by Theorem~\ref{thm:main}. We only treat the first term in detail (the other terms are treated similarly). By bounding the integral on $B^+$ by the sum of integrals on balls of radius $L$ and by H\"older's inequality, we have
\begin{multline*}
\fint_{B_L}\int_{B^+}\Big(|\nabla_x G_\mu(x^++t(x-x^+)+x',y)|^{\frac{q}{q-1}}
+|\nabla_xG_\mu(tx^++x',y)|^{\frac{q}{q-1}} \Big)dydx'
\\
\lesssim \, \sum_{i\in B^+\cap \frac{L}{\sqrt{d}}\Z^d} ((\nabla_x G_\mu)_L(x^++t(x-x^+),i))^{\frac{q}{q-1}}
+((\nabla_x G_\mu)_L(tx^+,i))^{\frac{q}{q-1}}.
\end{multline*}
We only treat the first RHS term.
By Jensen's inequality in probability it is enough to prove the claim for $p$ large enough, which we take 
such that $p\ge \frac{q}{q-1}$.
By Jensen's inequality on $\int_0^1dt$ and by the triangle inequality for $\expec{\int_0^1(\cdot)^{\frac{p(q-1)}{q}}}^{\frac{q}{p(q-1)}}$,
\begin{multline*}
\expec{\bigg(\int_0^1 \Big( \sum_{i\in B^+\cap \frac{L}{\sqrt{d}}\Z^d} ((\nabla_x G_\mu)_L(x^++t(x-x^+),i))^{\frac{q}{q-1}} \Big)^{\frac{q-1}{q}}dt\bigg)^p}
\\
\le \, \expec{\int_0^1 \Big( \sum_{i\in B^+\cap \frac{L}{\sqrt{d}}\Z^d} ((\nabla_x G_\mu)_L(x^++t(x-x^+),i))^{\frac{q}{q-1}} \Big)^{\frac{p(q-1)}{q}}dt}\\
\le \,
\bigg(\sum_{i\in B^+\cap \frac{L}{\sqrt{d}}\Z^d}  \expec{\int_0^1((\nabla_x G_\mu)_L(x^++t(x-x^+),i))^p dt}^{\frac{q}{p(q-1)}}\bigg)^{\frac{p(q-1)}{q}}.
\end{multline*}
Recall that by construction of $x^+$ and $B^+$, $|x^++t(x-x^+)-i|\sim |x-i|$ for all $t\in [0,1]$, so that 
by Theorem~\ref{thm:main}, 
$$
 \expec{\int_0^1((\nabla_x G_\mu)_L(x^++t(x-x^+),i))^p dt}^{\frac{1}{p}}\,\lesssim\, \frac{e^{-c\sqrt{\mu}|x-i|}}{|x-i|^{d-1}}.
$$
Giving up the exponential cut-off, this yields
\begin{multline*}
\expec{\bigg(\int_0^1 \Big( \sum_{i\in B^+\cap \frac{L}{\sqrt{d}}\Z^d} ((\nabla_x G_\mu)_L(x^++t(x-x^+),i))^{\frac{q}{q-1}} \Big)^{\frac{q-1}{q}}dt\bigg)^p}
\\
\lesssim \,
\bigg(\sum_{i\in B^+\cap \frac{L}{\sqrt{d}}\Z^d} |x-i|^{(1-d)\frac{q}{q-1}}\bigg)^{\frac{p(q-1)}{q}}
\,\lesssim\, \bigg(\int_{B_{2R}\setminus B_{L}(x)} |x-y|^{(1-d)\frac{q}{q-1}}dy\bigg)^{\frac{p(q-1)}{q}}
\,\lesssim \, 1
\end{multline*}
since $q>d$. This completes the proof of \eqref{eq:thm-de_giorgi1}.

\medskip

\step{6} Proof for $d=2$.

\noindent The proof for $d=2$ is identical as for $d>2$ except for Step~2. Indeed, if we proceed there as for $d>2$, the estimate fails optimality by a logarithm of $\mu$ due to the bound on the Green function $G_\mu$ in dimension 2 close to the singularity. Recall that $p>\frac{d}{2}=1$. 
To avoid this logarithmic correction, we follow the elegant argument by Avellaneda and Lin \cite{Avellaneda-Lin-87} and add a third dimension.
We denote by $G_\mu^{(2)}$ and $A^{(2)}$ the fields in dimension $2$ and consider the following extensions to dimension 3: $A^{(3)}(x_1,x_2,x_3):=\dig{A^{(2)}(x_1,x_2),1}$
and $G_\mu^{(3)}$ the Green function associated with $A^{(3)}$.
It is elementary to check using Definition~\ref{def:Green} that for all $x\neq y \in \R^2$,
\begin{equation*}
G_\mu^{(2)}(x,y)\,=\,\int_\R G_\mu^{(3)}((x,0),(y,t))dt,
\end{equation*}
and we rewrite the LHS of \eqref{eq:main-step2} as
\begin{multline}\label{eq:pro-d=2-sub0}
\int_{B_{3L}(x)\cup B_{3L}} (G_\mu^{(2)}(x+x',y)-G_\mu^{(2)}(x',y)) \eta(y) f(y) dy
\\=\, \int_{B_{3L}(x)\cup B_{3L}} \int_\R(G_\mu^{(3)}((x+x',0),(y,t))-G_\mu^{(3)}((x',0),(y,t))) \eta(y) f(y) dtdy.
\end{multline}
We then split the integral over $t$ into two parts: $|t|\le 1$ and $|t|>1$. We start by estimating the first part, and 
appeal to the deterministic pointwise estimate on $G_\mu^{(3)}$. By the triangle inequality,
\begin{multline*}
\Big|\int_{B_{3L}(x)\cup B_{3L}} \int_{|t|\le 1}(G_\mu^{(3)}((x+x',0),(y,t))-G_\mu^{(3)}((x',0),(y,t))) \eta(y) f(y) dtdy\Big|\\
\lesssim \,\int_{B_{3L}(x)\cup B_{3L}} \int_{|t|\le 1} \Big((|x+x'-y|^2+t^2)^{-\frac{1}{2}}+(|x'-y|^2+t^2)^{-\frac{1}{2}}\Big) |f(y)| dtdy
\end{multline*}
We first integrate in $y$ and use H\"older's inequality with exponents $(\frac{q}{q-1},q)$ for some $1<q\le p$ small enough so that $\frac{q}{q-1}>2$. This yields 
\begin{multline}\label{eq:pro-d=2-sub1}
\Big|\int_{B_{3L}(x)\cup B_{3L}} \int_{|t|\le 1}(G_\mu^{(3)}((x+x',0),(y,t))-G_\mu^{(3)}((x',0),(y,t))) \eta(y) f(y) dtdy\Big|\\
\lesssim \, \int_{|t|\le 1} |t|^{2\frac{q-1}{q}-1} dt \left(\int_{B_{3L}(x)\cup B_{3}}|f(y)|^qdy\right)^{\frac{1}{q}}
\,\lesssim \, \left(\int_{B_{3L}(x)\cup B_{3L}}|f(y)|^pdy\right)^{\frac{1}{p}},
\end{multline}
by Jensen's inequality since $L\sim 1$.

We turn to the second part of the integral. We bound the difference of the Green functions by the oscillation,
and appeal to the De Giorgi-Nash-Moser theory in the form of the deterministic estimate: For all $|t|>1$, and all $z,y\in \R^d$,
\begin{equation*}
\osc{z\in B_{2R}}{G_\mu^{(3)}((z,0),(y,t))} \,\lesssim \, |t|^{-{1+\alpha_0}} ,
\end{equation*}
for some $\alpha_0>0$ depending only on $\lambda$ (see \eqref{eq:BMO-type} in Step~2 of the proof of Lemma~\ref{lem:quenched} for details).
Since $x,x+x' \in B_{2R}$, this yields
\begin{multline}\label{eq:pro-d=2-sub2}
\Big|\int_{B_{3L}(x)\cup B_{3L}} \int_{|t|> 1}(G_\mu^{(3)}((x+x',0),(y,t))-G_\mu^{(3)}((x',0),(y,t))) \eta(y) f(y) dtdy\Big|\\
\le \,\int_{B_{3L}(x)\cup B_{3L}} \int_{|t|> 1}  \Big(\osc{z\in B_{2R}}{G_\mu^{(3)}((z,0),(y,t))} \Big)|f(y)| dtdy 
\\
\lesssim \, \int_{B_{3L}(x)\cup B_{3L}} \int_{|t|> 1} |t|^{-{1+\alpha_0}} |f(y)| dtdy \,\lesssim \, \int_{B_{3L}(x)\cup B_{3L}} |f(y)| dy.
\end{multline}
The desired estimate \eqref{eq:main-step2} for $d=2$ and $p>1$ follows from \eqref{eq:pro-d=2-sub0}, \eqref{eq:pro-d=2-sub1}, and \eqref{eq:pro-d=2-sub2}.

\subsection{Proof of Remark~\ref{cor:de_giorgi}}

Estimate \eqref{eq:thm-de_giorgi1} for all $x\in B_R$ is a straightforward combination of \eqref{eq:thm-de_giorgi1} for all $x\in B_R\setminus B_{2\ell}$ and of Schauder interior estimates.
We closely follow the corresponding proof in the discrete setting, cf. \cite[Corollary 4]{Marahrens-Otto-13}.

\medskip

\step{1}  Representation formula for solutions $u\in H^1(B_R)$ of 
$$
\mu u - \nabla\cdot A \nabla u = f \in L^p(B_{2R})
$$
in $B_{2R}$ for some $p>d$. 
Let $\eta$ be a smooth cutoff function for $B_{\frac{4R}{3}}$ in $B_{\frac{5R}{3}}$ such that $|\nabla \eta|\lesssim R^{-1}$.
We claim that for all $x\in B_{\frac{R}{2}}$,
\begin{equation}\label{repr_eta_u}
 \nabla u(x) = \int_{B_{2R}} \Big( \nabla_x G(x,y) \big( \eta(y) f(y) - \nabla \eta(y) \cdot A(y) \nabla u(y) \big) +  u(y) \nabla \nabla G(x,y) \cdot A(y) \nabla \eta(y)  \Big) dy
\end{equation}
Indeed the Leibniz rule yields
\begin{equation}\label{eta_u}
\mu \eta u - \nabla \cdot A \nabla (u\eta) = \mu \eta u - \eta \nabla\cdot A \nabla u - \nabla \eta \cdot A \nabla u - \nabla \cdot ( u A \nabla \eta ).
\end{equation}
The sum of the first two RHS terms equals $\eta f$ while the other two terms are in $H^{-1}(\R^d)$ and have compact support. Hence testing~\eqref{eta_u} with $G_\mu$ yields
\begin{equation}\label{repr_eta_u0}
 (\eta u)(x) = \int_{B_{2R}} \Big( G_\mu(x,y) \big( \eta(y) f(y) - \nabla \eta(y) \cdot A(y) \nabla u(y) \big) +  u(y) \nabla_y G_\mu(x,y) \cdot A(y) \nabla \eta(y)  \Big) dy
\end{equation}
and~\eqref{repr_eta_u} follows by taking the derivative w.~r.~t.~$x$. Note that the RHS of \eqref{repr_eta_u0} and \eqref{repr_eta_u} are well-defined for all $x\in B_{\frac{R}{2}}$ (so tha Green representation formula follows from mollifying the RHS). On the one hand, $G_\mu\in L^{\frac{d}{d-2+\e}}(B_{2R})$ and $\nabla G_\mu \in L^{\frac{d}{d-1+\e}}(B_{2R})$ for all $\e>0$ and $f\in L^p(B_R)$ for some $p>d$, so that the terms involving $f$ are well-defined. 
On the other hand, $\nabla G_\mu$, $G_\mu$, and $\nabla \nabla G_\mu$ are locally square-integrable away from the singularity,  and $\nabla \eta$ vanishes in $B_{\frac{4R}{3}}$ so that the terms involving $\nabla G_\mu$ or $\nabla \nabla G_\mu$ and $u$ or $\nabla u$ in \eqref{repr_eta_u0} and \eqref{repr_eta_u} are not singular and are integrable.

\medskip

\step{2} Proof that for $\alpha=1-\frac{d}{p}$,
\begin{multline}\label{eq:morrey}
 \Big( R^{\alpha} [u]_{C^\alpha(B_{R})} \Big)^p \lesssim R^{\alpha p} \int_{B_{R}}  \bigg( \int_{B_{\frac{5R}{3}}} |\nabla_x G_\mu(x,y)| |f(y)| \;dy \bigg)^p\;dx\\ + R^{p(\alpha-1)} \int_{B_{R}}  \bigg( \int_{\mathcal A_{\frac{4R}{3},\frac{5R}{3}}} \big( |\nabla \nabla G_\mu(x,y)| |u(y)| + |\nabla_x G_\mu(x,y)| |\nabla u(y)| \big) \;dy \bigg)^p \;dx.
\end{multline}
Indeed, in view of the definition of $\eta$, \eqref{repr_eta_u} in Step~1 yields for all $x\in B_{\frac{R}{2}}$
\begin{multline*}
| \nabla  u(x)|\, \lesssim \, \int_{B_{\frac{5R}{3}}} |\nabla_x G_\mu(x,y)| |f(y)|dy
\\+ R^{-1}\int_{\mathcal A_{\frac{4R}{3},\frac{5R}{3}}} \Big(|\nabla_x G_\mu(x,y)|| \nabla u(y)|+|u(y)||\nabla \nabla G_\mu(x,y)|\Big) dy,
\end{multline*}
where $\mathcal A_{\frac{4R}{3},\frac{5R}{3}}=\{\frac{4R}{3}<|y|\le \frac{5R}{3}\}$.
The desired estimate \eqref{eq:morrey} then follows from Morrey's inequality
$$
 [u]_{C^\alpha(B_{R})} = \sup_{\substack{x,y \in B_{R}\\ x\neq y}} \frac{|u(x)-u(y)|}{|x-y|^\alpha} \lesssim \bigg( \int_{B_{R}} |\nabla u|^p \;dy \bigg)^{\frac{1}{p}}$$
and the triangle inequality.

\medskip

\step{3} Proof of
\begin{multline}\label{eq:gamma_nabla_u}
  \expec{ \Bigg( \sup_{(u,f)} 
\frac{ R^\alpha \sup_{x, y \in B_R}{\frac{|u(x)-u(y)|}{|x-y|^\alpha}}}
{ \sup_{B_{2R}} |u| + (\fint_{B_{2R}} |R^2 f|^p )^{\frac{1}{p}} }  
\Bigg)^p }\\
 \lesssim 
  \bigg\langle R^{d(p-2)-p} \int_{B_{R}} \int_{\mathcal A_{\frac{4R}{3},\frac{5R}{3}}} |\nabla_x G_\mu(x,y)|^p \;dydx +  R^{d(p-2)} \int_{B_{R}} \int_{\mathcal A_{\frac{4R}{3},\frac{5R}{3}}} |\nabla \nabla G_\mu(x,y)|^p \;dydx\\
 + R^{-2p} \int_{B_{R}} \bigg( \int_{B_{2R}} |\nabla_x G_\mu(x,y)|^{\frac{p}{p-1}} \;dy \bigg)^{p-1} \;dx \bigg\rangle.
\end{multline}
The starting point is \eqref{eq:morrey} in Step~2, and we treat each of the three RHS terms separately.
For the first term, we use H\"older's inequality with exponents $(\frac{p}{p-1},p)$:
\begin{eqnarray*}
\lefteqn{ \int_{B_{R}} \bigg( \int_{B_\frac{5R}{3}} |\nabla_x G_\mu(x,y)| |f(y)| \;dy \bigg)^p \;dx}
\\ &\le& 
 R^{-2p} 
 \int_{B_{R}} \bigg( \int_{B_{2R}} |\nabla_x G_\mu(x,y)|^{\frac{p}{p-1}} \;dy \bigg)^{p-1} \;dx \, \int_{B_{2R}} |R^2f|^p\;dy .
\end{eqnarray*}
For the second term, we also use H\"older's inequality with exponents $(p,\frac{p}{p-1})$:
\begin{eqnarray*}
\lefteqn{ \int_{B_{R}} \Bigg( \int_{\mathcal A_{\frac{4R}{3},\frac{5R}{3}}} |\nabla \nabla G_\mu(x,y)| |u(y)| \;dy \Bigg)^p \;dx }
\\
&\le& \int_{B_{R}} \int_{\mathcal A_{\frac{4R}{3},\frac{5R}{3}}} |\nabla \nabla G_\mu(x,y)|^p \;dydx \Bigg( \int_{\mathcal A_{\frac{4R}{3},\frac{5R}{3}}} |u|^{\frac{p}{p-1}} \;dy \Bigg)^{p-1}
\\
 &\lesssim& R^{d(p-2)}\int_{B_{R}} \int_{\mathcal A_{\frac{4R}{3},\frac{5R}{3}}} |\nabla \nabla G_\mu(x,y)|^p \;dydx\,  R^d\Big( \sup_{B_{2R}} |u| \Big)^p.\nonumber
\end{eqnarray*}
Likewise, for the third term we have 
\begin{multline*}
 \int_{B_{R}} \Bigg( \int_{\mathcal A_{\frac{4R}{3},\frac{5R}{3}}} |\nabla_x G_\mu(x,y)| |\nabla u(y)| \;dy \Bigg)^p \;dx \le \int_{B_{R}} \int_{\mathcal A_{\frac{4R}{3},\frac{5R}{3}}} |\nabla_x G_\mu(x,y)|^p \;dydx \\ \times \Bigg( \int_{B_R} |\nabla u|^{\frac{p}{p-1}} \;dy \Bigg)^{p-1}.
\end{multline*}
Since $p>d\ge 2$, we have ${\frac{p}{p-1}}<2$, so that Jensen's inequality yields
\[
 \Bigg( \int_{B_R} |\nabla u|^{\frac{p}{p-1}} \;dy \Bigg)^{p-1} \lesssim R^{d(\frac{p}{2}-1)} \Bigg( \int_{B_R} |\nabla u|^2 \;dy \Bigg)^{\frac{p}{2}}.
\]
By Caccioppoli's inequality (cf. \eqref{eq:Caccio-useful}), 
\begin{equation*}
 \int_{B_R} |\nabla u|^2 \;dy\, \lesssim \,R^{-2} \int_{B_{2R}} |u|^2 \;dy+R^2 \int_{B_{2R}}f^2\;dy
 \\
 \lesssim \, R^{d-2}  \sup_{B_{2R}} |u|^2+R^2 \int_{B_{2R}}f^2\;dy;
\end{equation*}
and consequently, by Jensen's inequality on $f$ (using that $p>d\ge 2$),
\begin{equation*}
 \Bigg( \int_{B_R} |\nabla u|^{\frac{p}{p-1}} \;dy \Bigg)^{p-1} \lesssim R^{d(p-2)-p}  \bigg(R^d\sup_{ B_{2R}} |u|^p+ \int_{B_{2R}} |R^2f|^p\;dy\bigg).
\end{equation*}
Hence we have proved the following bound for the third RHS term of \eqref{eq:morrey}: 
\begin{multline*}
  \int_{B_{R}} \Bigg( \int_{\mathcal A_{\frac{4R}{3},\frac{5R}{3}}} |\nabla_x G_\mu(x,y)| |\nabla u(y)| \;dy \Bigg)^p \;dx
  \\ \lesssim \, R^{d(p-2)-p} \bigg( R^d\sup_{B_{2R}} |u|^p+\int_{B_{2R}} |R^2f|^p\;dy\bigg)\int_{B_{R}} \int_{\mathcal A_{\frac{4R}{3},\frac{5R}{3}}} |\nabla_x G_\mu(x,y)|^p \;dydx .
\end{multline*}
This concludes the proof of \eqref{eq:gamma_nabla_u} recalling that $R^{\alpha p}=R^{d}$.

\medskip

\step{4} Conclusion.

\noindent We bound each term of the r.~h.~s. of \eqref{eq:gamma_nabla_u} separately.
The first term is bounded by
\begin{equation*}
 R^{d(p-2)-p} \int_{B_{R}} \int_{\mathcal A_{\frac{4R}{3},\frac{5R}{3}}} \expec{ |\nabla_x G(x,y)|^p } \;dydx \,\lesssim \,R^{d(p-2)-p} \int_{B_{R}} \int_{\mathcal A_{\frac{4R}{3},\frac{5R}{3}}} |x-y|^{(1-d)p} \;dydx.
\end{equation*}
For $x\in B_{R}$ and $y\in\mathcal A_{\frac{4R}{3},\frac{5R}{3}}$, we have that $|x-y|\ge |y|-|x| \ge \frac{R}{3}$, so that
\begin{equation*}
 R^{d(p-2)-p} \int_{B_{R}} \int_{\mathcal A_{\frac{4R}{3},\frac{5R}{3}}} |x-y|^{(1-d)p} \;dydx \,\lesssim \,R^{d(p-2)-p+2d + (1-d)p} = 1.
\end{equation*}
Likewise, the second term is bounded by
\begin{align*}
 R^{d(p-2)} \int_{B_{R}} \int_{\mathcal A_{\frac{4R}{3},\frac{5R}{3}}} \expec{ |\nabla\nabla G_\mu(x,y)|^p} \;dydx &\lesssim \, R^{d(p-2)} \int_{B_{R}} \int_{\mathcal A_{\frac{4R}{3},\frac{5R}{3}}} |x-y|^{-dp} \;dydx\\
&\lesssim \,R^{d(p-2)+2d - dp} = 1.
\end{align*}
For the third term, we use the triangle inequality in form of
\begin{align*}
 \expec{ \bigg( \int_{B_{2R}} |\nabla_x G_\mu(x,y)|^\frac{p}{p-1} \;dy \bigg)^{p-1} } &\le\, \bigg( \int_{B_{2R}} \expec{ |\nabla_x G_\mu(x,y)|^p }^{\frac{1}{p-1}} \;dy \bigg)^{p-1}\\
 &\lesssim \,\bigg(\int_{B_{2R}} |x-y|^{(1-d)\frac{p}{p-1}} \;dy \bigg)^{p-1}
 \lesssim R^{d(p-1)+(1-d)p}\,=\,R^{p-d}.
\end{align*}
Hence,
\begin{equation*}
 \expec{ R^{-2p} \int_{B_{R}} \bigg( \int_{B_{2R}} |\nabla_x G(x,y)|^{\frac{p}{p-1}} \;dy \bigg)^{p-1} \;dx} \,\lesssim\, R^{-2p+d+p-d} = R^{-p}\,\lesssim \,1.
\end{equation*}
As before the bound on $\calY_R(x)$ 
is a simple reformulation. The proof of the remark is complete.


\section{Proofs of the auxiliary results}

\subsection{Proof of Lemma~\ref{lem:LSI_p}}

The proof is essentially identical to the proof in the discrete case. The only difference lies in the different form of the (LSI). We reproduce the proof for completeness. 

\medskip

\step{1} Result for $p=1$. 

\noindent We claim that for any $\delta > 0$ and all $\zeta(a)$:
\begin{equation}\label{L1.7}
 \langle \zeta^2 \rangle^{\frac{1}{2}} \le \bigg( \exp\Big(\frac{2}{\rho\delta^2}\Big) + \frac{\rho 
\delta^2}{2e}\bigg) \langle |\zeta| \rangle + \delta
\Big\langle \int_{\R^d} \Big(\oscDaniel_{A|_{B_\ell(z)}} \zeta \Big)^2 \;dz \Big\rangle^{\frac{1}{2}},
\end{equation}
where $\rho$ denotes the constant in the (LSI), see Definition~\ref{def:LSI}.
By homogeneity, we may assume $\langle\zeta^2\rangle = 1$.
For all real-valued $\zeta$ we have that
\begin{equation*}
 \zeta^2 \le \begin{cases}
            \exp(\frac{2}{\rho\delta^2}) |\zeta|& \quad\text{if $|\zeta| \le \exp{\frac{2}{\rho\delta^2}}$}\\
	    {\textstyle\frac{\rho\delta^2}{4}} \zeta^2 \log \zeta^2 &\quad\text{if
$|\zeta|\ge\exp{\frac{2}{\rho\delta^2}}$}
           \end{cases}\Bigg\}.
\end{equation*}
Since $x\log x$ is bounded from below by $\frac{1}{e}$, we have that $\frac{2}{e} |\zeta| + \zeta^2 \log \zeta^2 \ge 0$
for all $\zeta$. It follows that
\begin{equation*}
 \zeta^2 \le \bigg( \exp\Big(\frac{2}{\rho\delta^2}\Big) + \frac{\rho \delta^2}{2e} \bigg) |\zeta| +
\frac{\rho\delta^2}{4} \zeta^2 \log \zeta^2.
\end{equation*}
Hence taking the expectation $\langle \cdot \rangle$ yields
\begin{equation*}
 \langle \zeta^2 \rangle \le \bigg( \exp\Big(\frac{2}{\rho\delta^2}\Big) + \frac{\rho 
\delta^2}{2e} \bigg) \langle |\zeta| \rangle +
\frac{\rho\delta^2}{4}
\Big\langle \zeta^2 \log \zeta^2 \Big\rangle.
\end{equation*}
Since $\langle\zeta^2\rangle = 1$, Young's inequality yields
\begin{align*}
 \langle |\zeta| \rangle &\le \frac{1}{2} \bigg( \exp\Big(\frac{2}{\rho\delta^2}\Big) + \frac{\rho 
\delta^2}{2e} \bigg) \langle
|\zeta| \rangle^2 + \frac{1}{2} \bigg( \exp\Big(\frac{2}{\rho\delta^2}\Big) + \frac{\rho 
\delta^2}{2e} \bigg)^{-1}\\
 &= \frac{1}{2} \bigg( \exp\Big(\frac{2}{\rho\delta^2}\Big) + \frac{\rho 
\delta^2}{2e} \bigg) \langle |\zeta| \rangle^2 + \frac{1}{2}
\bigg( \exp\Big(\frac{2}{\rho\delta^2}\Big) + \frac{\rho 
\delta^2}{2e} \bigg)^{-1} \langle \zeta^2 \rangle.
\end{align*}
Combining the last two estimates, we deduce
\begin{equation*}
 \langle \zeta^2 \rangle \le \bigg( \exp\Big(\frac{2}{\rho\delta^2}\Big) + \frac{\rho 
\delta^2}{2e} \bigg)^2 \langle |\zeta| \rangle^2
+ \frac{\rho\delta^2}{2} \Big\langle \zeta^2 \log \frac{\zeta^2}{\langle\zeta^2\rangle} \Big\rangle.
\end{equation*}
Hence (LSI) yields
\begin{equation*}
 \langle \zeta^2 \rangle \le \bigg( \exp\Big(\frac{2}{\rho\delta^2}\Big) + \frac{\rho 
\delta^2}{2e} \bigg)^2 \langle |\zeta| \rangle^2
+ \delta^2 \Big\langle \int_{\R^d} \Big(\oscDaniel_{A|_{B_\ell(z)}} \zeta \Big)^2 \;dz \Big\rangle
\end{equation*}
and estimate \eqref{L1.7} follows from taking the square root and applying the inequality $\sqrt{\zeta+\xi} \le \sqrt{\zeta} + \sqrt{\xi}$ for all numbers $\zeta, \xi \ge 0$.

\medskip

\step{2} We finish the proof of \eqref{LSI_p}, i.e.\ we show that
\begin{equation*}
\langle \zeta^{2p}\rangle^\frac{1}{2p}\le C(\rho,\ell,p,\delta)\langle|\zeta|\rangle
+\delta\Bigg(\Big\langle\bigg( \int_{\R^d} \Big(\oscDaniel_{A|_{B_\ell(z)}} \zeta \Big)^2 \;dz \bigg)^p\Big\rangle\Bigg)^\frac{1}{2p}
\end{equation*}
for general $p \ge 1$.
To that end, we apply \eqref{L1.7} to $\zeta$ replaced by $|\zeta|^p$: 
\begin{equation}\nonumber
\langle |\zeta|^{2p}\rangle\le C(\rho,p,\delta)\langle |\zeta|^{p}\rangle^2
+\delta\Big\langle \int_{\R^d} \Big(\oscDaniel_{A|_{B_\ell(z)}} |\zeta|^p \Big)^2 \;dz \Big\rangle,
\end{equation}
where $C(\rho,p,\delta)$ denotes a generic constant only depending on $\rho$, $p$, and $\delta$.
Since $p<2p$, an application of H\"older's inequality in $\langle\cdot\rangle$
and Young's inequality on the first RHS term yields
\begin{equation}\label{L1.8}
\langle |\zeta|^{2p}\rangle\le C(\rho,p,\delta)\langle|\zeta|\rangle^{2p}
+2\delta\Big\langle \int_{\R^d} \Big(\oscDaniel_{A|_{B_\ell(z)}} |\zeta|^p \Big)^2 \;dz \Big\rangle.
\end{equation}
Now we use that
\begin{equation*}
 \oscDaniel_{A|_{B_\ell(z)}} |\zeta|^p \le C(p) \bigg( |\zeta|^{p-1} \oscDaniel_{A|_{B_\ell(z)}} \zeta + \Big( \oscDaniel_{A|_{B_\ell(z)}} \zeta \Big)^p \bigg)
\end{equation*}
which follows from the elementary inequality $|\zeta^p - \xi^p| \le C(p) (\zeta^{p-1} |\zeta-\xi| + |\zeta-\xi|^p)$ for all numbers $\zeta,\xi > 0$ and the triangle inequality in form of $\oscDaniel_{a(e)} |\zeta| \le \oscDaniel_{a(e)} \zeta$. Hence \eqref{L1.8} yields
\begin{equation}\label{L1.osc^p}
\langle |\zeta|^{2p}\rangle\le C(\rho,p,\delta)\langle|\zeta|\rangle^{2p} + 2 C(p) \delta\Big\langle |\zeta|^{2p-2} \int_{\R^d} \Big(\oscDaniel_{A|_{B_\ell(z)}}
\zeta \;dz \Big)^{2} \Big\rangle\\
+2 C(p) \delta\Big\langle \int_{\R^d} \Big(\oscDaniel_{A|_{B_\ell(z)}} \zeta \Big)^{2p} \;dz \Big\rangle.
\end{equation}
The last term on the right-hand side may be estimated by discreteness, using the argument developed in \cite[Proof of Lemma~2.3]{Gloria-Otto-10b}. Since every ball $B_{\ell}(z)$, $z \in \R^d$ is contained in the collection $(B_{2\ell}(z'))_{z'\in\frac{2\ell}{\sqrt{d}}\Z^d}$, we have that 
\begin{equation*}
 \Big\langle \int_{\R^d} \Big(\oscDaniel_{A|_{B_\ell(z)}} \zeta \Big)^{2p} \;dz \Big\rangle \le C \Big\langle \sum_{z\in\frac{2\ell}{\sqrt{d}}\Z^d} \Big(\oscDaniel_{A|_{B_{2\ell}(z)}} \zeta \Big)^{2p} \Big\rangle.
\end{equation*}
Hence, by discreteness, we find have
\begin{equation}\label{L1.discrete}
 \Big\langle \int_{\R^d} \Big(\oscDaniel_{A|_{B_\ell(z)}} \zeta \Big)^{2p} \;dz \Big\rangle \le C \Big\langle \bigg( \sum_{z\in\frac{2\ell}{\sqrt{d}}\Z^d} \Big(\oscDaniel_{A|_{B_{2\ell}(z)}} \zeta \Big)^2 \bigg)^p \Big\rangle.
\end{equation}
Furthermore, H\"older's inequality followed by Young's inequality yields
\begin{align}\nonumber
 \Big\langle |\zeta|^{2p-2} \int_{\R^d} \Big(\oscDaniel_{A|_{B_\ell(z)}} \zeta \Big)^{2p} \;dz \Big\rangle &\le \langle |\zeta|^{2p} \rangle^{1-\frac{1}{p}} \Big\langle
\bigg( \int_{\R^d} \Big(\oscDaniel_{A|_{B_\ell(z)}} \zeta \Big)^{2p} \;dz \bigg)^p \Big\rangle\\
 &\le \frac{1}{4 C(p) \delta} \langle |\zeta|^{2p} \rangle + (4C(p)\delta)^{p-1} \Big\langle \bigg( \int_{\R^d} \Big(\oscDaniel_{A|_{B_\ell(z)}} \zeta \Big)^{2p} \;dz \bigg)^p \Big\rangle.\label{L1.Young}
\end{align}
Hence collecting \eqref{L1.osc^p}, \eqref{L1.discrete} and \eqref{L1.Young} yields
\begin{equation*}
 \langle |\zeta|^{2p}\rangle\le C(\rho,p,\delta)\langle|\zeta|\rangle^{2p} + 2\big(2 C(p) \delta + (4 C(p) \delta)^p \big) \Big\langle
\bigg( \sum_{z\in\frac{2\ell}{\sqrt{d}}} \Big(\oscDaniel_{A|_{B_{2\ell}(z)}} \zeta \Big)^2 \;dz
\bigg)^p \Big\rangle,
\end{equation*}
where we have absorbed the second term of \eqref{L1.Young} in the LHS.
Since every ball $B_{2\ell}(z')$, $z' \in \frac{2\ell}{\sqrt{d}}\Z^d$ is contained in the collection $(B_{3\ell}(z))_{|z-z'|\le2ell}$, we also deduce
\begin{equation*}
 \sum_{z\in\frac{2\ell}{\sqrt{d}}\Z^d} \Big(\oscDaniel_{A|_{B_{2\ell}(z)}} \zeta \Big)^2 \le \frac{1}{\ell} \int_{\R^d} \Big(\oscDaniel_{A|_{B_{3\ell}(z)}} \zeta \Big)^2 \;dz
\end{equation*}
By redefining $\delta$, we obtain \eqref{LSI_p}.

\subsection{Proof of Lemma~\ref{lem:quenched}}\label{sec:quenched}

Estimate \eqref{quenched_ap1} is a Meyers' type estimate, for which we refer the reader to \cite[Lemma~2.9]{Gloria-Otto-10b}.
We split the rest of the proof into four steps.
For $d>2$, \eqref{quenched_ap} is a consequence of \eqref{quenched_ap1} and of Meyers' estimate,
see Step~1.
For $d=2$, however, we need sharper deterministic estimates on the decay of local averages of the gradient of the Green function. These are obtained using the De Giorgi-Nash-Moser theory and pointwise bounds on the Green function in Step~2.
We then prove \eqref{quenched_ap} for all $d\ge 2$ in Step~3.
We prove \eqref{loc_bound} in the fourth and last step. 

\medskip 

\step{1} Proof of
\begin{equation}
 \int_{R\le|x-y|<2R} \int_{|y|<L} |\nabla \nabla G_\mu(x,y)|^{2q_1} \;dydx \,\lesssim \, R^{-2q_1\alpha_1} e^{-c\sqrt{\mu} R} ,\label{quenched_ap_1}
\end{equation}
for some $q_1>1$ and $\alpha_1>0$ and all $R\ge 4L\sim 1$.

\noindent This follows from Meyers' estimate in the form of: There exists some $q_0>1$ depending only on $\lambda$ and $d$ such
that for all $1\le q \le q_0$ and all functions $u\in H^1(\R^d),g\in L^2(\R^d,\R^d),f\in L^2(\R^d)$ supported in $B_{2L}$ with $L\sim 1$ and related through
$$
-\nabla \cdot A \nabla u\,=\,\nabla \cdot g+f,
$$
we have 
$$
\Big(\int_{\R^d} |\nabla u|^{2q_0}dy\Big)^{\frac{1}{2q}}\,\lesssim\, \Big(\int_{\R^d} |g|^{2q}dy\Big)^{\frac{1}{2q}}+\Big(\int_{\R^d} f^{2}dy\Big)^{\frac{1}{2}}.
$$
For this estimate we refer the reader to the original article by Meyers \cite{Meyers-63} or to 
\cite[(4.31) in Proof of Lemma~2.9]{Gloria-Otto-09} (the proof of which is first presented in the continuum setting dealt with here).

\medskip

Let $\eta:\R^d \to \R$ be such that $\eta\equiv 1$ on $B_L$, $\eta\equiv 0$ on $\R^d\setminus B_{2L}$ and $|\nabla \eta|\lesssim 1$.
Assume momentarily that $A$ is smooth, so that $(x,y)\mapsto G_\mu(x,y)$ is smooth away from the diagonal $x=y$.
Let $i\in \{1,\dots,d\}$, we apply Meyers' estimate to the smooth function $u(y)=\eta(y)\nabla_{x_i} G_\mu(y,x)$ 
for $|x|\ge 4L$.
Indeed, the defining equation for $G_\mu$ yields
$$
-\nabla \cdot A(y) \nabla u(y)\,=\, -\mu \eta(y)\nabla_{x_i} G_\mu(y,x)-\nabla_y \eta(y)\cdot A(y)\nabla_y \nabla_{x_i} G_\mu(y,x)
-\nabla \cdot (A(y)\nabla \eta(y)\nabla_{x_i} G_\mu(y,x)),
$$
so that Meyers' estimate with exponent $q_0>1$ takes the form 
\begin{multline*}
\Big(\int_{B_L} |\nabla_y\nabla_{x_i} G_\mu(y,x)|^{2q_0}dy\Big)^{\frac{1}{2q_0}}
\\ \lesssim\, \Big(\int_{B_{2L}} |\nabla_{x_i} G_\mu(y,x)|^{2q_0}dy\Big)^{\frac{1}{2q_0}}+\Big(\int_{B_{2L}} \mu^2(\nabla_{x_i} G_\mu(y,x))^{2}+|\nabla_y \nabla_{x_i} G_\mu(y,x)|^2dy\Big)^{\frac{1}{2}}.
\end{multline*}
By Caccioppoli's inequality (cf. \eqref{eq:Caccio-useful} in the proof of Theorem~\ref{thm:main2}), since $L\sim 1$,
$$
\int_{B_{2L}}|\nabla_y \nabla_{x_i} G_\mu(y,x)|^2dy\,\lesssim \,\int_{B_{3L}} (\nabla_{x_i} G_\mu(y,x))^2dy,
$$
so that by H\"older's inequality,
\begin{equation*}
\Big(\int_{B_L} |\nabla_y\nabla_{x_i} G_\mu(y,x)|^{2q_0}dy\Big)^{\frac{1}{2q_0}}
\, \lesssim\, (1+\mu)\Big(\int_{B_{3L}} |\nabla_{x_i} G_\mu(y,x)|^{2q_0}dy\Big)^{\frac{1}{2q_0}}.
\end{equation*}
Taking the $(2q_0)^{th}$ power of this inequality, summing over $i=1,\dots,d$, and integrating over $\{R\le |x|<2R\}$ yield
combined with \eqref{quenched_ap1} and $L\sim 1$
\begin{equation}\label{eq:quenched_a_2}
\int_{R\le |x|<2R} \int_{|y|< L} |\nabla\nabla G_\mu(y,x)|^{2q_0} dydx\,\lesssim\, (1+\mu)^{2q_0}
R^d R^{2q_0(1-d)} \exp\big(-c\sqrt{\mu}R\big).
\end{equation}
Since $q_0>1$ and $d\ge 2$,
\begin{equation*}
 d+2q_0(1-d)\,=\,d(1-q_0)-q_0(d-2)\, \leq \, -2q_0 \frac{q_0-1}{q_0},
\end{equation*}
\eqref{eq:quenched_a_2} implies \eqref{quenched_ap_1} for $q_1=q_0>1$ and  $\alpha_1=\frac{q_0-1}{q_0}>0$.
This result carries over to general measurable coefficients $A$ by density.
(Note that for $d>2$, this already yields the desired result \eqref{quenched_ap} for all $1\le q\le q_0$ and $\alpha_0=\frac{1}{2}$. The following two steps
are forced upon us to deal with $d=2$.)

\medskip 

\step{2} Deterministic estimates on the gradient of the Green function.

\noindent In this step we show that there exists a H\"older exponent $\alpha_2>0$ 
such that for all $L\sim 1$ and $|x|\ge R\ge 4L\sim 1$,
\begin{equation}\label{eq:quenched-grad-subopt}
\Big(\int_{B_L} |\nabla_y G_\mu(x,y)|^2dy\Big)^{\frac{1}{2}}  \,\lesssim\, \frac{e^{-c\sqrt{\mu}|x|}}{|x|^{d-2+\alpha_2}}.
\end{equation}
Since $G_\mu(x,y;A)=G_\mu(y,x;A^*)$ (where $A^*$ is the transpose of $A$) and the bounds are uniform wrt $A\in \Omega$, it is enough to prove \eqref{eq:quenched-grad-subopt} with $\nabla_y G_\mu(x,y)$
replaced by $\nabla_y G_\mu(y,x)$.
We shall first prove \eqref{eq:quenched-grad-subopt} for $d>2$ and then deduce
it for $d=2$ from the result for $d=3$ following the argument by Avellaneda and Lin
already used in Step~6 of the proof of Theorem~\ref{thm:main2}.
By Caccioppoli's inequality, for all $K\in\R$, since $L\sim 1$,
$$
\int_{B_L}|\nabla_y G_\mu(y,x)|^2dy  \,\lesssim\, \int_{B_{2L}}(G_\mu(y,x)-K)^2dy +\mu K^2 ,
$$
so that 
\begin{equation}\label{eq:quenched-grad-subopt-1}
\int_{B_L}|\nabla_y G_\mu(y,x)|^2dy  \,\lesssim\, \Big(\osc{y\in B_{2L}}{G_\mu(y,x)}\Big)^2+\Big(\sqrt{\mu}\fint_{B_{2L}} G_\mu(y,x)dy\Big)^2.
\end{equation}
From \cite[Theorem~8.22]{Gilbarg-Trudinger-98}, since $\{y:|y|\le 2L\} \subset \{y:|y|\le |\frac{x}{2}|\}$ and 
$$
\mu G_\mu(y,x)-\nabla_y\cdot A(y)\nabla_y G_\mu(y,x)\,=\,0 \ \mbox{in} \ \{y:|y|\le |\frac{x}{2}|\},
$$
we learn that there exists $\alpha_2>0$ 
such that
$$
\osc{y\in B_{2L}}{G_\mu(x,y)} \,\lesssim \,L^{\alpha_2} |\frac{x}{2}|^{-\alpha_2}(1+|\frac{x}{2}|^2\mu) \sup_{|y|\le |\frac{x}{2}|} G_\mu(x,y).
$$
Appealing to the pointwise estimate \eqref{eq:ptwise-decay-estim} for $d>2$ to bound the supremum and using that $|x-y|\ge |\frac{x}{2}|$, this turns into
\begin{equation}\label{eq:BMO-type}
\osc{y\in B_{2L}}{G_\mu(x,y)} \,\lesssim \, |\frac{x}{2}|^{2-d-\alpha_2}  e^{-c\sqrt{\mu}|\frac{x}{2}|}.
\end{equation}
Likewise the pointwise estimate \eqref{eq:ptwise-decay-estim} for $d>2$ allows one to bound the average in the RHS of \eqref{eq:quenched-grad-subopt}  by
\begin{equation*}
\sqrt{\mu}\fint_{B_{2L}} G_\mu(y,x)dy \,\lesssim \, \sqrt{\mu} |\frac{x}{2}|^{2-d}  e^{-c\sqrt{\mu}|\frac{x}{2}|}\, \lesssim 
|\frac{x}{2}|^{1-d}  e^{-c\sqrt{\mu}|\frac{x}{2}|}
\end{equation*}
for some slightly smaller $c>0$ in the RHS. 
Hence, \eqref{eq:quenched-grad-subopt} follows from \eqref{eq:quenched-grad-subopt-1} for $d>2$.

\medskip

\noindent We now turn to $d=2$, which is the aim of this step, and prove the result by integrating the three-dimensional Green function.
Denote by $A^{(2)}$ the coefficients in $\R^{2\times 2}$, and let $A^{(3)}$ be the block diagonal matrix of  $\R^{3\times 3}$ given by $\dig{A^{(2)},1}$.
We denote by $G_\mu^{(3)}$ the Green function associated with $A^{(3)}$ and define
a function $G_\mu^{(2)}:\R^2\times \R^2 \setminus \{x=y\} \to \R^+, (x,y)\mapsto G_\mu^{(2)}(x,y)$ as follows:
$$
G_\mu^{(2)}(x,y)\,=\,\int_{\R} G_\mu^{(3)}(x,z,y,0)dz.
$$
Then, $G_\mu^{(2)}=G_\mu(\cdot,\cdot;A^{(2)})$.
By the triangle inequality,
\begin{eqnarray*}
\int_{B_L} |\nabla_y G_\mu^{(2)}(y,x)|^2dy
 &=&\int_{B_L}  \Big|\int_{\R}
\nabla_y G_\mu^{(3)}(y,z,x,0)dz \Big|^2 dy\\
&\le &\bigg( \int_{\R}\Big(\int_{B_L} |\nabla_y G_\mu^{(3)}(y,z,x,0)|^2dy\Big)^{\frac{1}{2}}dz\bigg)^2.
\end{eqnarray*}
Using Cauchy-Schwarz' inequality locally, this yields
\begin{equation}\label{eq:pr-Green-6}
\int_{B_L} |\nabla G_\mu^{(2)}(y,x)|^2dy 
\,\lesssim \, \bigg( \int_{\R}\Big(\int_{|(y,z')-(x,z)|\le \frac{5L}{4}} |\nabla_y G_T^{(3)}(y,z',x,0)|^2dydz'\Big)^{\frac{1}{2}}d z\bigg)^2.
\end{equation}
We then appeal to \eqref{eq:quenched-grad-subopt} for $d=3$, which yields
\begin{equation*}
\Big(\int_{|(y,z')-(x,z)|\le \frac{5L}{4}} |\nabla_y G_T^{(3)}(y,z',x,0)|^2dydz'\Big)^{1/2}
\, \lesssim \,
\frac{e^{-c\sqrt{\mu}(|x|^2+|z|^2)^{\frac{1}{2}}}}{(|x|^2+|z|^2)^{\frac{1+\alpha_2}{2}}}.
\end{equation*}
Estimating the $z$-integral as follows,
\begin{eqnarray*}
{\int_{\R} \frac{e^{-c\sqrt{\mu}(|x|^2+|z|^2)^{\frac{1}{2}}}}{(|x|^2+|z|^2)^{\frac{1+\alpha_2}{2}}}dz }
&\le &e^{-c\sqrt{\mu}|x|} \int_{\R} \frac{1}{(|x|^2+|z|^2)^{\frac{1+\alpha_2}{2}}} dz\\
&\lesssim &\frac{e^{-c\sqrt{\mu}|x|}}{|x|^{\alpha_2}},
\end{eqnarray*}
completes the proof of \eqref{eq:quenched-grad-subopt} for $d=2$.

\medskip 

\step{3} Proof of \eqref{quenched_ap} for all $1\le q \le q_0$.

\noindent We first prove that \eqref{quenched_ap} holds for $q=1$ using  Caccioppoli's inequality combined with \eqref{eq:quenched-grad-subopt}, and then conclude by interpolation using Step~1.
Assume that $A$ is smooth, so that $\nabla_{y_i} G_\mu(y,x)$ is smooth for $x\neq y$. Since for all $i\in \{1,\dots,d\}$
$$
\mu \nabla_{y_i} G_\mu(y,x)-\nabla_x \cdot A(x) \nabla_x\nabla_{y_i}G_\mu(y,x)\,=\,0 \ \mbox{in} \ \{\frac{R}{2}\le |x|<4R\},
$$
Caccioppoli's inequality yields
$$
 \int_{|y|\le L} \int_{R\le |x|<2R}|\nabla_x\nabla_{y_i} G_\mu(y,x)|^{2}dxdy \,\lesssim\,
 R^{-2}\int_{|y|\le L} \int_{\frac{R}{2}\le |x|<4R}(\nabla_{y_i} G_\mu(y,x))^{2}dxdy .
$$
Combined with \eqref{eq:quenched-grad-subopt} this turns into
\begin{equation}\label{eq:quenched_ap_3}
 \int_{|y|\le L} \int_{R\le |x|<2R}|\nabla\nabla G_\mu(y,x)|^{2}dxdy \,\lesssim\,
 R^{-2} R^dR^{2(2-d)-2\alpha_2}=R^{2-d}R^{-2\alpha_2}e^{-c\sqrt{\mu}R},
\end{equation}
that is \eqref{quenched_ap} for $q=1$ and exponent $\alpha_2$.
The case of measurable coefficients $A$ follows by density.

Set $\alpha_0=\min \{\alpha_1,\alpha_2\}$.
An elementary interpolation argument between \eqref{eq:quenched_ap_3} and \eqref{quenched_ap_1} then shows that for all $1\le q\le q_0$,
$$
 \int_{|y|\le L} \int_{R\le |x|<2R}|\nabla\nabla G_\mu(y,x)|^{2q}dxdy \,\lesssim\,
R^{-2q\alpha_0}e^{-c\sqrt{\mu}R},
$$
as desired.

\medskip

\step{4} Proof of \eqref{loc_bound}.

\noindent This is a consequence of Caccioppoli's inequality and \eqref{quenched_ap1}.
Indeed, for all $3L\le |x-y|< 6L$ with $L\sim 1$,
\begin{eqnarray*}
(\nabla \nabla G_\mu)_L(x,y)&=&\int_{B_L(y)}\int_{B_L(x)} |\nabla_{x'} \nabla_{y'} G_\mu(x',y')|^2dx'dy'\\
&\stackrel{L\sim1,\text{Caccioppoli}}{\le} & \int_{B_L(y)}\int_{B_{\frac{3}{2}L}(x)} |\nabla_{y'} G_\mu(x',y')|^2dx'dy' \\
&\le & \int_{B_{\frac{3}{2}L}(x)} \int_{\frac{L}{2} \le |x'-y'| < \frac{17L}{2}}    |\nabla_{y'} G_\mu(x',y')|^2dy'dx'\,\stackrel{\eqref{quenched_ap1}}{\lesssim} \, 1.
\end{eqnarray*}
%


\subsection{Proof of Lemma~\ref{lem:absorb}}

We only prove \eqref{absorb}, the proof of \eqref{absorb1} is similar and left to the reader.
We split the proof of \eqref{absorb} into three steps.
In the first step we estimate the oscillation of the mixed second derivative of the Green function. In the second step we control the RHS of this estimate using Lemma~\ref{lem:quenched}, and we conclude in the third step.

We let $\tilde A$ be a coefficient field which coincides with $A$ outside of $B_L(z)$, for $z\in \R^d$, and denote by $G_\mu$ and $\tilde G_\mu$ the Green functions associated with $A$ and $\tilde A$, respectively, for some $\mu>0$.
Set $\delta G_\mu:=\tilde G_\mu-G_\mu$.

\medskip

\step{1} Proof of
\begin{equation}\label{eq:oscillation}
(\nabla \nabla \delta G_\mu)_L(x,y)\,\lesssim \, (\nabla \nabla G_\mu)_L(z,y)
\left\{
\begin{array}{ll}
1 & \text{ if }|z-x|\le 6L,\\
(\nabla \nabla  G_\mu)_L(x,z) & \text{ if }|z-x|>6L.
\end{array}
\right.
\end{equation}
for all $x,y$ with $|z-y|>3L$ 	and $|x-y|>3L$.

\noindent
By density it is enough to take $A$ and $\tilde A$ smooth.
Estimate  \eqref{eq:oscillation} follows from the combination of a Green representation formula and an a priori estimate. We start with the former and proceed by regularization.
Let $(\rho_r)_{r>0}$ be a family of smooth non-negative approximations of the Dirac mass with total mass unity and support in $B_r$.
For all $r>0$ and $y'\in \R^d$, let $G_{\mu,r}(\cdot,y')$ be the unique weak solution in $H^1(\R^d)$ of 
$$
\mu G_{\mu,r}(x',y')-\nabla_{x'} \cdot A(x')\nabla_{x'} G_{\mu,r}(x',y')\,=\,\rho_{r}(y'-x').
$$
By standard elliptic regularity theory, $G_{\mu,r}$ is smooth on $\R^d\times \R^d$.
In addition, from the existence/uniqueness theory for the Green function, we learn that for all $y'\in \R^d$,
\begin{eqnarray}
&&G_{\mu,r}(\cdot,y')\,\stackrel{r\downarrow 0}{\longrightarrow} \,G_\mu(\cdot,y')\quad \text{ in }W^{1,1}(\R^d) \label{eq:approxW11}
\end{eqnarray}
Hence, for all $y'\in \R^d$,
\begin{equation}\label{eq:approx-delta-G}
\delta G_{\mu,r}(\cdot,y'):=\tilde G_{\mu,r}(\cdot,y')-G_{\mu,r}(\cdot,y')\,\stackrel{r\downarrow 0}{\longrightarrow}\, \delta G_\mu(\cdot,y') \quad \text{ in }W^{1,1}(\R^d).
\end{equation}
For all $y'\in \R^d$, $\delta G_{\mu,r}(\cdot,y')$ is a classical solution of
$$
\mu \delta G_{\mu,r}(x',y')-\nabla_{x'} \tilde A(x') \nabla_{x'} \delta G_{\mu,r}(x',y')\,=\,
\nabla_{x'}\cdot (\tilde A-A)(x')\nabla_{x'} G_{\mu,r}(x',y').
$$
Since the RHS has compact support, $\delta G_{\mu,r}(\cdot,y')$ satisfies the Green representation formula for all $x',y'\in \R^d$
\begin{equation}\label{eq:GRF-Gmur}
\delta G_{\mu,r}(x',y')\,=\,\int_{\R^d} \nabla_{z'} \tilde G_{\mu}(x',z')\cdot
(\tilde A-A)(z')\nabla_{z'} G_{\mu,r}(z',y')dz'.
\end{equation}
Provided $|z-x'|>2L$ and $|z-y'|>2L$,  standard deterministic estimates on the gradient of the Green function yield:
$$
\sup_{z'\in  B_L(z)} |\nabla_{z'} \tilde G_{\mu}(x',z')| \,\lesssim \, \sup_{z'\in  B_L(z)} |x'-z'|^{2-d} \,\le \, L^{2-d} \sim 1.
$$
Hence, using \eqref{eq:approxW11} and \eqref{eq:approx-delta-G}, as $r\downarrow 0$, the Green representation formula \eqref{eq:GRF-Gmur} turns into
\begin{equation}\label{eq:GRF-Gmu}
\delta G_{\mu}(x',y')\,=\,\int_{\R^d} \nabla_{z'} \tilde G_{\mu}(x',z')\cdot
(\tilde A-A)(z')\nabla_{z'} G_{\mu}(z',y')dz'
\end{equation}
for all  $|z-x'|>2L$ and $|z-y'|>2L$.
Since $G_\mu$ and $\tilde G_\mu$ are smooth away from the diagonal, we may differentiate twice \eqref{eq:GRF-Gmu}, which yields for all $|z-x'|>2L$ and $|z-y'|>2L$,
\begin{equation}\label{eq:GRF-GmuD2}
\nabla \nabla \delta G_{\mu}(x',y')\,=\,\int_{\R^d} \nabla \nabla \tilde G_{\mu}(x',z')\cdot
(\tilde A-A)(z')\nabla\nabla G_{\mu}(z',y')dz'.
\end{equation}
Recall that $|z-x|>3L$ and $|z-y|>3L$.
Integrating \eqref{eq:GRF-GmuD2} over $x'\in B_L(x)$ and $y'\in B_L(y)$, we obtain by Cauchy-Schwarz' inequality
\begin{equation}\label{eq:oscillation-1}
(\nabla \nabla \delta G_{\mu})_L(x,y)\,\lesssim\,(\nabla \nabla \tilde G_{\mu})_L(x,z)(\nabla \nabla G_{\mu})_L(z,y).
\end{equation}

We turn now to the a priori estimate.
Let $|y'-z|>2L$. Then, $\delta G_\mu(\cdot,y')$ is the unique distributional solution in $W^{1,1}(\R^d)$ of
$$
\mu \delta G_\mu(x',y')-\nabla_{x'} \tilde A(x') \nabla_{x'} \delta G_{\mu}(x',y')\,=\,
\nabla_{x'}\cdot (\tilde A-A)(x')\nabla_{x'} G_{\mu}(x',y').
$$
Since $G_{\mu}$ is smooth away from the diagonal, the RHS is smooth with compact support, so that $\delta G_\mu(\cdot,y')$ is a classical solution.
We then differentiate the equation with respect to $y'_i$ for $i\in \{1,\dots,d\}$:
$$
\mu \nabla_{y'_i}\delta G_\mu(x',y')-\nabla_{x'} \tilde A(x') \nabla_{x'} \nabla_{y'_i}\delta G_{\mu}(x',y')\,=\,
\nabla_{x'}\cdot (\tilde A-A)(x')\nabla_{x'} \nabla_{y'_i}G_{\mu}(x',y').
$$
Since the RHS is smooth and has compact support, $\nabla_{y'_i}\delta G_{\mu}(\cdot,y') \in H^1(\R^d)$, and we may test the weak formulation of the equation with the solution itself. This yields
$$
\int_{\R^d} |\nabla \nabla\delta G_{\mu}(x',y')|^2 \;dx' \,\lesssim\, \int_{B_L(z)} |\nabla \nabla \delta G_{\mu}(x',y')| |\nabla \nabla G_{\mu}(x',y')| \;dx',
$$
which, by Young's inequality,  turns into
\begin{equation}\label{eq:oscillation-2}
\int_{\R^d} |\nabla \nabla\delta G_{\mu}(x',y')|^2 \;dx' \,
\lesssim
\,\int_{B_L(z)} |\nabla \nabla  G_{\mu}(x',y')|^2 \;dx'.
\end{equation}

We are in position to conclude. 
On the one hand, integrating \eqref{eq:oscillation-2} over $y'\in B_L(y)$ yields
\begin{equation}\label{eq:oscillation-21}
(\nabla \nabla\delta G_{\mu})_L(x,y) \,
\lesssim
\,(\nabla \nabla  G_{\mu})_L(z,y).
\end{equation}
On the other hand, assume that $|z-x|>3L$.
Denote by $G^*_{\mu}$, $\tilde G^*_{\mu}$ and $\delta G^*_{\mu}$ the Green functions associated with $A^*$, $\tilde A^*$, and their difference. Estimate \eqref{eq:oscillation-2} takes the form
$$
\int_{\R^d} |\nabla \nabla\delta G_{\mu}^*(y',x)|^2 \;dx' \,
\lesssim
\,\int_{B_L(z)} |\nabla \nabla  G_{\mu}^*(y',x)|^2 \;dx',
$$
so that by integration over $y'\in B_L(x)$ and by the symmetry properties of the Green function,
$$
(\nabla \nabla\delta G_{\mu})_L(x,z) \,=\,(\nabla \nabla\delta G_{\mu}^*)_L(z,x)\,\lesssim \,(\nabla \nabla  G^*_{\mu})_L(z,x)
\,=\,(\nabla \nabla  G_{\mu})_L(x,z).
$$
Hence by the triangle inequality, the estimate \eqref{eq:oscillation-1} for $|z-x|>3L$ turns into 
\begin{equation}\label{eq:oscillation-11}
(\nabla \nabla \delta G_{\mu})_L(x,y)\,\lesssim\,(\nabla \nabla  G_{\mu})_L(x,z)(\nabla \nabla G_{\mu})_L(z,y).
\end{equation}
The claim \eqref{eq:oscillation} follows from the combination of \eqref{eq:oscillation-21} and \eqref{eq:oscillation-11}.

\medskip

\step{2} Proof of
\begin{equation}\label{eq:int-DDG}
\sup_{x'\in \R^d} \int_{\R^d} (|z-x'|+1)^{2q\alpha} 
\left\{
\begin{array}{ll}
1 & \text{ if }|z-x'|\le 6L\\
(\nabla \nabla  G_\mu)_L^{2q}(x',z) & \text{ if }|z-x'|>6L
\end{array}
\right\}
\;dz \,\lesssim \,1
\end{equation}
for all $1\le q\le q_0$ and $\alpha=\frac{\alpha_0}{2}$, where $q_0$ and $\alpha_0$ are as in Lemma~\ref{lem:quenched}.
For $|z-x'|$ small, we have
\begin{multline}\label{eq:ball-dyad}
 \int_{|z-x'|\le 6L} (|z-x'|+1)^{2q\alpha} \left\{
\begin{array}{ll}
1 & \text{ if }|z-x'|\le 6L\\
(\nabla \nabla  G_\mu)_L^{2q}(x',z) & \text{ if }|z-x'|>6L
\end{array}
\right\}\;dz 
\\
\le \int_{|z-x'|\le 6L} (|z-x'|+1)^{q\alpha_0} \;dz \lesssim
1.
\end{multline}
For larger $|z-x'|$, we decompose $\{z:|z-x'|>6L\}$ into dyadic annuli:
\begin{multline}\label{eq:dyad-sum}
 \int_{|z-x'|>6L} (|z-x'|+1)^{2q\alpha} 
\left\{
\begin{array}{ll}
1 & \text{ if }|z-x'|\le 6L\\
(\nabla \nabla  G_\mu)_L^{2q}(x',z) & \text{ if }|z-x'|>6L
\end{array}
\right\} \;dz\\ \le \sum_{n=0}^\infty \int_{2^n 6L < |z-x'| \le
2^{n+1} 6L} (|z-x'|+1)^{2q\alpha} (\nabla \nabla G_\mu)_L^{2q}(x',z) \;dz.
\end{multline}
On each dyadic annulus, 
\begin{multline*}
 \int_{2^n 6L < |z-x'| \le 2^{n+1} 6L}  (|z-x'|+1)^{2q\alpha} (\nabla \nabla G_\mu)_L^{2q}(x',z) \;dz\\
 \lesssim 2^{2q\alpha n} \int_{2^n 6L < |z-x'| \le 2^{n+1} 6L}  \bigg(\int_{B_L} \int_{B_L} |\nabla \nabla G_\mu(x'+x'',z+z')|^2 \;dx''dz'\bigg)^{q}
\;dz,
\end{multline*}
which we bound using Jensen's inequality and \eqref{quenched_ap} as
\begin{align}
 &2^{2q\alpha n} \int_{2^n 6L < |z-x'| \le 2^{n+1} 6L}  \bigg(\int_{B_L} \int_{B_L} |\nabla \nabla G_\mu(x'+x'',z+z')|^2 \;dx''dz'\bigg)^{q} \;dz\nonumber \\
 &\lesssim 2^{2q\alpha n} \int_{B_L} \int_{B_L} \int_{2^n 6L < |z-x'| \le 2^{n+1} 6L}  |\nabla \nabla G_\mu(x'+x'',z+z')|^{2q} \;dzdx''dz' \nonumber \\
 &\lesssim 2^{2q\alpha n} \int_{B_L} \int_{B_L} \int_{2^n 4L < |z-x'-x''| \le 2^{n+1} 8L}  |\nabla \nabla G_\mu(x'+x'',z)|^{2q} \;dzdx''dz' \nonumber \\
 &\stackrel{\eqref{quenched_ap}}{\lesssim} 2^{2q(\alpha-\alpha_0) n} = 2^{-q\alpha_0 n},\label{eq:annulus-dyad}
\end{align}
uniformly wrt $x'\in \R^d$.
The combination of \eqref{eq:ball-dyad}, \eqref{eq:dyad-sum}, and \eqref{eq:annulus-dyad} yields the claim \eqref{eq:int-DDG} since $\sum_{n\in\N} 2^{-q\alpha_0 n} \lesssim 1$.

\medskip

\step{3} Conclusion.

\noindent We first show that for all $|x-y|>6L$ and all $p$ large enough, we have
\begin{multline}\label{eq:pr-absorb-step3-1}
 \expec{\bigg(\int_{\R^d} \Big( \osc{B_L(z)}{(\nabla \nabla G_\mu)_L(x,y)}\Big)^2 |x-y|^{2d}e^{2c\sqrt{\mu}|x-y|}dz\bigg)^p}
 \\
  \lesssim\,
\sup_{z,y:|z-y|>3L} \Big\{ |z-y|^{2pd} e^{2c\sqrt{\mu}|z-y|}\expec{(\nabla\nabla G_\mu)_L^{2p}(z,y)} \Big\} .
\end{multline}
We claim that it is enough to prove that
\begin{multline}\label{eq:pr-absorb-step3-2}
 \expec{\bigg(\int_{|z-y|\ge |z-x|} \Big( \osc{B_L(z)}{(\nabla \nabla G_\mu)_L(x,y) \Big)^2}  |x-y|^{2d}dz\bigg)^p}
 \\
  \lesssim\,
\sup_{z,y:|z-y|>3L} \Big\{ |z-y|^{2pd} \expec{(\nabla\nabla G_\mu)_L^{2p}(z,y)} \Big\}.
\end{multline}
To this aim we have to prove that the corresponding integral on the LHS of \eqref{eq:pr-absorb-step3-2}, this time over $\{|z-y|\leq |z-x|\}$, is bounded by the
RHS of \eqref{eq:pr-absorb-step3-2}. Indeed, \eqref{eq:pr-absorb-step3-2} for $G^*_\mu$ with $x$ and $y$ switched takes the form after using the symmetry properties of the Green function
\begin{multline*}
 \expec{\bigg(\int_{|z-x|\ge |z-y|} \Big( \osc{B_L(z)}{(\nabla \nabla G_\mu)_L(x,y) \Big)^2}  |y-x|^{2d}dz\bigg)^p}
 \\
  \lesssim\,
\sup_{z,x:|z-x|>3L} \Big\{ |z-x|^{2pd} \expec{(\nabla\nabla G_\mu)_L^{2p}(x,z)} \Big\} .
\end{multline*}
The conclusion follows by stationarity since
\begin{eqnarray*}
\sup_{z,x:|z-x|>3L} \Big\{ |z-x|^{2pd} \expec{(\nabla\nabla G_\mu)_L^{2p}(x,z)} \Big\}&=&\,\sup_{z:|z|>3L} \Big\{ |z|^{2pd} \expec{(\nabla\nabla G_\mu)_L^{2p}(z,0)} \Big\}
\\
&=&\sup_{z,y:|z-y|>3L} \Big\{ |z-y|^{2pd} \expec{(\nabla\nabla G_\mu)_L^{2p}(z,y)} \Big\}.
\end{eqnarray*}
It is therefore enough to prove \eqref{eq:pr-absorb-step3-2}.

For $|z-y|\ge |z-x|$, we
have $|z-y|\ge \frac{|x-y|}{2} \ge 3L$, so that taking the supremum over $\tilde A$ (by a density argument the supremum can be taken on smooth fields $\tilde A$) in the estimate \eqref{eq:oscillation} of Step~1 yields
\begin{multline}\label{eq:pr-absorb-step3-3}
\int_{|z-y|\ge |z-x|} \Big( \osc{B_L(z)}{(\nabla \nabla G_\mu)_L(x,y)}\Big)^2 dz \\ \lesssim
 \int_{|z-y|\ge |z-x|}  (\nabla \nabla G_\mu)_L^2(z,y)
\left\{
\begin{array}{ll}
1 & \text{ if }|z-x|\le 6L\\
(\nabla \nabla  G_\mu)_L^2(x,z)& \text{ if }|z-x|>6L
\end{array}
\right\} dz.
\end{multline}
We smuggle in the weight $(|z-x|+1)^{\alpha}$ and apply H\"older's inequality 
with exponents $(p,q)$ for some $p> 1$ to be fixed below:
\begin{multline*}
\expec{\bigg(\int_{|z-y|\ge |z-x|} (\nabla \nabla G_\mu)_L^2(z,y) 
\left\{
\begin{array}{ll}
1 & \text{ if }|z-x|\le 6L\\
(\nabla \nabla  G_\mu)_L^2(x,z)& \text{ if }|z-x|>6L
\end{array}
\right\} dz\bigg)^p} \\
 \lesssim \, 
\,
\expec{\bigg(  \int_{|z-y| > 3L} (|z-x|+1)^{2q\alpha}  \left\{
\begin{array}{ll}
1 & \text{ if }|z-x'|\le 6L\\
(\nabla \nabla  G_\mu)_L^{2q}(x',z) & \text{ if }|z-x'|>6L
\end{array}
\right\}dz
 \bigg)^{\frac{p}{q}}}\\
\times \,\expec{\int_{|z-y|\ge |x-z|} (|z-x|+1)^{-2p\alpha} (\nabla \nabla G_\mu)_L^{2p}(z,y)dz}.
\end{multline*}
By \eqref{eq:int-DDG} in Step~2, the first term on the r.~h.~s. is bounded uniformly wrt $A$ as long as $1\le q\le q_0$, i.e.\ $p=\frac{q}{q-1}\ge \frac{q_0}{q_0-1} =: p_0$. Hence, using that $|z-y|\ge |x-y|/2$, this yields
\begin{multline*}
\expec{\bigg(|x-y|^{2d}e^{2c\sqrt{\mu}|x-y|} \int_{|z-y|\ge |z-x|}  (\nabla \nabla G_\mu)_L^2(z,y)
\left\{
\begin{array}{ll}
1 & \text{ if }|z-x|\le 6L\\
(\nabla \nabla  G_\mu)_L^2(x,z)& \text{ if }|z-x|>6L
\end{array}
\right\}dz\bigg)^p}^{\frac{1}{p}} \\
 \lesssim \,  \expec{\int_{|z-y|\ge |z-x|} (|z-x|+1)^{-2p\alpha} |z-y|^{2pd} e^{2pc\sqrt{\mu}|z-y|}(\nabla \nabla G_\mu)_L^{2p}(z,y) dz }^{\frac{1}{p}},
\end{multline*}
We then take the supremum of the last two factors of the integrand using that 
$|z-y|>3L$ and choose $p$ large enough so that $\int_{\R^d} (|x-z|+1)^{-2p\alpha} \;dz \lesssim 1$ (up to redefining $p_0$ accordingly) so that
\begin{multline}\label{eq:pr-absorb-step3-4}
\expec{\bigg(|x-y|^{2d} e^{2c\sqrt{\mu}|x-y|}\int_{|z-y|\ge |z-x|}   (\nabla \nabla G_\mu)_L^2(z,y)
\left\{
\begin{array}{ll}
1 & \text{ if }|z-x|\le 6L\\
(\nabla \nabla  G_\mu)_L^2(x,z)& \text{ if }|z-x|>6L
\end{array}
\right\}dz
\bigg)^p}^{\frac{1}{p}}  \\
 \lesssim \,  \sup_{z:|z-y|>3L}\Big\{|z-y|^{2pd} e^{2pc\sqrt{\mu}|z-y|}\expec{(\nabla \nabla G_\mu)_L^{2p}(z,y)}\Big\}^{\frac{1}{p}} .
\end{multline}
Estimate~\eqref{eq:pr-absorb-step3-2}, which implies \eqref{eq:pr-absorb-step3-1},  is now a consequence of \eqref{eq:pr-absorb-step3-3} and \eqref{eq:pr-absorb-step3-4}. 

Lemma~\ref{lem:absorb} then follows from \eqref{eq:pr-absorb-step3-1} combined with the local boundedness estimate~\eqref{loc_bound} in the form of 
\begin{multline*}
 \sup_{z:|z-y|>3L} \Big\{ |z-y|^{2pd} e^{2pc\sqrt{\mu}|z-y|}\expec{(\nabla \nabla G_\mu)_L^{2p}(z,y)}  \Big\}\\
\lesssim\, 1 + \sup_{z:|z-y|>6L} \Big\{ |z-y|^{2pd}e^{2c\sqrt{\mu}|z-y|}  \expec{(\nabla \nabla G_\mu)_L^{2p}(z,y)}  \Big\}.
\end{multline*}

\subsection{Proof of Lemma~\ref{lem:DGNM} for $p<\infty$}

The proof consists in a minor modification of the usual Moser iteration. We follow the proof of~\cite[Theorem~8.17]{Gilbarg-Trudinger-98} and mainly focus on the differences. 
Without loss of generality we may assume that $q<\frac{d}{2}$.
Up to multiplying the equation by $-1$ it is enough to prove the claim for the positive part $u^+=\max\{0,u\}$
of $u$.
Set $\bar u = u^+ + k$, where $k := \|f\|_{L^q(B_2)}$ with $q$ given in the statement. We test the equation~\eqref{u} with the test function $v = \eta^2 (\bar u^\beta - k^\beta)\ge 0$, where $\beta > 0$ and $\eta$ is a smooth cut-off function for $B_1$ in $B_2$ with $0\le\eta\le1$. 
In the following, we require that
\begin{equation}\label{beta_bound}
0< \beta < \frac{(q-1)d}{d-2q}.
\end{equation}
The derivative of $v$ is given by
$$
\nabla v = 2 \eta ( \bar u^\beta - k^\beta ) \nabla \eta + \eta^2 \beta \bar u^{\beta - 1} \nabla \bar u.
$$
Since by construction $\mu u(\bar u^\beta-k^\beta)\ge 0$
and either $\nabla \bar u$ and $\bar u^\beta-k^\beta$ vanish or $\nabla \bar u$ equals $\nabla u$, equation~\eqref{u} with test-function $v$ yields
\begin{align*}
 0 &= \int_{\R^d} \big(\mu vu+ \nabla v \cdot A \nabla u - v f \big) \;dx\\
 &= \int_{\R^d} \Big( \mu \eta^2 (\bar u^\beta-k^\beta)u+\beta \eta^2 \bar u^{\beta-1} \nabla \bar u\cdot A \nabla \bar u + \big(2\eta (\bar u^\beta-k^\beta) \nabla \eta \big)\cdot A\nabla \bar u   - v f \Big) \;dx\\
 &\ge \int_{\R^d} \Big( \lambda \beta \eta^2 \bar u^{\beta-1} |\nabla \bar u|^2 - 2 |\nabla \eta| (\bar u^\beta -k^\beta)\eta |\nabla \bar u| - \eta^2(\bar u^{\beta}-k^\beta) |f| \Big) \;dx \\
&\ge \int_{\R^d} \Big( \lambda \beta \eta^2 \bar u^{\beta-1} |\nabla \bar u|^2 - 2 |\nabla \eta| \bar u^\beta\eta |\nabla \bar u| - \eta^2 \bar u^{\beta}|f| \Big) \;dx.
\end{align*}
By Young's inequality,
$$
 \int_{\R^d} 2 |\nabla \eta| \bar u^\beta \eta |\nabla \bar u| \;dx \le \frac{\lambda \beta}{2} \int_{\R^d} \eta^2 \bar u^{\beta-1} |\nabla \bar u|^2 \;dx + \frac{2}{\lambda \beta} \int_{\R^d} |\nabla \eta|^2 \bar u^{\beta+1} \;dx,
$$
so that
\begin{equation}\label{eq:moser-1}
\frac{\lambda \beta}{2} \int_{\R^d} \eta^2 \bar u^{\beta-1} |\nabla \bar u|^2 \;dx\,\leq \,
 \frac{2}{\lambda \beta} \int_{\R^d} |\nabla \eta|^2 \bar u^{\beta+1} \;dx+\int_{\R^d} \eta^2 \bar u^{\beta}|f| \;dx.
\end{equation}
So far, the computations are identical to the usual Moser iteration. 
Here comes the difference: Let $\chi = \frac{d}{d-2}$ if $d>2$ (or fix any $1<\chi < +\infty$ if $d=2$) and let $s\ge 1$ be such that $1 = \frac{1}{s} + \frac{\beta}{(\beta+1)\chi}$.
Then, the choice~\eqref{beta_bound} implies that $1\le s<q$.
Indeed,
\begin{equation*}
\frac{1}{s}-\frac{1}{q}\,=\,1-\frac{\beta}{(\beta+1)\chi}-\frac{1}{q} 
\,>\, \frac{q-1}{q}-\frac{\beta}{\beta+1}\frac{d-2}{d}
\,=\,\frac{-\beta(d-2q)+d(q-1)}{q(\beta+1)d}\,\stackrel{\eqref{beta_bound}}{>}\,0.
\end{equation*}
We now treat the second RHS term of \eqref{eq:moser-1}.
H\"older's inequality on $\eta^2\bar u^\beta|f|=(\eta^{2\frac{\beta}{\beta+1}}\bar u^\beta)(\eta^{2\frac{1}{\beta+1}}|f|)$ with exponents $(\chi\frac{\beta+1}{\beta},s)$ yields
\begin{equation*}
 \int_{\R^d} \eta^2 \bar u^\beta|f| \;dx\,\le \, \bigg( \int_{\R^d} \eta^{2\chi} \bar u^{(\beta+1)\chi}\;dx \bigg)^{\frac{\beta}{(\beta+1)\chi}} \bigg( \int_{\R^d} \eta^{\frac{2 s}{\beta+1}} |f|^s \;dx \bigg)^{\frac{1}{s}}.
\end{equation*}
Let $C$ denote a generic constant depending only on $d$, $\lambda$ and $q$ (but which can change from line to line) --- note that since $1\le\beta<\frac{(q-1)d}{d-2q}$, constants depending on $\beta$ are also bounded  by $C$.
Since $0\le \eta \le 1$ and $s<q < \frac{d}{2}$, it follows by Jensen's inequality that
$$
\int_{\R^d} \eta^2 \bar u^\beta|f| \;dx \le C \ k \; \bigg( \int_{\R^d} \eta^{2\chi} \bar u^{(\beta+1)\chi}\;dx \bigg)^{\frac{\beta}{(\beta+1)\chi}},
$$
where we recall that $k=\|f\|_{L^q(B_2)}$. By Young's inequality we thus have for all $\e > 0$:
$$
\int_{\R^d} \eta^2 \bar u^\beta|f| \;dx \le \e \bigg( \int_{\R^d} \eta^{2\chi} \bar u^{(\beta+1)\chi}\;dx \bigg)^{\frac{1}{\chi}} + Ck^{\beta+1},
$$
where $C$ depends now in addition on $\e$.
Combined with \eqref{eq:moser-1} this yields
$$
 \frac{\lambda \beta}{2} \int_{\R^d} \eta^2 \bar u^{\beta-1} |\nabla \bar u|^2 \;dx \le \frac{2}{\lambda \beta}\int_{\R^d}  |\nabla \eta|^2 \bar u^{\beta+1} \;dx + Ck^{\beta+1}+\e \bigg( \int_{\R^d} \eta^{2\chi} \bar u^{(\beta+1)\chi}\;dx \bigg)^{\frac{1}{\chi}}.
$$
Next we introduce another function $w:=\bar u^{\frac{\beta+1}{2}}$ and rewrite this inequality as
\begin{equation}\label{moser2}
 \lambda\int_{\R^d} \eta^2 |\nabla w|^2 \;dx \le \frac{2}{\lambda \beta}\int_{\R^d}  |\nabla \eta|^2w^2 \;dx + Ck^{\beta+1}+ \e \bigg( \int_{\R^d} |\eta w|^{2\chi} \;dx \bigg)^{\frac{1}{\chi}}.
\end{equation}
This yields
$$
 \int_{\R^d} |\nabla (\eta w)|^2 \;dx \le  \frac{2}{\lambda^2 \beta}\int_{\R^d}|\nabla \eta|^2  w^2 \;dx + Ck^{\beta+1}+\frac{\e}{\lambda} \bigg( \int_{\R^d} |\eta w|^{2\chi} \;dx \bigg)^{\frac{1}{\chi}}.
$$
By the Sobolev embedding, this turns into
$$
C_{\mathrm{Sob}}\bigg( \int_{\R^d} |\eta w|^{2\chi} \;dx \bigg)^{\frac{1}{\chi}} \le \frac{2}{\lambda^2 \beta}\int_{\R^d}|\nabla \eta|^2  w^2 \;dx + Ck^{\beta+1}+\frac{\e}{\lambda} \bigg( \int_{\R^d} |\eta w|^{2\chi} \;dx \bigg)^{\frac{1}{\chi}},
$$
so that for $\e$ small enough (and only depending on $d$, $\lambda$ and $q$), we have
$$
 \bigg( \int_{\R^d} |\eta w|^{2\chi} \;dx \bigg)^{\frac{1}{\chi}} \le C\int_{\R^d}  |\nabla \eta|^2 w^2 \;dx
+Ck^{\beta+1}.
$$
This corresponds to the usual Moser iteration (albeit the dependence of the constants on $\beta$ is worse),
and yields the desired result for $p=\chi(\beta+1)$.
We can then iterate by increasing $\beta$ to yield bounds as long as $\beta < \frac{(q-1)d}{d-2q}$. In this case any exponent of the form $p = (\beta+1)\chi$ can be attained, which yields
$$
 \frac{1}{p} > \frac{d-2q}{((q-1)d+d-2q)\chi} = \frac{d-2q}{dq} = \frac{1}{q} - \frac{2}{d},
$$
as claimed.
Note that (unlike the usual Moser iteration) the dependence of the constants on $\beta$ does not matter since we only need to iterate finitely many times in order to reach $p<+\infty$. 

\subsection{Proof of Lemma~\ref{lem:osc_u}}

Let $u$ and $\tilde u$ be solutions of~\eqref{u} with coefficient fields $A$ and $\tilde A$, respectively, 
where the coefficients coincide outside the ball $B_\ell(z)$.
Their difference $\delta u$ solves
\begin{eqnarray}
\mu\delta u -\nabla\cdot \tilde A \nabla \delta u &=& -\nabla \cdot ( \tilde A - A ) \nabla u,\label{eq:pr-lem-osc00}\\
\mu\delta u -\nabla\cdot A \nabla \delta u &=& -\nabla \cdot ( \tilde A - A ) \nabla \tilde u.\label{eq:pr-lem-osc01}
\end{eqnarray}

\medskip

\step{1} Preliminary result and proof of
$$
\sup_{x,x'\in \R^d} \sup_{\tilde A\in \Omega} \bigg\| \fint_{B_{\frac{3\ell}{2}}(x)} \nabla_{y'} \tilde G_\mu(y,y') \;dy \bigg\|_{L^2_{y'}(B_{\frac{3\ell}{2}}(x'))} \lesssim 1.
$$
To see this, we note that 
$$
 v:y'\mapsto \fint_{B_{\frac{3\ell}{2}}(x)} \tilde G_\mu(y,y') \;dy \quad \text{solves}\quad \mu v-\nabla\cdot \tilde A \nabla v = \frac{1}{|B_{\frac{3\ell}{2}}(x)|}\chi_{B_{\frac{3\ell}{2}}(x)},
$$
where $\chi_D$ denotes the characteristic function of the set $D\subseteq \R^d$, that is, a regularized version
of the defining equation for the Green function without singularity.
The proof that 
$\int_{B_{\frac{3\ell}{2}}(x')}|\nabla v|^2dy'$ is bounded and only depends on 
$\ell$ and $\lambda$ is similar to the corresponding proof of \cite[Corollary~2.3]{Gloria-Otto-09} in the discrete case (since there is no singularity to be taken care of).

\medskip

\step{2} Proof of \eqref{eq:osc-estim-u} for $|x-z| > 6\ell$.

\noindent The Green function representation formula associated with \eqref{eq:pr-lem-osc01} yields
$$
 u(x) - \tilde u(x) = \int_{B_\ell(z)} \nabla_z G_\mu(x,z') \cdot (\tilde A(z') - A(z')) \nabla \tilde u(z') \;dz'.
$$
Hence, by the triangle inequality and H\"older's inequality,
$$
 \| u - \tilde u \|_{L^{\lambda_1'}(B_\ell(x))} \lesssim \| \nabla_2  G_\mu \|_{L_x^{\lambda_1'}(B_\ell(x), L_z^2(B_\ell(z)))} (\nabla \tilde u)_\ell(z),
$$
where we recall that $ (\nabla \tilde u)_\ell(z)=\| \nabla\tilde u \|_{L^2(B_\ell(z))}$.
Since $|x-z| > 6\ell$, for all $i\in \{1,\dots,d\}$ the function $x\mapsto \nabla_{z_i} G_\mu(x,z)$ is in the kernel of $(\mu-\nabla \cdot A \nabla)$ in $B_{2\ell}(x)$ for all $z\in B_\ell(z)$ and Lemma~\ref{lem:DGNM} implies that
$$
 \| \nabla_2 G_\mu \|_{L_x^{\lambda_1'}(B_\ell(x),L_z^2(B_\ell(z)))} \lesssim \| \nabla_2  G_\mu \|_{L_x^{2}(B_{2\ell}(x)\times B_\ell(z))} \le (\nabla_z G_\mu)_{2\ell}(x,z).
$$
On the other hand, an energy estimate based on \eqref{eq:pr-lem-osc00} yields
$$
(\nabla \tilde u)_\ell(z) \lesssim (\nabla u)_\ell(z).
$$
Estimate~\eqref{eq:osc-estim-u} for $|x-z|>6\ell$ is proved.

\medskip

\step{3} Proof of \eqref{eq:osc-estim-u} for $|x-z|\le6\ell$. 

\noindent Let $x$ be fixed such that $|x-z|\le6\ell$.
We shall consider a third coefficient field $A_0\in \Omega$ such that $A_0|_{\R^d\setminus B_{9\ell}(z)}=A|_{\R^d\setminus B_{9\ell}(z)}$, $A_0|_{B_{8\ell}(z)}=\Id$, and denote by $u_0$ the associated solution of \eqref{u} with coefficient fields $A_0$.
We denote the local averages of $u$, $\tilde u$, and $u_0$ around $x$ by
\begin{equation*}
 \bar u= \fint_{B_{\frac{3\ell}{2}}(x)} u\;dy ,\quad\bar{\tilde u} = \fint_{B_{\frac{3\ell}{2}}(x)} \tilde u\;dy, \quad \text{and }
\bar u_0 = \fint_{B_{\frac{3\ell}{2}}(x)} u_0\;dy.
\end{equation*}
The triangle inequality yields
\begin{equation}\label{eq:triang-u-tilde-0}
 \| u - \tilde u \|_{L^{\lambda_1'}(B_\ell(x))} \lesssim \| u - \bar u_0 \|_{L^{\lambda_1'}(B_\ell(x))} + \| \tilde u - \bar u_0 \|_{L^{\lambda_1'}(B_\ell(x))} .
\end{equation}
By the De Giorgi-Nash-Moser estimate of Lemma~\ref{lem:DGNM} with $p=\lambda_1'$ and $q=\lambda_2$ (note $\frac{1}{\lambda_2} < \frac{2}{d} + \frac{1}{\lambda_1'}$ and $u -\bar u_0$ solves the same equation as $u$ with the addition of $-\mu \bar u_0$ on the RHS), the triangle inequality, and Poincar\'{e}'s inequality on $B_{\frac{3\ell}{2}}(x)$, the first term yields
\begin{eqnarray}
 \label{eq:u-baru0} \| u - \bar{u}_0 \|_{L^{\lambda_1'}(B_\ell(x))} &\lesssim& \Big(\int_{B_{\frac{3\ell}{2}}(x)} |u(y) - \bar u_0|^2 \;dy\Big)^{\frac{1}{2}}  + \|\mu \bar u_0+ f \|_{L^{\lambda_2}(B_{\frac{3\ell}{2}}(x))} 
\\ 
\nonumber &\lesssim &\Big(\int_{B_{\frac{3\ell}{2}}(x)} |u(y) - \bar u|^2 \;dy\Big)^{\frac{1}{2}} + |\bar u_0-\bar u|+ \mu |\bar u_0|+\| f \|_{L^{\lambda_2}(B_{\frac{3\ell}{2}}(x))} 
\\
&\lesssim& \Big(\int_{B_{\frac{3\ell}{2}}(x)} |\nabla u|^2 \;dy\Big)^{\frac{1}{2}}  +|\bar u_0-\bar u|+ \mu |\bar u_0|+ \| f \|_{L^{\lambda_2}(B_{\frac{3\ell}{2}}(x))}.
\nonumber
\end{eqnarray}
Likewise,
\begin{eqnarray}
\label{eq:tildeu-baru0}  \| \tilde u - \bar{u}_0 \|_{L^{\lambda_1'}(B_\ell(x))} 
&\lesssim& \Big(\int_{B_{\frac{3\ell}{2}}(x)} |\nabla \tilde u|^2 \;dy\Big)^{\frac{1}{2}}  +|\bar u_0-\bar{\tilde u}|+ \mu |\bar u_0|+ \| f \|_{L^{\lambda_2}(B_{\frac{3\ell}{2}}(x))}.
\end{eqnarray}
On the one hand, an energy estimate based on  \eqref{eq:pr-lem-osc00} yields
\begin{equation}\label{eq:energy_nablau-tilde-0}
 \int_{B_{9\ell}(z)} |\nabla \tilde u|^2 \;dy \lesssim \int_{B_{9\ell}(z)} |\nabla u|^2 \;dy,
\quad \int_{B_{9\ell}(z)} |\nabla  u_0|^2 \;dy \lesssim \int_{B_{9\ell}(z)} |\nabla u|^2 \;dy.
\end{equation}
It remains to bound $\mu |\bar u_0|$ and $|\bar u_0-\bar u|$ and $|\bar u_0-\bar{\tilde u}|$.
We start with the two differences. The Green representation formula yields
$$
 \bar u -\bar u_0 = \fint_{B_{\frac{3\ell}{2}}(x)} (u - u_0)\;dy = \fint_{B_{\frac{3\ell}{2}}(x)} \int_{B_{9\ell}(z)} \nabla_{y'}  G_{\mu,0}(y,y') \cdot (A_0 - A)(y') \nabla u(y') \;dy'dy,
$$
so that
\begin{equation*}
 |\bar u -\bar u_0| \,\lesssim \, \int_{B_{9\ell}(z)} \Big| \fint_{B_{\frac{3\ell}{2}}(x)} \nabla_{y'} G_{\mu,0}(y,y') \;dy \Big| |\nabla u(y')| \;dy'.
\end{equation*}
Proceeding also the same way for $|\bar u_0-\bar{\tilde u}|$, we conclude by Cauchy-Schwarz' inequality and Step~1 that
\begin{equation}\label{bar_diff}
 |\bar u -\bar u_0|+ |\bar {\tilde u} -\bar u_0|\, \lesssim \, \| \nabla u \|_{L^2(B_{9\ell}(z))} = (\nabla u)_{9\ell}(z).
\end{equation}
We turn to the estimate of $\mu \bar u_0$ and recall that by the choice of $A_0$, $u_0$ solves in $B_{8\ell}(z)$
$$
\mu u_0-\triangle u_0\,=\,f.
$$
Hence the function $u_\ell:y \mapsto \fint_{B_{\frac{3\ell}{2}}(y)} u_0dy'$ solves in $B_{\frac{\ell}{2}}(x)$ the equation
$$
\mu u_\ell-\triangle u_\ell\,=\,f_\ell,
$$
where $f_\ell(y):=\fint_{B_{\frac{3\ell}{2}}(y)} fdy'$. Testing this equation with test-function $\eta^2u_\ell$ 
with $\eta$ supported in $B_{\frac{\ell}{2}}(x)$ yields
$$
\mu\int_{B_{\frac{\ell}{2}}(x)} \eta^2 u_\ell^2dy+\int_{B_{\frac{\ell}{2}}(x) }\eta^2|\nabla u_\ell|^2dy \,=\,
\int_{B_{\frac{\ell}{2}}(x)} f_\ell \eta^2 u_\ell dy-2\int_{B_{\frac{\ell}{2}}(x)} u_\ell\eta \nabla \eta \cdot \nabla u_\ell dy,
$$
which turns, by Young's inequality, into
$$
\mu^2 \int_{B_{\frac{\ell}{2}}(x)} \eta^2 u_\ell^2dy \,\lesssim \, 
\int_{B_{\frac{\ell}{2}}(x)} f_\ell^2dy+\int_{B_{\frac{\ell}{2}}(x)} |\nabla u_\ell|^2dy.
$$
With $\eta$ a cut-off for $B_{\frac{\ell}{4}}(x)$ in $B_{\frac{\ell}{2}}(x)$, Lemma~\ref{lem:DGNM} with $p=\infty$  yields for $q=d$
$$
\mu^2 \sup_{B_{\frac{\ell}{4}}(x)} u_\ell^2 \,\lesssim \, \Big(\int_{B_{\frac{\ell}{2}}(x)} |f_\ell|^{d}dy\Big)^{\frac{2}{d}}+\int_{B_{\frac{\ell}{2}}(x)} |\nabla u_\ell|^2dy,
$$
and therefore by definition of $u_\ell$ and $f_\ell$ and Cauchy-Schwarz' inequality,
\begin{equation}\label{eq:estim-u_0}
\mu^2\bar u_0^2\,=\,\mu^2 \Big(\fint_{B_{\frac{3\ell}{2}}(x)} u_0dy\Big)^2
\,\lesssim \,\Big(\int_{B_{2\ell}(x)} |f|dy\Big)^2+\int_{B_{2\ell}(x)} |\nabla u_0|^2dy.
\end{equation}
The combination of \eqref{eq:triang-u-tilde-0}, \eqref{eq:u-baru0}, \eqref{eq:tildeu-baru0}, 
\eqref{eq:energy_nablau-tilde-0}, \eqref{bar_diff}, and \eqref{eq:estim-u_0} then yields
\begin{equation*}
 \| u - \tilde u \|_{L^{\lambda_1'}(B_\ell(x))} \lesssim (\nabla u)_{9\ell}(z)+ \| f \|_{L^{\lambda_2}(B_{2\ell}(x))},
\end{equation*}
which proves \eqref{eq:osc-estim-u} for $|x-z|\le 6\ell$.



\section*{Acknowledgments}

The first author acknowledges financial support from the European Research
Council under the European Community Seventh Framework Programme
(FP7/2014-2019 Grant Agreement QUANTHOM 335410).

\bibliographystyle{plain}

\end{document}